\documentclass[11pt]{article}

\usepackage[
  left=1in,
  right=1in,
  top=1in,
  bottom=1in
]{geometry}
\usepackage[utf8]{inputenc}
\usepackage[T1]{fontenc}

\usepackage{amsmath, amsfonts, amssymb, amsthm}
\usepackage{mathtools}
\usepackage{bm}
\usepackage{bbm}
\usepackage{algorithm}
\usepackage{algpseudocode}
\usepackage{tikz}
\usepackage{booktabs}
\usepackage{multirow}
\usepackage{makecell}
\usepackage{tabularray}
\usepackage{adjustbox}
\usepackage{subcaption}
\usepackage{rotating}

\usepackage[title]{appendix}
\usepackage[subnum]{cases}
\usepackage{accents}
\usepackage[normalem]{ulem}
\usepackage{cancel}
\usepackage{inconsolata}
\usepackage{tikz}
\usetikzlibrary{shapes.geometric, shapes.symbols, positioning, arrows.meta, calc}

\usepackage[colorlinks=true, allcolors=blue]{hyperref}
\usepackage{natbib}
\bibpunct[, ]{(}{)}{,}{a}{}{,}%
\usepackage[title]{appendix}
\newtheorem{theorem}{Theorem}
\newtheorem{lemma}{Lemma}
\newtheorem{proposition}{Proposition}
\newtheorem{corollary}{Corollary}
\newtheorem{definition}{Definition}
\newtheorem{example}{Example}
\newtheorem{assumption}{Assumption}
\newtheorem{notation}{Notation}
\newtheorem{remark}{Remark}

\newcommand{\abs}[1]{\left\lvert #1 \right\rvert}

\newcommand{\inprod}[1]{\left\langle #1 \right\rangle }
\newcommand{\norm}[1]{ \left\|#1\right\|}

\newcommand{\doi}[1]{}
\newcommand{\Halmos}{\hfill\ensuremath{\square}} 
 
\title{\fontsize{14}{12}\selectfont\textbf{Concave Certificates: Geometric Framework \\ for Distributionally Robust Risk and Complexity Analysis}}

\author{
    Hong T.M. Chu \\
    \small College of Engineering and Computer Science, VinUniversity \\
    \small {hong.ctm@vinuni.edu.vn}
}

\date{\today}

\begin{document}

\maketitle

\begin{abstract}

Distributionally Robust (DR) optimization aims to certify worst-case risk within a Wasserstein uncertainty set. Current certifications typically rely either on global Lipschitz bounds, which are often conservative, or on local gradient information, which provides only a first-order approximation. This paper introduces a novel geometric framework based on the least concave majorants of the growth rate functions. Our proposed concave certificate establishes a tight bound on DR risk that remains applicable to non-Lipschitz and non-differentiable losses. We extend this framework to complexity analysis, introducing the worst-case generalization bound that complements the standard statistical generalization bound. Furthermore, we utilize this certificate to bound the gap between adversarial and empirical Rademacher complexity, demonstrating that dependencies on input diameter, network width, and depth can be eliminated. For practical application in deep learning, we introduce the adversarial score as a tractable relaxation of the concave certificate that enables efficient and layer-wise analysis of neural networks. We validate our theoretical results in various numerical experiments on classification and regression tasks using real-world data.

\end{abstract}

\vspace{1em}

\section{Introduction}\label{sec:intro}

Given feature data \(X\) and label \(Y\), we seek a parameterized network \(f_{\theta} \) to model their relationship \(Y\approx f_{\theta}(X) \). This is typically formulated as minimizing the expected loss 
\begin{equation}\label{eq:ERM}
    \inf_{\theta\in\Theta} \mathbb{E}_{Z\sim\mathbb{P}_{\rm true}}[\bm{l}(Z;\theta)],
\end{equation}
where \(Z=(X,Y)\), \(\Theta\) is the set of feasible parameters and \(\bm{l} \) is the loss function. The true data distribution \(\mathbb{P}_{\rm true}\) in \eqref{eq:ERM} is often unknown and approximated by the empirical distribution  \(\mathbb{P}_N = \frac{1}{N}\sum_{i=1}^{N} \bm{\chi}_{\{Z^{(i)}\}} \), which can lead to over-fitting. To mitigate this issue, the robust counterpart of \eqref{eq:ERM} aims to minimize the worst-case loss within a neighborhood of \(\mathbb{P}_N\) by solving
\begin{equation*}\label{eq:inf-sup}
    \inf_{\theta\in\Theta}\sup_{\mathbb{P} : \mathcal{D}(\mathbb{P},\mathbb{P}_N) \leq \epsilon }  \ \mathbb{E}_{\mathbb{P}} [\bm{l} (Z;\theta) ],
\end{equation*}
where \(\mathcal{D}\) is a discrepancy on the probability space \(\mathcal{P}(\mathcal{Z})\). In this work, we focus on the Wasserstein discrepancy (Definition~\ref{def:Wass}), which intuitively represents the minimum cost to transport the mass of \(\mathbb{P}\) to that of \(\mathbb{P}_N\). The inner supremum problem is referred to as the distributionally robust (DR) risk:
\begin{equation}\label{eq:DR-loss}
    \text{(DR risk)}\quad  \mathcal{R}_{p}(\epsilon) = \sup_{\mathbb{P} \colon \mathcal{W}_{p}(\mathbb{P},\mathbb{P}_N) \leq \epsilon }   \mathbb{E}_{\mathbb{P}} [\bm{l} (Z;\theta) ].
\end{equation}

\textbf{Quantify Robustness.} Essentially, the DR risk \eqref{eq:DR-loss} quantifies the sensitivity of the loss value under distributional shifts bounded by a budget $\epsilon$. In general, computing \(\mathcal{R}_{p}(\epsilon)\) is intractable. To bypass this, two primary approaches are used: the {\textit{Lipschitz certificate}} and the {\textit{gradient certificate}}. The Lipschitz certificate estimates an upper bound of  \(\mathcal{R}_{p}(\epsilon)\). For instance, if \(\bm{l}\) is \(L_{\theta}\)-Lipschitz then \(\mathcal{R}_{p}(\epsilon)\leq \mathbb{E}_{\mathbb{P}_N} [\bm{l} (Z;\theta) ] + L_{\theta}\epsilon \), and this bound is known to be tight for linear hypothesis \citep{goh2010distributionally,blanchet2019quantifying,blanchet2019robust,an2021generalization,gao2022wasserstein,gao2023distributionally}. We refer reader to a recent survey  \cite{zuhlke2025adversarial} on how Lipschitz calculus can be used to study robustness. Estimating the global Lipschitz constant \(L_{\theta}\) for deep networks often reduces to a layer-wise estimation of Lipschitz constants \citep{virmaux2018lipschitz,shafieezadeh2019regularization,Latorre2020Lipschitz}. This raises some fundamental questions: why do functions like the entropy loss or the square-root loss \citep{belloni2011square} exhibit robustness despite being non-Lipschitz?  Furthermore, even though the \(4\times\)Sigmoid, Tanh and ReLU activations share a Lipschitz modulus of 1, why  do they not possess identical robustness properties?  Do modern architectures with  LayerNorm or Attention exhibit these common  robustness properties as well? Alternatively, the other approach of gradient certificate \citep{bartl2021sensitivity,gao2022finite,bai2024wasserstein} approximates \(\mathcal{R}_{p}(\epsilon)\)  using first-order information as  \(\mathcal{R}_{p}(\epsilon) \approx \mathbb{E}_{\mathbb{P}_N} [\bm{l} (Z;\theta) ] + \operatorname{grad}_{*}\epsilon  \) where \(1/p+1/q=1 \) and \(\operatorname{grad}_{*}= \left( \mathbb{E}_{\mathbb{P}_N} [ \norm{\nabla_x \bm{l} (Z;\theta)}^q ] \right) ^{1/q}\).  However, this first-order estimation is asymptotic and holds only as the budget \(\epsilon\rightarrow0\). Moreover, it requires  \(\bm{l}\) to be differentiable and  does not provide a true upper bound of \(\mathcal{R}_{p}(\epsilon)\). Tighter bounds on $\mathcal{R}_{p}$ are also emerging from a separate line of work \citep{pal2023adversarial,pal2024certified} on $(\epsilon,\delta)$-robust classifiers, where it is crucial to exactly characterize when \(\mathcal{R}_p(\epsilon)\leq\delta\).

\textbf{Generalization Capability.} To understand the model's generalization capability in this robust setting, we recall the notion of Rademacher complexity \citep{bartlett2002model,koltchinskii2002empirical}. 
For a class of loss functions \(\mathcal{L} \coloneqq \{z \mapsto \bm{l}(z; \theta) \mid \theta \in \Theta\}\), the empirical Rademacher Complexity (RC) measures the richness of $\mathcal{L}$ by its ability to correlate with random noise on the  sample $\mathcal{Z}_N $:
\begin{equation} \label{eq:RC} 
    \text{(RC)}\quad  \hat{\mathfrak{R}}_{\mathcal{Z}_N}(\mathcal{L}) = \mathbb{E}_{\sigma} \left[ \sup_{\theta \in \Theta} \frac{1}{N} \sum_{i=1}^N \sigma_i \bm{l}(Z^{(i)}; \theta) \right],
\end{equation}
where $\sigma = (\sigma_1, \dots, \sigma_N)$ are independent Rademacher variables taking values in $\{-1, +1\}$ with equal probability. Intuitively, if RC is large, then \(\mathcal{L}\) has the capacity to fit arbitrary noise \(\sigma\), leading to overfitting.  Conversely, a small RC indicates that $\mathcal{L}$ learns meaningful patterns. Standard results by \cite{bartlett2002rademacher,koltchinskii2002empirical} show that \(\hat{\mathfrak{R}}_{\mathcal{Z}_N}(\mathcal{L}) \) directly bounds the generalization gap:
\begin{equation} \label{eq:stat-gen-bound}
    \mathbb{E}_{\mathbb{P}_{\rm true}}[\bm{l}(Z;\theta)] \lesssim \mathbb{E}_{\mathbb{P}_{N}}[\bm{l}(Z;\theta)] + \operatorname{const}\times\hat{\mathfrak{R}}_{\mathcal{Z}_N}(\mathcal{L}) + \operatorname{conf}(\delta),
\end{equation}
where \(\operatorname{conf}(\delta)\) is a confidence term. Therefore, a small RC also implies a tight generalization bound. In the context of DRO, we focus on the class of worst-case loss functions $\tilde{\mathcal{L}}_{\epsilon} \coloneqq \{z \mapsto \tilde{\bm{l}}_{\epsilon}(z;\theta)= \sup_{z' : d(z', z) \le \epsilon} \bm{l}(z'; \theta) \mid \theta \in \Theta\}$, which leads to the definition of Adversarial Rademacher Complexity (ARC) \citep{khim2018adversarial,yin2019rademacher,awasthi2020adversarial,xiao2022adversarial}:
\begin{equation} \label{eq:ARC}
     \text{(ARC)}\, \hat{\mathfrak{R}}_{\mathcal{Z}_N}(\tilde{\mathcal{L}}_{\epsilon}) = \mathbb{E}_{\sigma} \left[ \sup_{\theta \in \Theta} \frac{1}{N} \sum_{i=1}^N \sigma_i   \tilde{\bm{l}}_{\epsilon}(Z^{(i)};\theta) \right].
\end{equation}
Deriving tight bounds for ARC is significantly more challenging than for RC because the inner supremum operator destroys  structural properties that are typically exploited in traditional complexity analysis. For linear hypotheses, \citet{khim2018adversarial, yin2019rademacher} establish that the gap between ARC and RC scales linearly with the weight norm \(\norm{\theta}\), which serves as the Lipschitz constant. For deep neural networks, however, \citet{awasthi2020adversarial} and \citet{xiao2022adversarial} show these bounds grow not just with weights, but also with the network’s depth, width, and data diameter $\mathcal{Z}_N$. This predicted surge is counter-intuitive when compared to the DR risk \(\mathcal{R}_p\) mentioned earlier that the adversarial risk of a feedforward network is controlled by its Lipschitz constant, suggesting that the actual complexity should not blow up simply because a network becomes larger.

\begin{figure}[htbp]
    \centering
    \resizebox{0.9\textwidth}{!}{%
        \begin{tikzpicture}[
            % -------------------------------------------------------------------------
            % GLOBAL SETTINGS
            % -------------------------------------------------------------------------
            font=\tiny,
            % -------------------------------------------------------------------------
            % SHAPE STYLES
            % -------------------------------------------------------------------------
            thm/.style={rectangle, draw=gray!80!black, thick, fill=gray!10, align=center, rounded corners=3pt, minimum height=1.4cm, inner xsep=0.2cm},
            prop/.style={ellipse, draw=blue!80!black, thick, fill=blue!10, align=center, minimum height=1.4cm, inner xsep=0.1cm},
            % The 'signal' shape perfectly creates a flat top/bottom with pointed left/right sides
            exam/.style={signal, signal to=east and west, draw=orange!80!black, thick, fill=orange!10, align=center, minimum height=1.4cm, inner xsep=0.15cm},
            def/.style={rectangle, draw=gray!80!black, dashed, thick, fill=gray!10, align=center, rounded corners=3pt, minimum height=1.4cm, inner xsep=0.2cm},
            arrow/.style={->, >=Stealth, thick, draw=black!80},
            darrow/.style={->, >=Stealth, thick, dashed, draw=black!80}
        ]

            % -------------------------------------------------------------------------
            % THE ROOT
            % -------------------------------------------------------------------------
            \node[thm] (thm1) at (0, 0) {\textbf{Theorem~\ref{thm:main0}}\\(certify DR risk $\mathcal{R}_p$ \eqref{eq:DR-loss} via geometric majorants)};

            % -------------------------------------------------------------------------
            % LEFT COLUMN: DR Risk Properties & Classifiers
            % -------------------------------------------------------------------------
            \node[prop] (cor1)  at (-3.6, -2.0) {Corollary~\ref{cor:p-dynamic}\\($p$-dynamic)};
            \node[prop] (cor2)  at (-6.2, -2.0) {Corollary~\ref{cor:finite-Rp}\\(finiteness)};
            \node[exam] (ex1)   at (-4.9, -4.0) {Example~\ref{exam:unbounded}\\(tightness)};

            \node[prop] (cor3)  at (-4.9, -6.0) {Corollary~\ref{cor:Lip-ub}\\(Lipschitz)};
            \node[thm] (prop1) at (-4.9, -8.0) {\textbf{Proposition~\ref{prop:ExactCD}}\\(Exact CD)};

            % -------------------------------------------------------------------------
            % CENTER COLUMN (SWAPPED IN): Neural Network Adversarial Scores
            % -------------------------------------------------------------------------
            \node[thm]  (def6)  at (0, -2.0) {\textbf{Definition}~\ref{def:adv-score}\\(adversarial score for DNN)};

            \node[prop] (lem5)  at (0, -4.0) {Lemma~\ref{lem:layer-rule}\\(layer rule)};

            \node[thm] (prop4) at (-1.2, -6.0) {\textbf{Proposition}~\ref{prop:classification}\\(classification)};
            \node[thm] (prop5) at ( 1.2, -6.0) {\textbf{Proposition}~\ref{prop:regression}\\(regression)};

            \node[exam] (ex6)   at ( 0, -8.0)  {Example~\ref{exam:activation}+\ref{exam:gamma}\\(DNN)};

            % -------------------------------------------------------------------------
            % RIGHT COLUMN (SWAPPED IN): Generalization & Complexity Gaps
            % -------------------------------------------------------------------------
            \node[prop] (cor4)  at (5.0, -2.0)  {Corollary~\ref{cor:det-gen-bound}\\(gen. bound via \(\hat{\mathfrak{C}}_N\))};
            \node[thm] (prop2) at (3.2, -4.0)  {\textbf{Proposition}~\ref{prop:calculus}+\ref{prop:contraction}\\(\(\hat{\mathfrak{C}}_N\) v.s. \(\hat{\mathfrak{R}}_N\))};

            \node[exam] (ex2)   at (7.2, -4.0)   {Example~\ref{exam:concave-complexity}\\(\(\hat{\mathfrak{C}}_N\) of linear\&DNN)};

            \node[thm]  (thm2)  at (5.0, -6.0) {\textbf{Theorem~\ref{thm:ACG}}\\(adversarial complexity gap ACG)};

            \node[exam] (ex3)   at (5.0, -8.0) {Example~\ref{exam:ACG-linear}+\ref{exam:ACG-MLP}\\(ACG of linear\&MLP )};

            % -------------------------------------------------------------------------
            % ROUTING: ARROWS & LOGIC FLOW
            % -------------------------------------------------------------------------

            % Thm 1 Outbound
            \draw[arrow] (thm1) -- (def6);
            \draw[arrow, rounded corners] (thm1.west) -| (cor1.north);
            \draw[arrow, rounded corners] (thm1.east) -| (cor4.north);
            \draw[arrow, rounded corners] (thm1.west) -| (cor2.north);

            % Long vertical bus on the far left for Cor3 and Prop 1
            \draw[arrow, rounded corners] (thm1.west) -- (-7.6, 0) |- (cor3.west);
            \draw[arrow, rounded corners] (thm1.west) -- (-7.6, 0) |- (prop1.west);

            % Left Column Flow
            \draw[arrow] (cor1) -- (ex1);
            \draw[arrow] (cor2) -- (ex1);

            % Center Column Flow (MLP)
            \draw[arrow] (def6) -- (lem5);
            \draw[arrow] (lem5) -- (prop4);
            \draw[arrow] (lem5) -- (prop5);
            \draw[arrow] (prop4) -- (ex6);
            \draw[arrow] (prop5) -- (ex6);

            % Right Column Flow (Gen/ACG)
            \draw[arrow] (cor4) -- (prop2);
            \draw[arrow] (prop2) -- (ex2);
            \draw[arrow] (prop2.south) --  (thm2.north);
            \draw[arrow] (thm2) -- (ex3);

        \end{tikzpicture}%
    }
    \caption{Logical structure and reading flow of the paper. Building upon the foundational bounds established in Theorem~\ref{thm:main0}, the first branch connects our results to duality, robust certificates, and the theoretical analysis of $\mathcal{R}_p$. The second branch extends our framework to deep neural networks (DNN) via tractable adversarial scores. Finally, the third branch derives deterministic generalization bounds and adversarial complexity gaps (ACG).  }
    \label{fig:summary}
\end{figure}

\textbf{Main Contributions.} We summarize our main contributions and organize our paper as follows. (See Figure~\ref{fig:summary} for a diagram of summary.)
\begin{itemize}
    \item We introduce a novel framework  to estimate the distributionally robust risk \( \mathcal{R}_{p}(\epsilon)\). This geometric framework establishes an elegant bound showing that the robust-empirical risk gap \( \mathcal{R}_{p}(\epsilon) -\hat{\mathcal{R}} \) stays between the average of the \textit{least star-shaped majorants} \(\operatorname{lb}_p(\epsilon) \)  and the \textit{least concave majorant} \(\operatorname{cc}_p(\epsilon) \)  of the loss's growth functions  (Theorem~\ref{thm:main0}). Notably, we do not require the loss  \(\bm{l}\) to be convex/differentiable/Lipschitz, or the cost  \(d\) to be a metric, or the domain  \(\mathcal{Z}\) to be bounded. Our analysis reveals how the DR risk \( \mathcal{R}_{p}(\epsilon)\) evolves as the exponent \(p\) changes (Corollary~\ref{cor:p-dynamic}) and provides exact conditions for determining whether \( \mathcal{R}_{p}(\epsilon)\) is finite (Corollary~\ref{cor:finite-Rp}). In addition, our framework leads to a necessary and sufficient condition for the existence of robust classifiers (Proposition~\ref{prop:ExactCD}).

    \item To facilitate practical implementation for Euclidean domains, we introduce the adversarial score $\mathcal{A}_{\theta}$. This relaxation enables layer-wise  calculation rules via a composition and product maps, providing explicit certificates for classification and regression in deep learning. It successfully accounts for modern architectures (e.g., LayerNorm, Attention) and yields tighter bounds for non-Lipschitz/non-differentiable losses  (Figure~\ref{fig:adv-score}).
    
    \item We introduce the Concave Complexity (CC) \(\hat{\mathfrak{C}}_{\mathcal{Z}_N}(\mathcal{L},\epsilon)\) \eqref{eq:CC}  that complements the standard Rademacher Complexity (RC). Besides deriving a worst-case version of the generalization bound,  CC reveals that the adversarial-empirical complexity gap is strictly controlled by the complexity of the rate class. Although CC trades the standard $\mathcal{O}(1/\sqrt{N})$ statistical rate of RC for the optimal transport budget $\epsilon \approx \mathcal{O}(N^{-1/n})$, it remains beneficial for analyzing certain over-parameterized networks.

    \item We validate our theoretical results through two experiments using real-world datasets. In the regression task (Section~\ref{sec:madrid}), we numerically demonstrate that our adversarial score is strictly tighter and more informative than traditional Lipschitz and gradient-based certificates. In the classification task  (Section~\ref{sec:mnist}), we verify that the adversarial-empirical Rademacher gap does not blow up with the depth, width or data dimension. We conclude our paper in Section~\ref{sec:con}.

\end{itemize}

\section{Preliminaries and Notations} \label{sec:pre} 

Let the indicator function \(\bm{\delta}_{S}  \colon  \mathcal{Z}  \rightarrow \mathbb{R}\) of a set $S\subset\mathcal{Z} $ be defined as \(\bm{\delta}_{S}(z) =  0\) if \(z\in S\), and \(\infty\) otherwise.  Let the point mass function (Dirac measure) $\bm{\chi}_{\{\hat{z}\}}\in\mathcal{P}(\mathcal{Z})\colon \bm{A} \rightarrow \mathbb{R}  $ at point $\hat{z}\in \mathcal{Z}$ be defined as   \(\bm{\chi}_{\{\hat{z}\} }(A) = 1\) if \(\hat{z} \in A\), and \(0\) otherwise. We adopt the convention of extended arithmetic such that \(0\cdot\infty=0 \). The Rademacher random variable is \(\sigma = \pm 1 \) where \(P(\sigma=-1)=P(\sigma=1)=\frac{1}{2} \). For any real number $t$, the sign function is defined as $\operatorname{sgn}(t)=-1$ if $t<0$, and $\operatorname{sgn}(t)=1$ otherwise. For any positive integer \(n\), we denote \([n]\coloneqq \{1,2,\dots,n\}\). Denote the inner product on \(\mathbb{R}^n\) by \(\inprod{x,y} = \sum_{i=1}^{n}x_i y_i  \) for any \(x,y\in\mathbb{R}^n\). Let \(\norm{\cdot}\) be an arbitrary norm on \(\mathbb{R}^{n}\) and \(\norm{\cdot}_{*} \) be its dual norm defined as \(\norm{x}_{*}\coloneqq  \max_{y\in\mathbb{R}^n} \left\{ \inprod{x,y} \mid \norm{y}_{\mathbb{R}^n} =1 \right\}  \).

A set \(\emptyset\ne\Omega\subset\mathbb{R}^n\) is convex if \(\eta x +(1-\eta)x'\in\Omega \) for any \(x,x'\in\Omega \) and \(\eta\in[0,1]\). A function  \(f\colon \Omega \rightarrow \mathbb{R}\) is concave if its hypograph \( \operatorname{hypo} f= \left\{ (x,y) \in \Omega\times\mathbb{R} \mid y \leq f(x) \right\} \) is convex. A set \(\emptyset\ne\Omega\subset\mathbb{R}^n\) is star-shaped (with respect to origin \(\bm{0}_n\)) if \(\eta x \in\Omega \) for any \(x\in\Omega \) and \(\eta\in[0,1]\).  A function  \(f\colon \Omega \rightarrow \mathbb{R}\) where \(\bm{0}_n\in\Omega \) and \(f(0)\geq 0\) is star-shaped if its hypograph \( \operatorname{hypo} f\) is star-shaped. (Note that this notion mirrors \citet[16.B.9]{marshall1979inequalities} in which \(f\) is star-shaped if \(f(0)\leq 0\) and its epigraph is star-shaped.)  Obviously, a concave function is star-shaped.

In this work, we are interested in the smallest concave/star-shaped upper bound of a non-negative univariate function on \([0,\infty)\). These concepts have been used  to analyze the magnitude of  Brownian motion or behaviors of regressors  \citep{pitman1982remarks,groeneboom1983concave,bennett1988interpolation}.

\begin{definition}[least concave and star-shaped majorants] \label{def:least-majorant}
Given \(f\colon [0,\infty) \rightarrow [0,\infty)  \), define the least concave majorant \(\mathcal{C}_{f} \colon [0,\infty)  \rightarrow \mathbb{R}\cup\{+\infty\}\) and the least star-shaped majorant \(\mathcal{S}_{f} \colon [0,\infty)  \rightarrow \mathbb{R}\cup\{+\infty\}\) of \(f\) as follows.
\begin{equation*}
\begin{array}{rl}
     \mathcal{C}_{f} (t) &\coloneqq \operatorname{inf}\left\{ H(t) \mid H(t) \geq f(t),\, H \text{ is concave} \right\},\\
     \mathcal{S}_{f} (t) &\coloneqq \operatorname{inf}\left\{ H(t)  \mid  H(t)  \geq  f(t),\, H \text{ is star-shaped} \right\} .
\end{array}   
\end{equation*}
\end{definition}

It is worth noting that the definitions of least majorant are valid, since the infimum of a collection of functions is equivalent to the intersection of their hypographs; thus, concavity or star-shapedness is induced immediately. The following lemma is derived directly from the definitions and \citet{rockafellar-convex-analysis,marshall1979inequalities,hardy1988inequalities,groeneboom2014nonparametric}.

\begin{lemma} \label{lem:calculation}
Suppose that \(f\colon [0,\infty) \rightarrow [0,\infty)  \).  
\begin{itemize}
    \item \(\mathcal{S}_f(t) \leq \mathcal{C}_f(t) \leq \inf_{a,b\in\mathbb{R}}\{at+b \mid au+b \geq f(u) \forall u \geq 0 \} \); 
      \(\mathcal{S}_f(t) = \sup_{u\in[t,\infty)} \frac{t f(u)}{u}\) for any \(t>0 \).
    \item If \(f_1\leq f_2\) then \(\mathcal{S}_{f_1}\leq \mathcal{S}_{f_2}\) and \(\mathcal{C}_{f_1}\leq \mathcal{C}_{f_2}\).
    \item If \(f\) is non-decreasing then \(\mathcal{S}_f\) and \(\mathcal{C}_f\) are non-decreasing as well.
    \item Since  \(\mathcal{C}_{f_{\theta}} \leq \mathcal{C}_{\sup_{\theta\in\Theta}f_{\theta}} \) for any \(\theta\in\Theta\), thus \(\sup_{\theta\in\Theta}\mathcal{C}_{f_{\theta}} \leq \mathcal{C}_{\sup_{\theta\in\Theta}f_{\theta}} \) and \(\mathcal{C}_{f_1+f_2}\leq \mathcal{C}_{f_1}+\mathcal{C}_{f_2}\).
\end{itemize}
\end{lemma}

Finally, we recall definition of the Wasserstein discrepancy, which serves as a metric to measure the difference between two probability distributions.

\begin{definition}[Wasserstein discrepancy]\label{def:Wass}
    Given two probability distributions \(\mathbb{P},\mathbb{Q}\in\mathcal{P}(\mathcal{Z}) \) and a non-negative function \(d\colon \mathcal{Z}\times\mathcal{Z}\rightarrow [0,\infty] \), the Wasserstein discrepancy with respect to \(d\) and an exponent  \(p\in[1,\infty]\) is defined via the Kantorovich problem \citep{villani2009optimal,peyre2017computational} as follows.
    \begin{itemize}
        \item If \(p\in[1,\infty)\), then \(\mathcal{W}_{p}(\mathbb{P},\mathbb{Q}) \triangleq \left(\inf_{\pi \in \Pi(\mathbb{P},\mathbb{Q})} \int_{\mathcal{Z}\times\mathcal{Z}} d^{p}(z',z)\mathrm{d}\pi(z',z)\right)^{1/p}\).
\item If \(p=\infty\), then \(\mathcal{W}_{\infty}(\mathbb{P},\mathbb{Q}) \triangleq \inf_{\pi \in \Pi(\mathbb{P},\mathbb{Q})}\operatorname{ess.sup}_{{\pi}}(d).\)
    \end{itemize}
Here \(\Pi(\mathbb{P},\mathbb{Q})\) \cite[Definition 1.1]{villani2009optimal} denotes the set of all couplings  (joint probability distributions) between \(\mathbb{P}  \) and \(\mathbb{Q}\), i.e., the set of all \(\pi \in\mathcal{P}(\mathcal{Z}\times\mathcal{Z}) \) such that \(\pi (A\times \mathcal{Z}) = \mathbb{P}(A)\) and \(\pi(\mathcal{Z}\times B) = \mathbb{Q}(B) \) for all measurable sets \(A,B\subset\mathcal{Z}\). 
\end{definition}

The following notation is adopted throughout this paper.
\begin{notation}\label{nota:main}
    Let \(\mathcal{Z}\) be a measurable space,  \(d\colon\mathcal{Z}\times\mathcal{Z}\rightarrow[0,\infty]\) be a cost function  on \(\mathcal{Z}\) (that is, \(d\) is measurable and \(d(z,z)=0\) for any \(z\in\mathcal{Z}\)),  and  \(\bm{l}:\mathcal{Z}\times \Theta \rightarrow \mathbb{R} \) be a  loss function.  Let \(\mathcal{Z}_N\coloneqq\{Z^{(1)},\dots,Z^{(N)}\}\subset\mathcal{Z} \) be a (finite) empirical dataset and \(\mathbb{P}_N\coloneqq\sum_{i=1}^{N}\mu_i\bm{\chi}_{\{Z^{(i)}\}}\in\mathcal{P}(\mathcal{Z})\) be the corresponding empirical distribution. Denote the empirical loss as \(\hat{\mathcal{R}}\coloneqq \mathbb{E}_{\mathbb{P}_N}[\bm{l}(Z;\theta)] \). Given a parameter \(\theta\in\Theta\), a positive budget \(\epsilon>0\)  and an extended-value number \(p\in[1,\infty]\), we define distributionally robust risk (DR risk), Rademacher complexity (RC) and adversarial Rademacher complexity (ARC) as in \eqref{eq:DR-loss}, \eqref{eq:RC} and \eqref{eq:ARC}, respectively.
\end{notation}
\begin{lemma}[\(\mathcal{W}_p\) order]\label{lem:Wp-ineq} (See Remark 6.6, \cite{villani2009optimal})  Given Notation~\ref{nota:main}, if \(1\leq p\leq p_2\leq\infty\) then \(\mathcal{W}_{p}(\mathbb{P},\mathbb{P}_N) \leq \mathcal{W}_{p_2}(\mathbb{P},\mathbb{P}_N) \) and  \(\mathcal{R}_{p}(\epsilon)\geq \mathcal{R}_{p_2}(\epsilon) \). 
\end{lemma}

\begin{lemma}
    [Point-wise RO is \(\mathcal{W}_{p=\infty}\)DRO] (See Appendix~\ref{proof:p-infinity}) \label{lem:p-infinity} Given Notation~\ref{nota:main}, then
    \begin{equation*}
        \begin{array}{cc}
            \sum_{i=1}^{N}\mu_i \sup_{\tilde{z}\in B_{d,\epsilon}^{(i)}} \bm{l}(\tilde{z};\theta) \leq \mathcal{R}_{\infty}(\epsilon) \leq \sum_{i=1}^{N}\mu_i \sup_{\tilde{z}\in B_{d,\epsilon+\rho}^{(i)}} \bm{l}(\tilde{z};\theta),
        \end{array}
    \end{equation*}
    for any \(\rho>0\), where \(B_{d,\epsilon}^{(i)} = \{ \tilde{z} \colon d(\tilde{z},Z^{(i)})\leq \epsilon \} \).
\end{lemma}

\section{Theoretical Analysis}\label{sec:theory}

We begin by formally defining key quantities that measure how the loss changes in response to localized perturbations of the data.  In fact, this is a generalization of the (scalar) growth rate notion proposed in \cite{gao2023distributionally} and connects directly to the adversarial loss \(\tilde{\bm{l}}(\hat{z};\theta)  \) \eqref{eq:ARC} \citep{khim2018adversarial,yin2019rademacher,xiao2022adversarial}. 

\begin{figure}[h]
    \centering
    \includegraphics[width=0.9\linewidth]{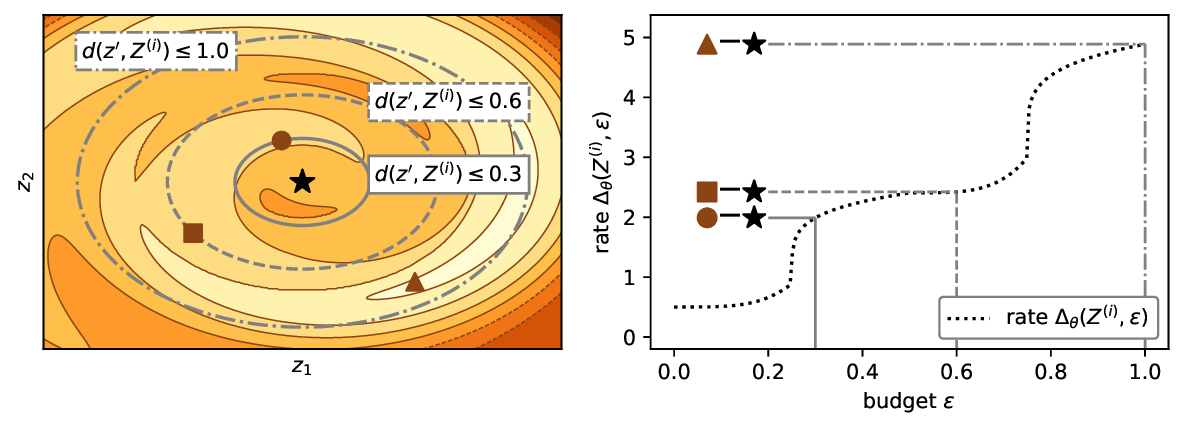}
    \caption{Illustration of individual rate (Definition~\ref{def:rate}). Given an empirical point $Z^{(i)}$ at $\star$, three brown points ($\color{brown}\bullet,\,\blacksquare$ and $\color{brown}\blacktriangle$) are the maximizers (lightest contours) of the loss within the set $\{z' \colon d(z',Z^{(i)})\leq \epsilon \} $ where \(\epsilon=0.3,0.6\) and \(1.0\). The individual rate $\Delta_{\theta}(Z^{(i)},\epsilon)$ is defined as the  difference in loss between each $\color{brown}\bullet,\,\blacksquare,\,\blacktriangle$ and $\star$. Note that we do not require \(d\) to be a metric.}
    \label{fig:rate}
\end{figure}

\begin{definition}[Rate Function - Figure~\ref{fig:rate}]\label{def:rate} Given Notation~\ref{nota:main}, we define the \textbf{individual} rate \(\Delta_{\theta}\) of the loss \(\bm{l}\) at \(z\) as 
 \begin{equation}\label{eq:ind-rate}
    \Delta_{\theta}(z,t) \triangleq \sup_{z'\in\mathcal{Z}} \left\{ \bm{l}(z';\theta) - \bm{l}(z;\theta) \mid d(z',z) \leq  t \right\},
\end{equation}
for any \(z\in\mathcal{Z}_N\) and \(t\geq0\). We define the (empirical) \textbf{maximal} rate as
\begin{equation}\label{eq:max-rate}
        \Delta_{\theta}^{\max}(t) = \max_{Z^{(i)}\in\mathcal{Z}_N}\Delta_{\theta}(Z^{(i)},t).
    \end{equation}
\end{definition}

 By shifting the focus from the complex loss function \(\bm{l}\) to these univariate rate functions, our first theoretical result proposes to bound the adversarial-empirical risk gap through the geometric construction of least concave and star-shaped majorants.

\begin{theorem}[Distributional Robustness Certificates via Least Majorants]\label{thm:main0}
    Given Notation~\ref{nota:main}, define  the least concave majorant \(\mathcal{C}_f \), the least star-shaped majorant \(\mathcal{S}_f\), the individual rate \(\Delta_{\theta}(Z^{(i)},t) \), and maximal rate \(\Delta_{\theta}^{\max}(t) \) as in Definitions~\ref{def:least-majorant} and \ref{def:rate}. Then for any \(\epsilon>0\),
	\begin{equation} \label{eq:main-L}
		\begin{array}{cc}
		     \mathcal{R}_{p}(\epsilon) - \hat{\mathcal{R}}  \geq \operatorname{lb}_{p}(\epsilon) =  \sum_{i=1}^{N}\mu_i   s^{(i)}(\epsilon), 
		\end{array}
	\end{equation}
	where  \( s^{(i)}(\epsilon) = \mathcal{S}_{f^{(i)}} (\epsilon^p) \) with \(f^{(i)}\colon t \mapsto \Delta_{\theta}(Z^{(i)},t^{1/p} )\) if \(p<\infty \), and \( s^{(i)}(\epsilon) = \Delta_{\theta}(Z^{(i)},\epsilon)  \) if \(p=\infty\). In addition,
    \begin{equation} \label{eq:main-U}
        \mathcal{R}_{p}(\epsilon)-  \hat{\mathcal{R}} \leq \operatorname{cc}_{p}(\epsilon),
    \end{equation}
    where \(\operatorname{cc}_{p}(\epsilon) = \mathcal{C}_{f^{\max}}(\epsilon^p) \) with \(f^{\max}\colon t \mapsto \Delta_{\theta}^{\max}(t^{1/p} )\) if \(p<\infty \), and  \(\operatorname{cc}_{p}(\epsilon) = \sum_{i=1}^{N}\mu_i \lim\limits_{t\rightarrow\epsilon^{+}}\Delta_{\theta}(Z^{(i)},t)\) if \(p=\infty\).
	
\end{theorem}

\noindent\textit{Sketch of Proof.} The cases of \(\operatorname{lb}_{\infty}\) and \(\operatorname{cc}_{\infty}\) are immediate results by Lemma~\ref{lem:p-infinity}. We now consider \(p\in[1,\infty) \). For \(p<\infty\), we first derive the lower bound \(\operatorname{lb}_{p}(\epsilon)\) by constructing a discrete perturbation \(\tilde{\mathbb{P}} \) such that \(\mathcal{W}_{p}(\tilde{\mathbb{P}},\mathbb{P}_N)\leq \epsilon \) and show that the risk  gap  \(\mathbb{E}_{\tilde{\mathbb{P}}}[\bm{l}(Z;\theta)] - \mathbb{E}_{{\mathbb{P}_N}}[\bm{l}(Z;\theta)] \geq \sum_i \mu_i\sup_{u\in[\epsilon^p,\infty)}  
\left\{ \frac{\epsilon^{p} \Delta_{\theta}(Z^{(i)},u^{1/p})}{u}  \right\}  \), which is equal to \(\sum_{i=1}^{N} \mu_i \mathcal{S}_{f^{(i)}}(\epsilon^p)\) by Lemma~\ref{lem:calculation}. We then derive the upper bound by rewriting the risk gap as \(\int_{\mathcal{Z}\times\mathcal{Z}_N} \left(\bm{l}(\tilde{z};\theta)-\bm{l}({z};\theta)\right)\mathrm{d}\tilde{\pi}(\tilde{z},z)\) where \(\tilde{\pi}\) is the optimal coupling between \(\tilde{\mathbb{P}}\) and \(\mathbb{P}_N\). Note that the loss change \(\bm{l}\) is upper bounded by the maximal rate \(\Delta_{\theta}^{\max}\), which is in turn upper bounded by its least concave majorant \(\mathcal{C}_{f^{\max}}\). Since \(\mathcal{C}_{f^{\max}}\) is concave, we can apply Jensen's inequality to move the majorant outside the integral \(\int_{\mathcal{Z}\times\mathcal{Z}_N} \) to get the desired conclusion. The full proof of Theorem~\ref{thm:main0} is given in Appendix~\ref{proof:main0}. Notably, the nature of \(u\in[\epsilon^p,\infty) \) also emerges from another line of related work \cite{liu2025wasserstein,chu2026wasserstein}. \hfill\Halmos

\begin{figure}[h]
    \centering
    \includegraphics[width=0.9\linewidth]{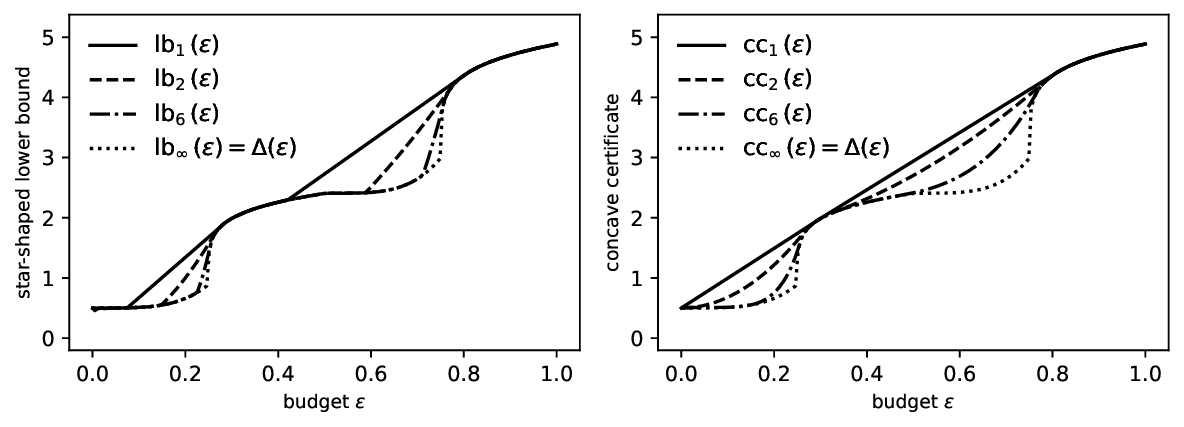}
    \caption{Illustration of Theorem~\ref{thm:main0}. Given a rate \(\Delta(t)\) (dotted curve) in Figure~\ref{fig:rate}, this plot visualizes the geometric construction of the proposed lower bound \eqref{eq:main-L} (least star-shaped majorant - Left) and upper bound \eqref{eq:main-U} (least concave majorant - Right). The \(p\)-Wasserstein perturbation elevates the risk gap $\mathcal{R}_{p}(\epsilon)-\hat{\mathcal{R}}$ by at least $\operatorname{lb}_{p}(\epsilon)$ and at most $\operatorname{cc}_{p}(\epsilon)$.}
    \label{fig:p-dymanic}
\end{figure}

Roughly speaking, Theorem~\ref{thm:main0} establishes that the robust-empirical risk gap \( \mathcal{R}_{p}(\epsilon) -\hat{\mathcal{R}} \) is bounded between the average of the \textit{least star-shaped majorants} \(\operatorname{lb}_p(\epsilon) \) of individual growth rates  and the \textit{least concave majorant} \(\operatorname{cc}_p(\epsilon) \)  of the maximal growth rate of the loss function. We emphasize that this geometric framework does not require the loss function \(\bm{l}\) to be convex/differentiable/Lipschitz, or the cost function \(d\) to be a metric, or the domain  \(\mathcal{Z}\) to be bounded. We illustrate Theorem~\ref{thm:main0} in a special case where \(N=1\) and \(\Delta_{\theta}(Z^{(i)},t)=\Delta_{\theta}^{\max}(t) = \Delta \) being continuous in Figure~\ref{fig:p-dymanic}.

\subsection{Dynamic of the DR risk \(\mathcal{R}_p\)} \label{sec:DR}

In Figure~\ref{fig:p-dymanic}, we can see that when \(p\) increases, both lower bound and upper bound decrease, which is proved in the following Corollary~\ref{cor:p-dynamic}. To the best of our knowledge, this \(p\)-dynamic has not been previously explored in the literature for 
DR risk estimations.

\begin{corollary}
    [\(p\)-dynamic of \(\mathcal{R}_{p}(\epsilon) \)] \label{cor:p-dynamic} Take \(1\leq p\leq p_2\leq \infty\), it is known that \(p\)-Wasserstein uncertainty set is larger than  the  \(p_2\)-Wasserstein uncertainty set (Lemma~\ref{lem:Wp-ineq}). Consequently, \(\mathcal{R}_{p}(\epsilon)\geq \mathcal{R}_{p_2}(\epsilon) \). Our proposed geometric bounds preserve this ordering as well. That is,
    \begin{itemize}
        \item \(\operatorname{lb}_1(\epsilon)\geq \operatorname{lb}_p(\epsilon) \geq \operatorname{lb}_{p_2}(\epsilon) \geq \operatorname{lb}_{\infty}(\epsilon)=\sum_{i=1}^{N}\mu_i  \Delta_{\theta}(Z^{(i)},\epsilon) \), and
        \item \(\operatorname{cc}_1(\epsilon)\geq \operatorname{cc}_p(\epsilon) \geq \operatorname{cc}_{p_2}(\epsilon) \geq \operatorname{cc}_{\infty}(\epsilon) =\sum_{i=1}^{N}\mu_i \lim\limits_{t\rightarrow\epsilon^{+}}\Delta_{\theta}(Z^{(i)},t) \geq \operatorname{lb}_{\infty}(\epsilon)  \).        
    \end{itemize}
\end{corollary}

\noindent\textit{Proof.} By Lemma~\ref{lem:calculation}, \(\mathcal{S}_{f^{(i)}}(\epsilon^p) = \sup_{u\geq \epsilon^p}\frac{\epsilon^p f^{(i)}(u)}{u} =  \sup_{t\geq \epsilon}\frac{\epsilon^p \Delta_{\theta}(Z^{(i)},t) }{t^p}\). For any fixed $t \ge \epsilon$, the function $ p\mapsto (\frac{\epsilon}{t})^p$ is non-increasing since $\frac{\epsilon}{t} \le 1$. Therefore, \(\frac{\epsilon \Delta_{\theta}(Z^{(i)},t)}{t} \geq  \frac{\epsilon^p \Delta_{\theta}(Z^{(i)},t)}{t^p} \geq \frac{\epsilon^{p_2} \Delta_{\theta}(Z^{(i)},t)}{t^{p_2}} \). Taking supremum on \(t\in[\epsilon,\infty) \), we have the first conclusion. Next, let \(\mathcal{C}_f \) and \(\mathcal{C}_{f_2} \) be the least concave majorants of \(f\colon t\mapsto \Delta_{\theta}^{\max}(t^{1/p}) \) and \(f_2\colon t\mapsto \Delta_{\theta}^{\max}(t^{1/p_2}) \), respectively. By Lemma~\ref{lem:concave-function} and Lemma~\ref{lem:univariate-majorant},  \(\mathcal{C}_f (t^{p/p_2}) \) is concave. In addition, \(\mathcal{C}_f(t^{p/p_2}) \geq \Delta_{\theta}^{\max}(t^{p/p_2 \times 1/p}) = \Delta_{\theta}^{\max}(t^{1/p_2})\). Thus, \(\mathcal{C}_f(t^{p/p_2}) \geq \mathcal{C}_{f_2}(t) \). Choose \(t=\epsilon^{p_2} \), one has \(\mathcal{C}_f(\epsilon^p) \geq \mathcal{C}_{f_2}(\epsilon^{p_2}) \). Therefore, \(\operatorname{cc}_{p}\geq \operatorname{cc}_{p_2} \). \hfill\Halmos

Beyond characterizing the magnitude of the robust gap, our analysis allows us to identify the conditions under which the distributionally robust loss remains finite: either the robustness certificate is finite across the entire domain, or it diverges everywhere. Existing literature typically addresses the finiteness of the DR risk through the lens of strong duality \citep{zhang2022short, zhen2025unified} or equilibrium theory \citep{shafiee2025wasserstein}.

\begin{corollary}
    [Finiteness of \(\mathcal{R}_p\)] \label{cor:finite-Rp} Given \(1\leq p<\infty\), then exactly one of the following two cases must occur.
    \begin{itemize}
        \item \(\operatorname{lb}_p(\epsilon)=\operatorname{cc}_p(\epsilon)=\infty\) for any \(\epsilon>0\).
        \item \(\operatorname{lb}_p(\epsilon)<\infty\) and \(\operatorname{cc}_p(\epsilon)<\infty\) for any \(\epsilon>0\).
    \end{itemize}
    (This claim is not true for \(p=\infty\) as \(\lim\limits_{t\rightarrow\epsilon^{+}}\Delta_{\theta}(Z^{(i)},t)\) could be infinite even though \(\Delta_{\theta}(Z^{(i)},\epsilon)<\infty\).)
\end{corollary}

\noindent\textit{Proof.} Suppose that \(\operatorname{lb}_p(\epsilon)=\infty\) for any \(\epsilon>0\). Since \(\operatorname{lb}_p(\epsilon)\leq \operatorname{cc}_p(\epsilon)\), it implies that \(\operatorname{cc}_p(\epsilon)=\infty\) for any \(\epsilon>0\). Suppose otherwise that \(\operatorname{lb}_p(\epsilon) = \sum_i \mu_i \mathcal{S}_{f^{(i)}}(\tilde{\epsilon}^p)<\infty \) for some \(\tilde{\epsilon}>0\).  Since \(f^{\max}=\max_i f^{(i)}\), this implies \(\mathcal{S}_{f^{\max}}(\tilde{\epsilon}^p)  <\infty\). As \(\mathcal{S}_{f^{\max}}(t)\) is star-shaped, \(\frac{\mathcal{S}_{f^{\max}}(t)}{t} \) is non-increasing on \((0,\infty) \). Then \(\frac{\mathcal{S}_{f^{\max}}(t)}{t} \leq \frac{\mathcal{S}_{f^{\max}}(\tilde{\epsilon}^p)}{\tilde{\epsilon}^p} = \tilde{a} <\infty  \) and thus \(f^{\max}(t) \leq \mathcal{S}_{f^{\max}}(t) \leq \tilde{a}t \) for any \(t\geq \tilde{\epsilon}^p\). Note that \(f^{\max}\) is non-decreasing on \([0,\infty) \), it follows that  \(f^{\max}(t) \leq  \tilde{a}t + f^{\max}(\tilde{\epsilon}^p) \) for any \(t\geq 0\). That is to say \(f^{\max}\) is upper bounded by an affine (concave) function, thus \(\mathcal{C}_{f^{\max}}(t)<\infty \) for any \(t\geq 0\). \hfill\Halmos

Corollary~\ref{cor:p-dynamic} and Corollary~\ref{cor:finite-Rp} demonstrate that our proposed bounds not only track how the exponent \(p\) (as in the \(p\)-Wasserstein) dictates the magnitude of the robust risk \(\mathcal{R}_p\), but also provide a rigorous criterion for its finiteness. By evaluating whether the loss growth is  compatible with the exponent \(p\), we illustrate this transition and demonstrate the superiority of our approach over traditional convexity, differentiability or Lipschitz  certificates in the following example.

\begin{figure}[H]
    \centering
    \begin{subfigure}[b]{0.3\linewidth}
        \centering
        \includegraphics[width=\linewidth]{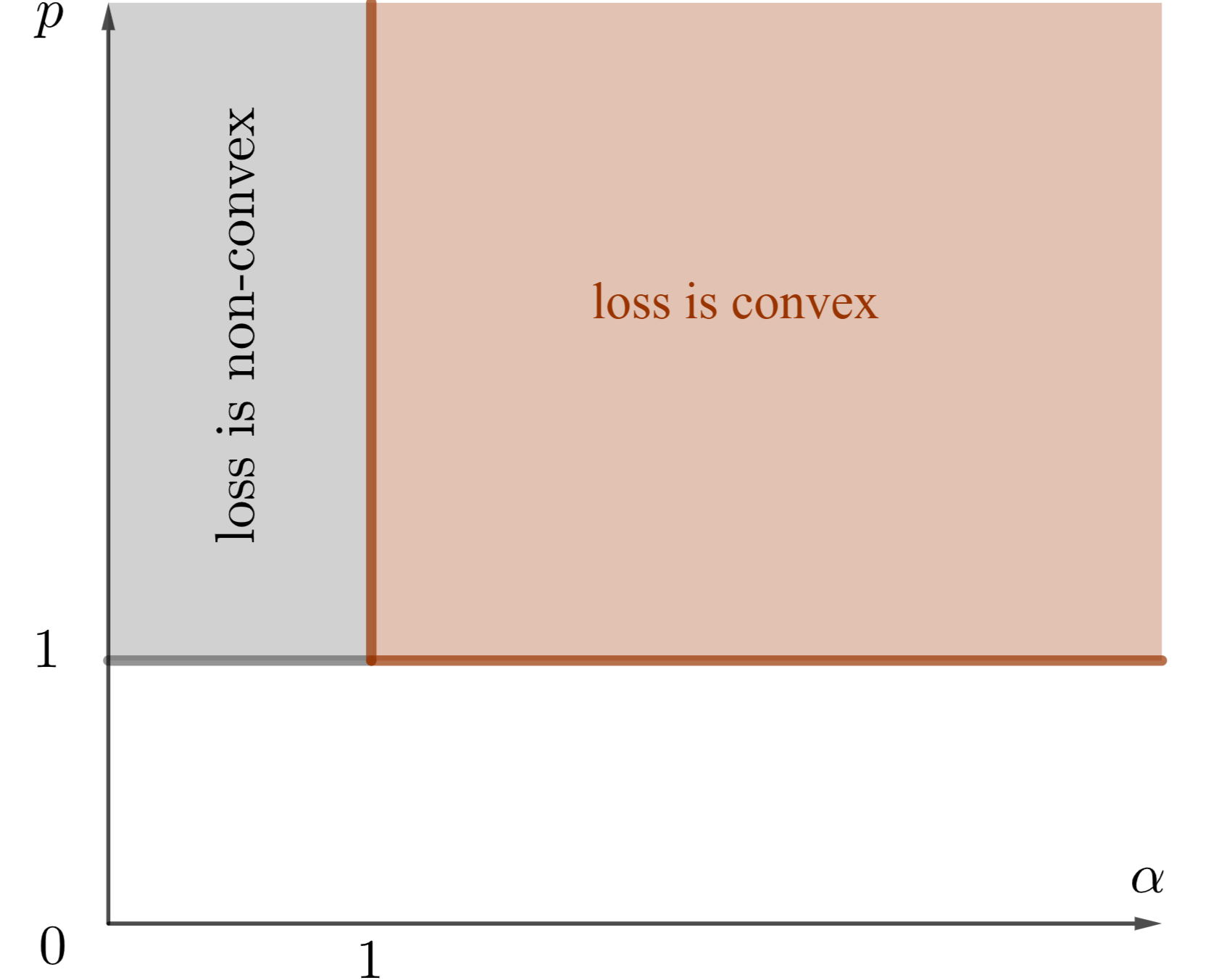}
        \caption{existing convexity certificate}
        \label{fig:conv}
    \end{subfigure}
    \hfill 
    \begin{subfigure}[b]{0.3\linewidth}
        \centering
        \includegraphics[width=\linewidth]{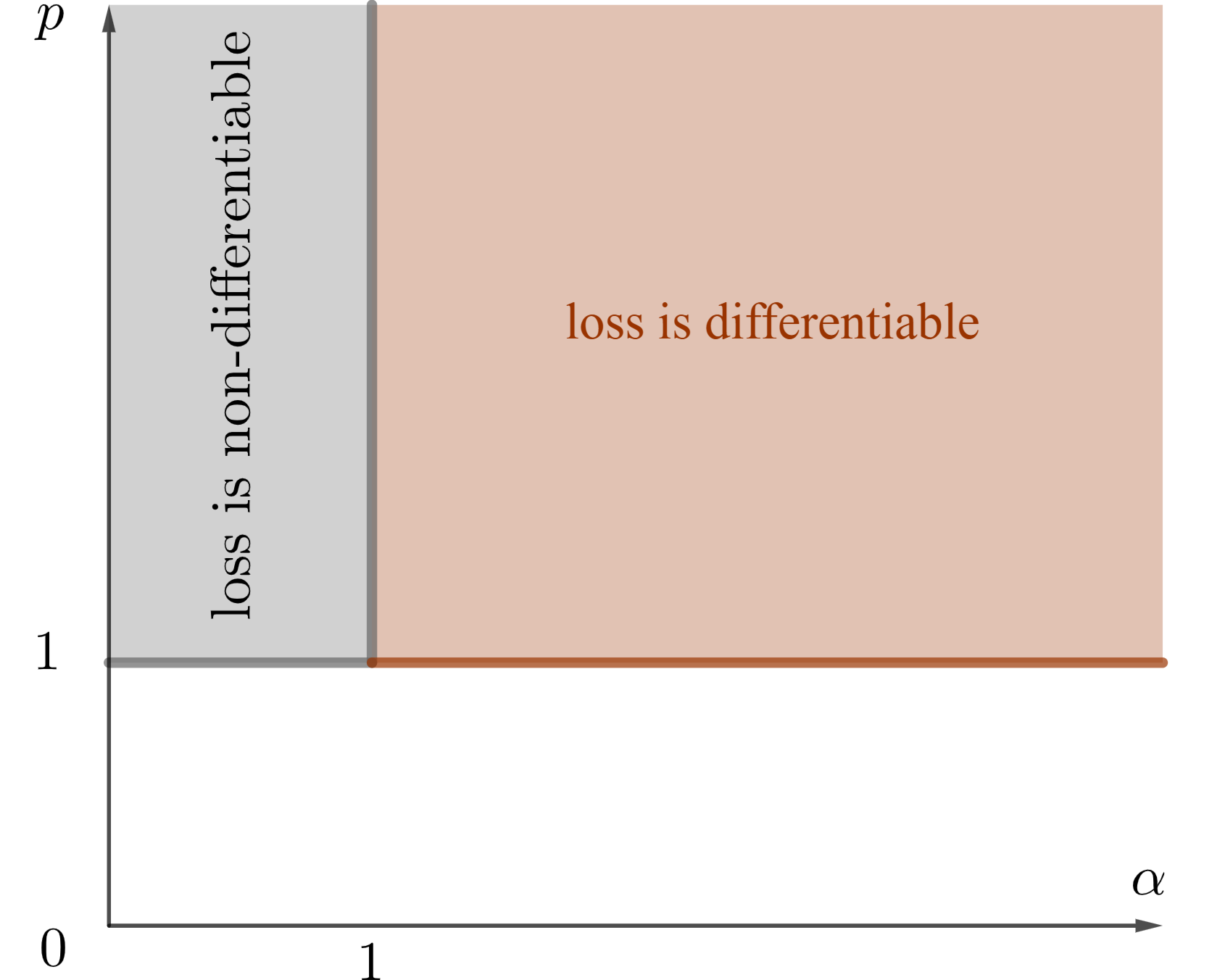}
        \caption{existing differentiability cert.}
        \label{fig:lip}
    \end{subfigure}
    \hfill
    \begin{subfigure}[b]{0.3\linewidth}
        \centering
        \includegraphics[width=\linewidth]{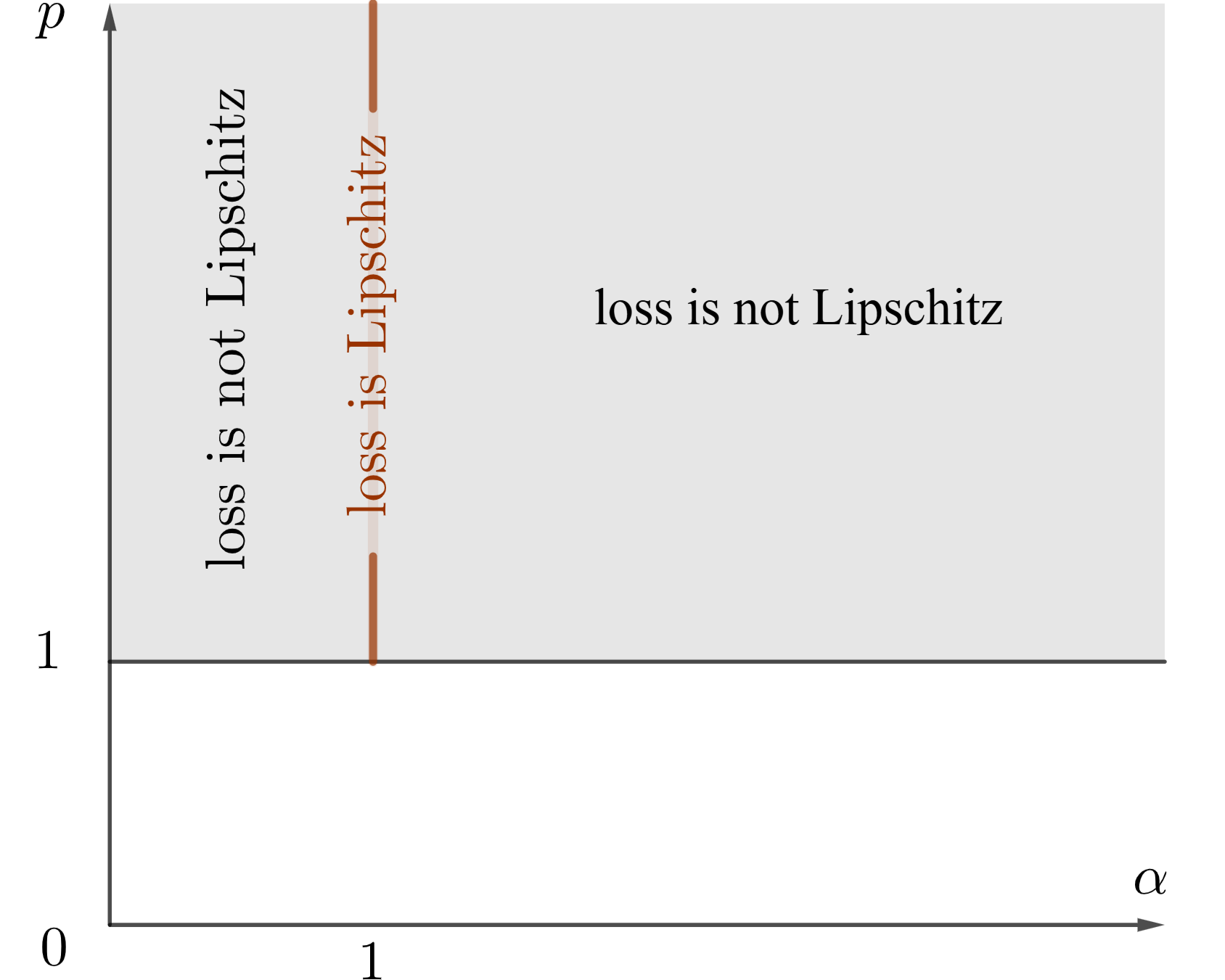}
        \caption{existing Lipschitz certificate}
        \label{fig:diff}
    \end{subfigure}
    \begin{subfigure}[b]{0.3\linewidth}
        \centering
        \includegraphics[width=\linewidth]{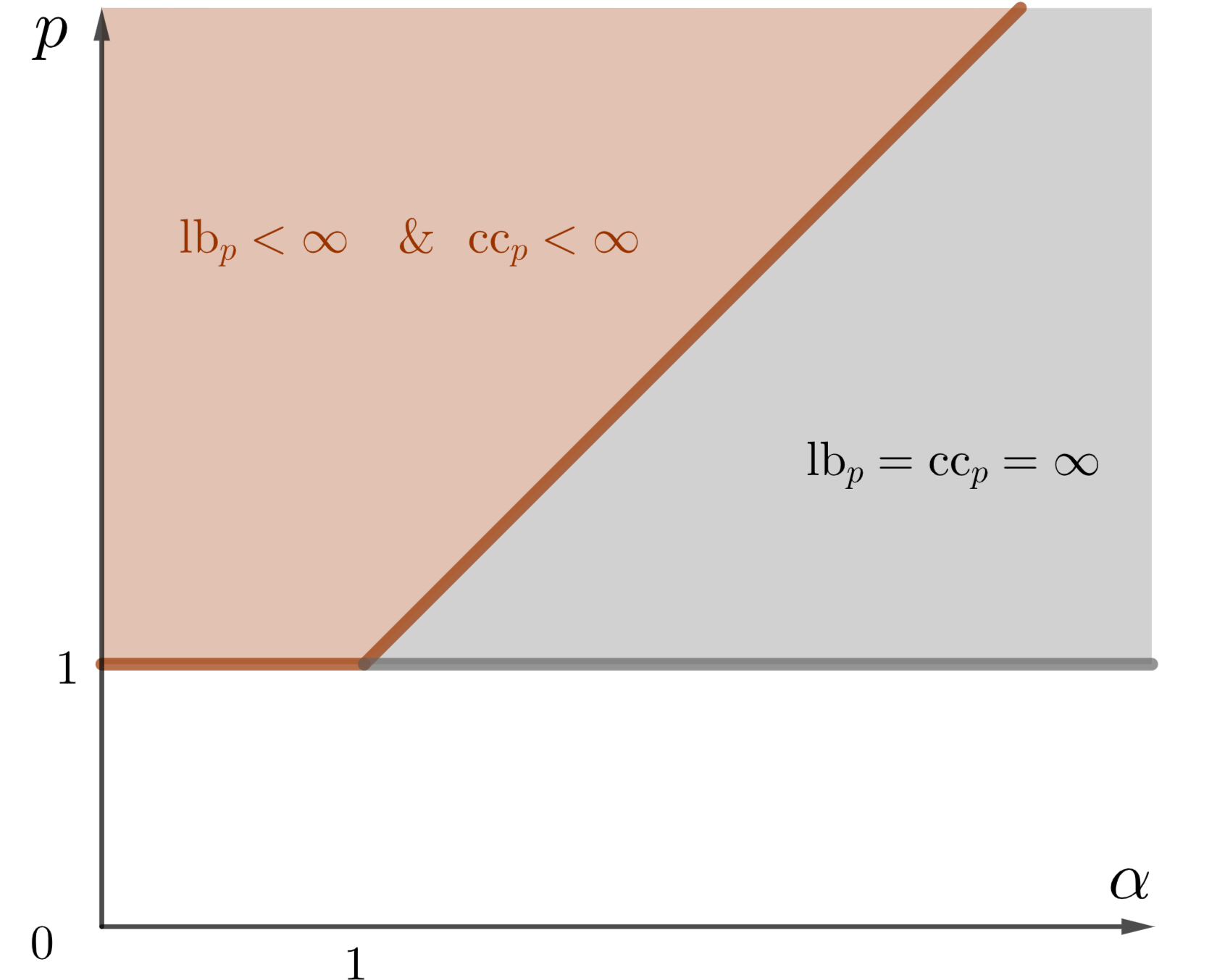}
        \caption{proposed concave certificate}
        \label{fig:delta}
    \end{subfigure}
    \begin{subfigure}[b]{0.3\linewidth}
        \centering
        \includegraphics[width=\linewidth]{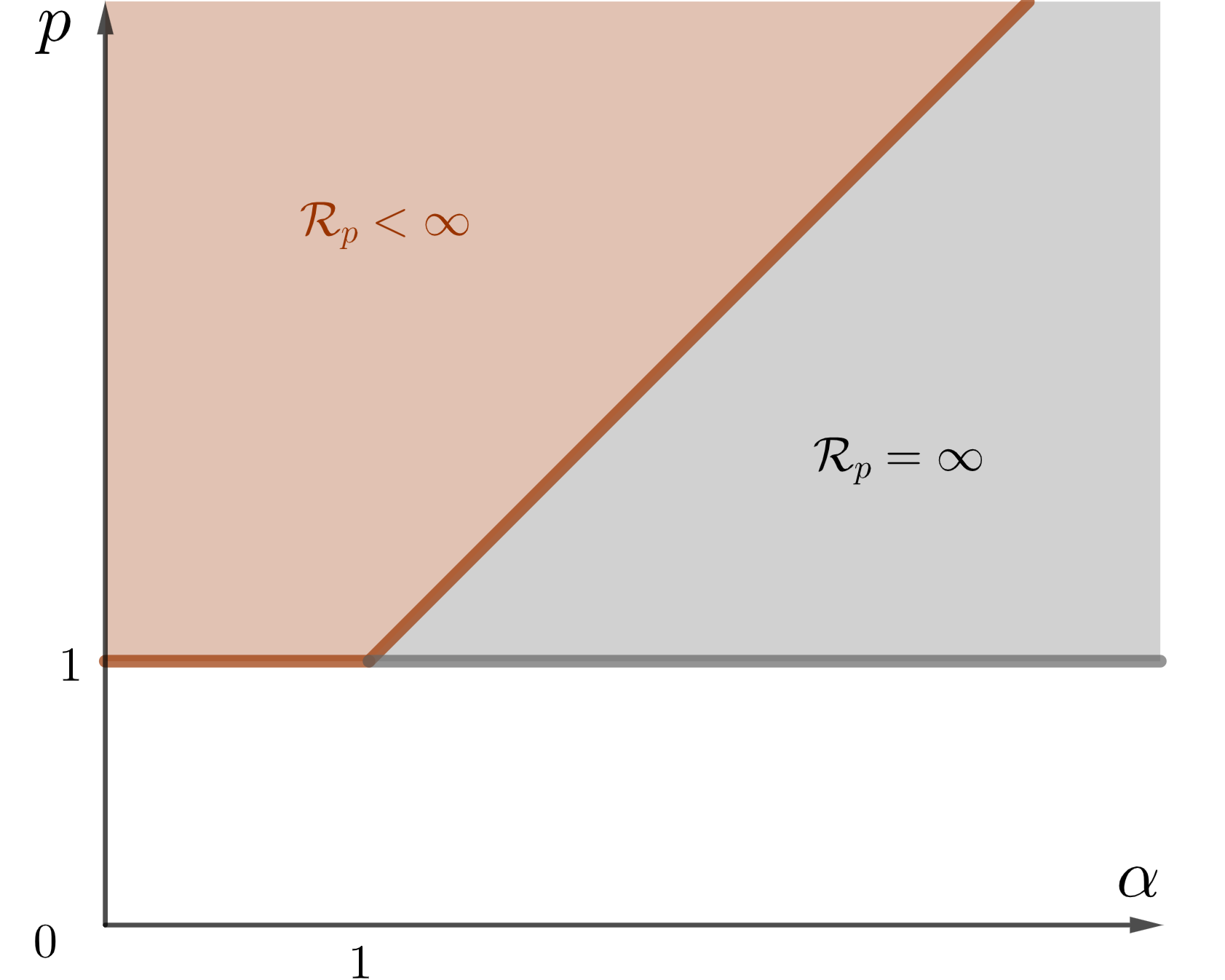}
        \caption{actual robustness}
        \label{fig:radius}
    \end{subfigure}
    \caption{Illustration of Example~\ref{exam:unbounded} certificating \(\bm{l}(z;\theta) = \abs{y-\inprod{x,\theta}}^{\alpha}\). Compared to existing convexity (a), differentiability (b), and Lipschitz (c) certificates, only our proposed concave certificate (d) is able to characterize the domain of robustness (e) exactly.}
    \label{fig:unbounded}
\end{figure}

 \begin{example}[Tightness of Theorem~\ref{thm:main0}.]\label{exam:unbounded}
    (Figure~\ref{fig:unbounded}) Suppose that the loss function \(\bm{l}\colon \mathbb{R}^{n+1} \times \mathbb{R}^{n}\rightarrow\mathbb{R} \) is defined by \(\bm{l}(z;\theta) = \abs{y-\inprod{x,\theta}}^{\alpha}\) for some given \(\alpha\in(0,\infty) \) and \(z=(x,y)\); and the cost function \(d\colon \mathbb{R}^{n+1}\times\mathbb{R}^{n+1}\rightarrow[0,\infty] \) is defined by \( d(z',z) = \norm{x'-x}_{r} + \infty\abs{y'-y} \) where \(r\geq 1\) with the convention \(\infty\cdot 0 = 0\). Denote \(\hat{c}_{i}\coloneqq Y^{(i)} - \inprod{X^{(i)},\theta}  \), then the individual rate \(\Delta_{\theta}(Z^{(i)},t)\) satisfies that
    \begin{equation*}
        t^{\alpha}\norm{\theta}_s^{\alpha}\leq \Delta_{\theta}(Z^{(i)},t) \leq  \left( \abs{\hat{c}_{i}} + t\norm{\theta}_s \right)^{\alpha}-\abs{\hat{c}_i}^{\alpha},
    \end{equation*}
    where \(1/r + 1/s = 1 \).Therefore for any \(\epsilon>0\),
    \begin{itemize}
    \item if \(p\in[1,\infty)\cap [1,\alpha)  \) then   \(\operatorname{lb}_{p}(\epsilon)=\operatorname{cc}_{p}(\epsilon)=\mathcal{R}_{p}(\epsilon) =\infty \); and
    \item if \(p \in [1,\infty)\cap[\alpha,\infty) \) or \(p=\infty \) then  \(\operatorname{lb}_{p}(\epsilon)\leq \mathcal{R}_{p}(\epsilon) \leq \operatorname{cc}_{p}(\epsilon)\leq \infty \).
    \end{itemize}
\end{example}

\noindent\textit{Proof.} By Definition~\ref{def:rate}, the individual rate is computed by
\begin{equation*}
    \Delta_{\theta}(z,t) =       \sup_{x'\colon\norm{x'-x}_{r} \leq  t}       \left\{\abs{y-\inprod{x',\theta}}^{\alpha} - \abs{y-\inprod{x,\theta}}^{\alpha}  \right\}.
\end{equation*}
Since \(\abs{y-\inprod{x',\theta}}\leq  \abs{y-\inprod{x,\theta}} + \norm{x'-x}_{r}\norm{\theta}_s\), we obtain \(\Delta_{\theta}(Z^{(i)},t) \leq  \left( \abs{\hat{c}_{i}} + t\norm{\theta}_s \right)^{\alpha} -\abs{\hat{c}_{i}}^{\alpha} \) for any \(t>0\). On the other hand, choose \(X'= X^{(i)} - \operatorname{sign}(\hat{c}_i) t\xi \) where \(\xi\coloneqq\arg\max_{\norm{\xi}_{r}=1} \inprod{\xi,\theta} \). Then \(\Delta_{\theta}(Z^{(i)},t)  \geq \abs{\hat{c}_i + \operatorname{sign}(\hat{c}_i) t\norm{\theta}_s }^{\alpha} - \abs{\hat{c}_i}^{\alpha} \geq t^{\alpha}\norm{\theta}_s^{\alpha}  \).

Suppose that \(p\in[1,\alpha)\cap [1,\infty) \), then \(p<\alpha \). Let \(g^{(i)}(t)= (t^{1/p}) ^{\alpha}\norm{\theta}_s^{\alpha}\), then \(g^{(i)}(t)\leq f^{(i)}(t) \) where \(f^{(i)}=\Delta_{\theta}(Z^{(i)},t^{1/p})\). By Lemma~\ref{lem:calculation}, \( \mathcal{S}_{g^{(i)}}(t)= \sup_{u\in[t,\infty)} \frac{t g^{(i)}(u) }{u}=\sup_{u\in[t,\infty)} t u^{\alpha/p-1}\norm{\theta}_s^{\alpha} = \infty \) for any \(t>0\). Therefore, \(\mathcal{S}_{f^{(i)}}(t) \geq \mathcal{S}_{g^{(i)}}(t) = \infty \) and \(\operatorname{lb}_p=\operatorname{cc}_p= \mathcal{R}_p=\infty\).

Otherwise, if \(p \in [\alpha,\infty) \cap [1,\infty) \) then \(\Delta_{\theta}^{\max}(t) = \sup_{\hat{z}\in\mathcal{Z}_N}\Delta_{\theta}(\hat{z},t) \leq  \left(\hat{C}+ t\norm{\theta}_s \right)^{\alpha} \)  where \(\hat{C}\coloneqq\max_{i=1}^{N}\{\hat{c}_{i}\} \). Let \(f^{\max}(t) = \Delta_{\theta}^{\max}(t^{1/p})\), then   \(f^{\max}(t)\leq \left(\hat{C}+ t^{1/p} \norm{\theta}_s \right)^{\alpha}\), which is concave. Thus \(\mathcal{C}_{f^{\max}}(t)  \leq\left(\hat{C}+ t^{1/p} \norm{\theta}_s \right)^{\alpha} <\infty \). By Theorem~\ref{thm:main0}, we have that  \( \operatorname{lb}_p \leq\mathcal{R}_{p} -\hat{R}\leq \operatorname{cc}_p\leq \infty \). The case of \(p=\infty\) trivially follows since \(\Delta_{\theta}^{\max}\) is finite. \hfill\Halmos

We conclude this section by noting that the existing Lipschitz certificate \citep{blanchet2019quantifying,blanchet2019robust,an2021generalization,gao2022wasserstein} is a direct consequence of Theorem~\ref{thm:main0}. Furthermore, our analysis allows us to remove the boundedness assumption on the domain \(\mathcal{Z}\) required by \citet[Lemma 2]{gao2023distributionally} while showing that the DR risk remains lower bounded by its scalar growth rate.

\begin{corollary}[Lipschitz Certificate]\label{cor:Lip-ub}
    Given Notation~\ref{nota:main}, if  \(\abs{\bm{l}(z';\theta) - \bm{l}(z;\theta) }\leq \operatorname{Lip}\times d(z',z) \) for any \(z',z\in\mathcal{Z}\)  then \(\mathcal{R}_{p}(\epsilon)\leq \hat{\mathcal{R}}+ \operatorname{Lip}\times \epsilon \) for any \(p\in[1,\infty]\). Besides, if \(\sup_{t\in[\epsilon,\infty) } \frac{\Delta_{\theta}(Z^{(i)},t)}{t^p} \geq \kappa \) then \(\mathcal{R}_{p}(\epsilon)\geq  \hat{\mathcal{R}}+ \kappa\times \epsilon^p \).
\end{corollary}

\subsection{Distributional Robust Classifier}
We recall the standard definition of (point-wise) robust classifier \cite{pal2023adversarial,pal2024certified}, which is conceptually originated from adversarial studies \cite{madry2018towards,schmidt2018adversarially,cullina2018pac,diochnos2018adversarial}.

\begin{definition}[Robust Classifier] \label{def:robust-classifier}    
Given $\mathcal{X}\subseteq\mathbb{R}^{n}$ and $\mathcal{Y} = \{1,2,\dots,m\}$, a classifier $f_{\theta}\colon \mathcal{X}\rightarrow\mathcal{Y}$ is called $(\epsilon,\delta)$-robust with respect to $\mathbb{P}_N$ if 
    \begin{equation*}
        \operatorname{Prob}_{(\hat{x},\hat{y})\sim \mathbb{P}_N} \left( \exists \tilde{x} \text{ s.t. } d_{\mathcal{X}}(\tilde{x},\hat{x})\leq   \epsilon \text{ and } f_{\theta}(\tilde{x})\ne\hat{y}  \right) \leq \delta.
    \end{equation*}
\end{definition}

Roughly speaking, this inequality says that the probability of an empirical data point being vulnerable to an \(\epsilon\)-adversarial perturbation is at most \(\delta\). Interestingly, this concept coincides with our above lower bound $\operatorname{lb}_{\infty}(\epsilon)$  in Theorem~\ref{thm:main0} plus empirical loss \(\hat{\mathcal{R}}\) under the zero-one loss setting. 

\begin{lemma} Given Notation~\ref{nota:main}, let $\bm{l}$ be the zero-one loss given by $\bm{l}\big(z = (x,y);\theta\big) = \bm{1} \left( {f_{\theta}(x) \ne  y }\right) \in \{0,1\}$ and $d$ given by $d(z',z) = d_{\mathcal{X}}(x',x) + \infty\cdot\abs{y'-y}$. Then \(f_{\theta}\) is \((\epsilon,\delta)\)-robust with respect to  \(\mathbb{P}_N\)  if and only if \(\operatorname{lb}_{\infty}(\epsilon) +\hat{\mathcal{R}} \leq \delta.\)
\end{lemma}

\noindent\textit{Proof.} A direct manipulation of $\operatorname{lb}_{\infty}(\epsilon)$ gives
\begin{equation*}
    \begin{array}{ll}
        \operatorname{lb}_{\infty}(\epsilon) +\hat{\mathcal{R}} 
        &= \sum_{i=1}^N \mu_i\Delta_{\theta}(Z^{(i)},\epsilon) + \bm{l}({Z}^{(i)};\theta) = \sum_{i=1}^N \mu_i\sup_{d(\tilde{Z}^{(i)},Z^{(i)})\leq   \epsilon} \bm{l}(\tilde{Z}^{(i)};\theta)\\
        &= \sum_{i=1}^N \mu_i\sup_{d_{\mathcal{X}}(\tilde{X}^{(i)},X^{(i)})\leq   \epsilon} \bm{l}(\tilde{X}^{(i)}, Y^{(i)};\theta)\\
        &= \sum_{i=1}^N \mu_i\bm{1} \left( \exists \tilde{X}^{(i)} \text{ s.t. } d_{\mathcal{X}}(\tilde{X}^{(i)},X^{(i)})\leq   \epsilon \text{ and } f_{\theta}(\tilde{X}^{(i)})\ne Y^{(i)}   \right) \\
        & = \operatorname{Prob}_{\mathbb{P}_N} \left( \exists \tilde{X}^{(i)} \text{ s.t. } d_{\mathcal{X}}(\tilde{X}^{(i)},X^{(i)})\leq   \epsilon \text{ and } f_{\theta}(\tilde{X}^{(i)})\ne Y^{(i)}  \right). \hfill\Halmos
    \end{array}
\end{equation*}

Through the lens of \(\operatorname{lb}_{\infty}\), we derive the necessary and sufficient condition for a classifier  \(f_\theta\) to be $(\epsilon,\delta)$-robustness as follows, with the proof provided in Appendix~\ref{proof:ExactCD}. We emphasize that we do not require \(d_{\mathcal{X}}\) to satisfy the triangle inequality, nor \(\mathcal{X}\) to be bounded.

\begin{proposition}[Exactly Concentrated Distribution] \label{prop:ExactCD} 
For any set $\Omega\subseteq\mathcal{X}$, define its $\epsilon$-expansion as $\Omega^{+\epsilon}\coloneqq \{ x \in \mathcal{X} \colon \exists \bar{x} \in\Omega  \text{ s.t. }  d_{\mathcal{X}}(x,\bar{x}) \leq \epsilon \}  $. We say that $\mathbb{P}_N$ is $(\epsilon,\delta)$-exactly concentrated (with respect to classifier $f_\theta$) if there exists  $\{\Omega_k\}_{k=1}^m$ with $\Omega_k\subseteq \{X^{(i)} \colon Y^{(i)}=k \}$ such that
\begin{itemize}
    \item[\(\operatorname{(i)}\)]   $\mathbb{P}_N(\cup_k \{X \in \Omega_k, Y=k\})=  \sum_{k=1}^m \sum_{i \colon X^{(i)} \in \Omega_k} \mu_i  \geq 1 - \delta$, and \hfill ($\delta$-coverage)
    \item[\(\operatorname{(ii)}\)] $\Omega_k^{+\epsilon} \subseteq \mathcal{D}_{f_\theta,k} $ for any \(k=1,2,\dots,m\) where $\mathcal{D}_{f_\theta,k}  \coloneqq \{x \in \mathcal{X} \colon  f_\theta(x) = k\}$. \hfill ($\epsilon$-immunity)
\end{itemize}
 Then $f_{\theta}$ is $(\epsilon,\delta)$-robust \textbf{if and only if} $\mathbb{P}_N$ is $(\epsilon,\delta)$-exactly concentrated with respect to classifier $f_\theta$. In particular, if \(\mu_i=\frac{1}{N}\) for \(i=1,\dots,N\), then \(\operatorname{(i)}\) becomes $ \frac{1}{N} \sum_{k=1}^m \#\Omega_k \geq 1 - \delta$.
\end{proposition}
\begin{figure}[htbp]
    \centering
    % Scales the diagram to exactly 0.8 of the text width
    \resizebox{0.7\textwidth}{!}{%
    \begin{tikzpicture}[
        x=1cm, y=1cm,
        % Base font applied to every node text (even when nodes specify their own size/style)
        basefont/.style={execute at begin node=\small},
        every node/.append style={basefont}
    ]

        % 2. Enforce a strict 2:1 aspect ratio (Width = 12, Height = 6)
        \clip (0,0) rectangle (12,6);

        % Backgrounds for continuous decision regions D1 and D2 (Non-class regions in very faint gray)
        \fill[blue!2] (0,0) rectangle (6,6);
        \fill[brown!2] (6,0) rectangle (12,6);

        % Define epsilon radius
        \def\eps{0.8}

        % ---------------------------------------------------------
        % Epsilon Expansions (\Omega_k^{+\epsilon})
        % Drawn FIRST so they sit behind the points.
        % ---------------------------------------------------------

        % Class 1 Halos (Blue) - 3 robust points
        \foreach \p in {(2.5, 4), (3.5, 5), (3.5, 2.5)} {
            \fill[blue!15] \p circle (\eps);
            \draw[blue!30, thick] \p circle (\eps);
        }

        % Class 2 Halos (Brown) - 4 robust points
        \foreach \p in {(8.5, 4.5), (9.5, 3), (10.5, 4.8), (11, 1.5)} {
            \fill[brown!15] \p circle (\eps);
            \draw[brown!30, thick] \p circle (\eps);
        }

        % ---------------------------------------------------------
        % Model Decision Boundary & Top Labels (Strictly Black)
        % ---------------------------------------------------------
        \draw[ultra thick, dashed, black] (6,-1) -- (6,7);

        % Vertical decision boundary label
        \node[black, font=\bfseries, rotate=90] at (5.7, 1.8) {decision boundary};

        % Top labels next to the dashed line
        \node[black, font=\bfseries, anchor=east] at (5.7, 5.5) {$\mathcal{D}_{f_\theta, 1}$};
        \node[black, font=\bfseries, anchor=west] at (6.3, 5.5) {$\mathcal{D}_{f_\theta, 2}$};

        % Epsilon radius indicator (Black)
        \draw[thick, <->, >=stealth, black] (2.5, 4) -- +(180:\eps) node[midway, above, black] {$\epsilon$};

        % ---------------------------------------------------------
        % Empirical Data Points: Class 1 (Blue)
        % ---------------------------------------------------------
        % Robust Points (Discrete set \Omega_1) - 3 points
        \foreach \p in {(2.5, 4), (3.5, 5), (3.5, 2.5)} {
            \node[circle, fill=blue, inner sep=1.5pt] at \p {};
        }

        % Vulnerable Point 1 (Outside \Omega_1, crossing boundary)
        \node[circle, draw=blue, fill=white, thick, inner sep=1.5pt] (v1) at (5.4, 4.0) {};
        \draw[dashed, black, thick, fill=black!5, fill opacity=0.3] (v1) circle (\eps);

        % Vulnerable Point 2 (Added extra blue vulnerable point)
        \node[circle, draw=blue, fill=white, thick, inner sep=1.5pt] (v1b) at (5.5, 2.0) {};
        \draw[dashed, black, thick, fill=black!5, fill opacity=0.3] (v1b) circle (\eps);

        % ---------------------------------------------------------
        % Empirical Data Points: Class 2 (Brown)
        % ---------------------------------------------------------
        % Robust Points (Discrete set \Omega_2) - 4 points
        \foreach \p in {(8.5, 4.5), (9.5, 3), (10.5, 4.8), (11, 1.5)} {
            \node[circle, fill=brown, inner sep=1.5pt] at \p {};
        }

        % Vulnerable Point (Outside \Omega_2, crossing boundary)
        \node[circle, draw=brown, fill=white, thick, inner sep=1.5pt] (v2) at (6.6, 2.5) {};
        \draw[dashed, black, thick, fill=black!5, fill opacity=0.3] (v2) circle (\eps);

        % ---------------------------------------------------------
        % Bottom Legends (Visual Replicas)
        % ---------------------------------------------------------

        % Class 1 Legend (Blue)
        \node[circle, fill=blue, inner sep=1.5pt] (leg_dot1) at (2, 1.2) {};
        \node[anchor=west, black, font=\bfseries] at (2.2, 1.2) {$= \Omega_1$};

        \fill[blue!15] (1.8, 0.4) rectangle (2.2, 0.8);
        \draw[blue!30, thick] (1.8, 0.4) rectangle (2.2, 0.8);
        \node[anchor=west, black, font=\bfseries] at (2.2, 0.6) {$= \Omega_1^{+\epsilon}$};

        % Class 2 Legend (Brown)
        \node[circle, fill=brown, inner sep=1.5pt] (leg_dot2) at (8, 1.2) {};
        \node[anchor=west, black, font=\bfseries] at (8.2, 1.2) {$= \Omega_2$};

        \fill[brown!15] (7.8, 0.4) rectangle (8.2, 0.8);
        \draw[brown!30, thick] (7.8, 0.4) rectangle (8.2, 0.8);
        \node[anchor=west, black, font=\bfseries] at (8.2, 0.6) {$= \Omega_2^{+\epsilon}$};

    \end{tikzpicture}%
    }
    \caption{Illustration of the Exact CD condition when \(\mu_i=\frac{1}{N} = \frac{1}{10}\). The \(\delta\)-coverage property means that there are at least \(N(1-\delta) \) points being covered. The \(\epsilon\)-immunity means that all covered points are not flipped by \(f_{\theta}\) under \(\epsilon\)-perturbation. In this figure, \(f_{\theta}\) is \((\epsilon,0.3)\)-robust but not \((\epsilon,0.299)\)-robust.}
    \label{fig:ECD}
\end{figure}

In the remainder of this section, we shall show that our Exact CD aligns with existing results in \cite{pal2023adversarial,pal2024certified}. Following these results, we assume that  \(d_{\mathcal{X}}\) satisfies the triangle inequality.

\subsubsection{$(\epsilon,\delta-\rho,\rho)$-Strong CD implies  $(\epsilon,\delta)$-Exact CD}

    Given $\rho\in(0,\delta)$, Recall that $\mathbb{P}_N$ is $(\epsilon,\delta-\rho,\rho)$-strongly concentrated (regardless of $f_\theta$) \cite{pal2023adversarial,pal2024certified} if there exists $\{\Omega_k\}_{k=1}^m$ with $\Omega_k\subseteq \{X^{(i)} \colon Y^{(i)}=k \}$ such that
    \begin{itemize}
        \item[\(\operatorname{(i')}\)] $\mathbb{P}_N^{k}(\Omega_k) = \mathbb{P}_N(X\in\Omega_k \mid Y= k) \geq 1 - \delta +\rho$ for any $k=1,2,\dots,m$, and 
        \item[\(\operatorname{(ii')}\)] $\mathbb{P}_N^{k}\left(\cup_{k'\ne k} \Omega_{k'}^{+2\epsilon} \right) = \mathbb{P}_N(\cup_{k'\ne k} \Omega_{k'}^{+2\epsilon} \mid Y= k) \leq \rho$ for any $k=1,2,\dots,m$.
    \end{itemize}
    We shall show that Strong CD implies Exact CD for some suitable classifier. Let $V_k \coloneqq \Omega_k \cap \left( \bigcup_{k' \neq k} \Omega_{k'}^{+2\epsilon} \right)$ be the set of points being covered by \(\Omega_k\) but also being a \(2\epsilon\)-perturbation of some \(k'\).  Define the refined subset $\tilde{\Omega}_k \coloneqq \Omega_k \setminus V_k$, then \(\mathbb{P}_N^k(V_k) \leq \rho \) and \(\mathbb{P}_N^k(\tilde{\Omega}_k) \geq (1-\delta+\rho) - \rho = 1-\delta\). Hence, \(\{\tilde{\Omega}_k\}_{k=1}^{N}\) satisfies the $\delta$-coverage condition $\operatorname{(i)}$  that $\mathbb{P}_N(\cup_k \{X \in \tilde{\Omega}_k, Y=k\}) = \sum_{k=1}^m \mathbb{P}_N^k(X\in\tilde{\Omega}_k) \mathbb{P}_N(Y= k) \geq 1- \delta$.

    To find a suitable \(f_{\theta}\) such that \(\{\tilde{\Omega}_k\}_{k=1}^{N}\) is \(\epsilon\)-immunity, suppose by contradiction there exists a point $x \in \tilde{\Omega}_k^{+\epsilon} \cap \tilde{\Omega}_{k'}^{+\epsilon}$ for $k \neq k'$. By the triangle inequality, there must be two empirical data points $x_k \in \tilde{\Omega}_k$ and $x_{k'} \in \tilde{\Omega}_{k'}$ such that $d_{\mathcal{X}}(x_k,x_{k'}) \leq 2\epsilon$. This implies $x_k \in \Omega_{k'}^{+2\epsilon}\subseteq V_k$, which  contradicts the construction of $\tilde{\Omega}_k$ that excluded $V_k$. Since the $\epsilon$-expansions $\tilde{\Omega}_k^{+\epsilon}$ are strictly disjoint across different classes, there exists a classifier $f_\theta$ such that $\tilde{\Omega}_k^{+\epsilon} \subseteq \mathcal{D}_{f_\theta, k}$ for all $k$, satisfying $\operatorname{(ii)}$. Finally, it is  worth noting that the slack budget $+2\epsilon$ arrived from the usage of the triangle inequality of \(d_{\mathcal{X}}\).

\subsubsection{\((\epsilon,\delta)\)-Exact CD implies \((\Phi(\epsilon,n,d_{\mathcal{X}}),\delta)\)-Standard CD}  Recall that \(\mathbb{P}\in\mathcal{P}(\mathcal{X})\) is \((\Phi,\delta) \)-concentrated (or localized) if there exists a subset \(\Omega\subseteq\mathcal{X}\) such that \(\mathbb{P}(\Omega)\geq 1-\delta \) and \(\operatorname{Vol}(\Omega)\leq  \Phi\).  It has been shown in \citet{pal2023adversarial,pal2024certified} that if \(f_{\theta}\) is \((\epsilon,\delta)\)-robust then 
\begin{itemize}
    \item[\(\operatorname{(i^*)}\)] there is at least one class \(\bar{k}\) such that the conditional distribution \(\mathbb{P}_N^{\bar{k}} = \mathbb{P}_N(\cdot\mid Y = \bar{k}) \) is \((\Phi_{\bar{k}},\delta) \)-concentrated, where \(\Phi_{\bar{k}}\) depends on \(\epsilon,n\) and \(d_{\mathcal{X}}\). \\
    In addition, if \( \mathbb{P}_N(Y=k) \) are the same for all \(k\), then  all \(\mathbb{P}_N^{{k}}\) are \((\max_{k}\Phi_{{k}},\delta) \)-concentrated.
\end{itemize}

We shall show that Exact CD implies Standard CD. By \(\operatorname{(i)}\), there exists a coverage \(\{\Omega_k\}_{k=1}^m\) such that \( \sum_{k=1}^{m}\mathbb{P}_N^{k}(\Omega_k) \times \mathbb{P}_N(Y=k) \geq 1-\delta \). Since \(\sum_{k=1}^{m} \mathbb{P}_N(Y=k) = 1\), there must exist at least one class \(\bar{k}\) such that \(\mathbb{P}_N^{\bar{k}}(\Omega_k)\geq 1 - \delta \). By \(\operatorname{(ii)}\), we have \(\Omega_{\bar{k}}^{+\epsilon} \subseteq \mathcal{D}_{f_{\theta},\bar{k}} \) and then \(\operatorname{Vol} (\Omega_{\bar{k}}^{+\epsilon}) \subseteq \operatorname{Vol} (\mathcal{D}_{f_{\theta},\bar{k}} )\). Finally, by applying the Burnn-Monkowski inequality \citep{gardner2002brunn} (for \(d_{\mathcal{X}}\) being \(l_2\) or \(l_{\infty}\))/concentration theorem \citep{talagrand1995concentration} (for \(d_{\mathcal{X}}\) being \(l_1\))/\citep{mcdiarmid1989method} (for \(d_{\mathcal{X}}\) being  \(l_{0}\)), we can shrink down the budget \(+\epsilon\) and obtain an upper bound \(\Phi_{\bar{k}}\) of \(\operatorname{Vol} (\Omega_{\bar{k}})\).

\begin{remark}
    In the original results, \citet{pal2023adversarial,pal2024certified} stated the Strong-CD and Standard-CD for any empirical distribution. By focusing  on the finite empirical \(\mathbb{P}_N\), the necessary and sufficient condition of adversarial robustness can be exactly characterized by our proposed Exact CD condition, noting that searching for \(\Omega_k\) is NP-hard. This translates the \(f_{\theta}\)-independent sufficient condition \(\operatorname{(i')}+\operatorname{(i'')}\) and  volume-based necessary condition \(\operatorname{(i^*)}\) into a countable geometric property directly tied to the training dataset \(\operatorname{(i)}+\operatorname{(ii)}\). In addition, the proposed framework allow us to extend the above Definition~\ref{def:robust-classifier} of (point-wise) $(\epsilon, \delta)$-robust classifier to its distributional counterpart by requiring 
    \begin{equation}\label{eq:DRC}
        \sup_{\mathbb{P} \colon \mathcal{W}_{p}(\mathbb{P},\mathbb{P}_N) \leq \epsilon } \mathbb{P} \left(f_{\theta}(X)\ne Y  \right) \leq \delta.
    \end{equation}
    As a consequence of Theorem~\ref{thm:main0}, a sufficient condition of \eqref{eq:DRC} is \(\operatorname{cc}_p(\epsilon)+\hat{\mathcal{R}}\leq\delta \); and a necessary condition of \eqref{eq:DRC} is  \(\operatorname{lb}_p(\epsilon)+\hat{\mathcal{R}}\leq\delta \). This successfully connects this notion of a robust classifier with the popular concept of Wasserstein distributionally robust risk.
\end{remark}

\subsection{Extension for Deep Neural Networks}\label{sec:AS}

Despite the theoretical tightness of the concave certificates discussed in previous sections, calculating the rates exactly for a complex network is often intractable. To bridge the gap between these powerful theoretical bounds and the practical requirements of deep learning, we dedicate this entire section to focus on the Euclidean assumption and derive the corresponding results explicitly for modern architectures, including LayerNorm and Attention maps. 

\begin{assumption} \label{asmp:Zdl}
The data space is given as
    \(\mathcal{Z}=\mathcal{X}\times\mathcal{Y}\), where  \(\mathcal{X}\subseteq\mathbb{R}^{n} \) is the space of features and \(\mathcal{Y}\subseteq\mathbb{R}^{m}\) is the space of labels. The cost function \(d\) is  given by 
    \[d(z',z) = \norm{x'-x}_{r} + \kappa\norm{y'-y}_{1},\]
    where \(\norm{\cdot}_{r}\) is \(r\)-norm defined on \(\mathbb{R}^{n}\) with \(r\in[1,\infty]\) and \(\kappa\in(0,\infty) \cup \{\infty\}\).
\end{assumption}

To compensate vector-to-vector maps in  deep neural networks, we first extend the notion of rate function in Definition~\ref{def:rate} for any \(f\colon\mathcal{X}\subseteq\mathbb{R}^{n}\rightarrow\mathbb{R}^{\tilde{n}} \) where \(\tilde{n}>1\) as 
\begin{equation*}
    \Delta_{f}(x,t) \triangleq \sup_{x'\in\mathcal{X}}\left\{\norm{f(x') - f(x)}_r \colon \norm{x'-x}_r \leq t \right\}. 
\end{equation*}
In fact, this is precisely the  modulus of continuity \citep{timan1963theory} of \(f\). We now define the \textit{adversarial score} as a relaxation of the \(\operatorname{cc}_1\) when exact least concave majorant \(\mathcal{C}\) is not available.
\begin{definition}[adversarial score] \label{def:adv-score}
    Given Assumption~\ref{asmp:Zdl} and \(f\colon\mathcal{X}\subseteq\mathbb{R}^{n}\rightarrow\mathbb{R}^{\tilde{n}} \), then \(F\) is called an adversarial score  of \(f\) if \(F\) is  \textbf{non-decreasing concave} and \(\sup_{x\in\mathcal{X}}\Delta_{f}(x,t) \leq F\left(t \right)\) for any \(t\geq0\). By Theorem~\ref{thm:main0}, if \(\mathcal{A}_{\theta}\) is an adversarial score of \(\bm{l}(\cdot;\theta)\), then for any \(p\in[1,\infty]\),
    \begin{equation*}
        \mathcal{R}_{p}(\epsilon)-\hat{R} \leq  \operatorname{cc}_{1}(\epsilon) \leq \mathcal{A}_{\theta}(\epsilon).
    \end{equation*}
\end{definition}
As demonstrated in the later part of this section, $\mathcal{A}_{\theta}$ can be derived explicitly for both classification (Proposition~\ref{prop:classification}) and regression (Proposition~\ref{prop:regression}).

\subsubsection{Hypothesis Function} 
A  deep neural network is a hypothesis function \(f_{\theta}\colon\mathcal{X}\rightarrow\mathcal{Y}\) parameterized by its weights \(\theta\in\Theta \), where \( f_{\theta} = f_{\theta}^{(K)}\circ f_{\theta}^{(K-1)} \circ \dots \circ f_{\theta}^{(1)} \) is  a composition of \( K \) component functions (or layers) \( f^{(k)} \). The following Lemma~\ref{lem:layer-rule} shows that analysis of the entire network \(f_{\theta}\) can be reduced to each individual layer. For example, the  LayerNorm \(x\mapsto \frac{x - \operatorname{Mean}(x) }{\sqrt{\operatorname{Var}(x)+c}}\cdot w + b \) can be seen as \( {\operatorname{Linear}}\circ {\operatorname{Normalization}} \circ {\operatorname{Centering}} \) or the Attention \(X \mapsto  \operatorname{Softmax}(XWX^T) XV \) can be seen as \(\operatorname{A-map}\times\operatorname{Linear} \).
\begin{lemma}[Layer Rule]\label{lem:layer-rule}
    Suppose that  \(\overset{\circ}{f}(x) = g(h(x))\) and \(\overset{\times}{f}(x) = g(x) \times h(x)\). Then
    \begin{itemize}
        \item[\(\operatorname{(a)}\)] \( \sup_{x\in\mathcal{X}} \Delta_{\overset{\circ}{f}}(x,t) \leq \sup_{x\in\mathcal{X}} \Delta_{g}(h(x),\Delta_{h}(x,t)) \). \hfill (composition)
        \item[\(\operatorname{(b)}\)] If \(\sup_{x\in\mathcal{X}}\norm{g(x)}_r \leq M_g \) and \(\sup_{x\in\mathcal{X}}\norm{h(x)}_r \leq M_h \) then \( \sup_{x\in\mathcal{X}} \Delta_{\overset{\times}{f}}(x,t) \leq M_h\sup_{x\in\mathcal{X}}  \Delta_{g}(x,t) +  M_g\sup_{x\in\mathcal{X}}  \Delta_{h}(x,t) \). \hfill (product)
        \item[\(\operatorname{(c)}\)] Composition and Product rules still hold when replacing rate \(\Delta\) with adversarial score \(F\).
    \end{itemize}
\end{lemma}
By Lemma~\ref{lem:layer-rule}, it is sufficient to calculate adversarial scores of each layer function. It is worth noting that in the following Example~\ref{exam:activation}, many layers have adversarial score \(F(t)\) strictly lower than their Lipschitz certificates \(\operatorname{Lip}_ft\).

\begin{figure}[h]
    \centering
    \includegraphics[width=0.9\linewidth]{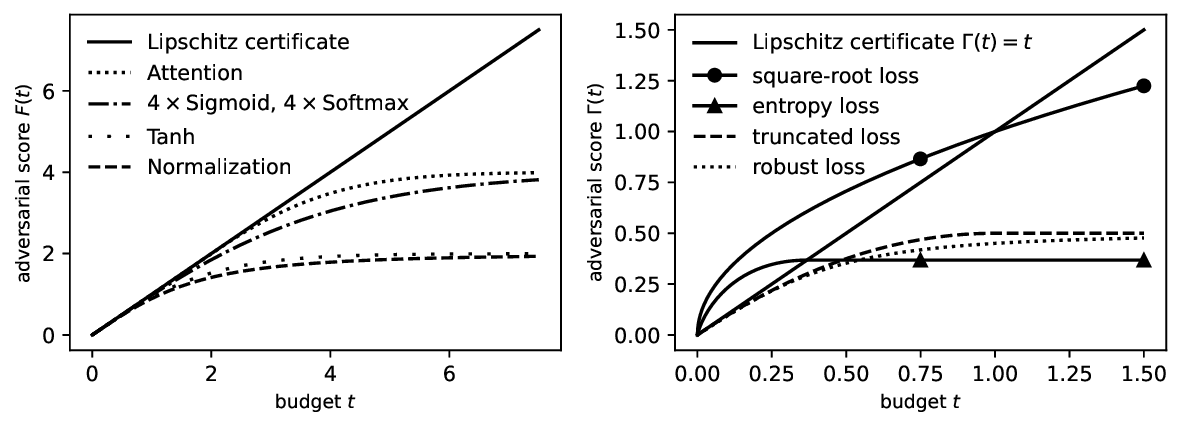}
    \caption{Adversarial scores for various layer  (Example~\ref{exam:activation} - Left) and loss  (Example~\ref{exam:gamma} - Right) functions. Notably, LayerNorm exhibits robustness similar to Tanh, while Attention is similar to Sigmoid, despite being scaled to the same Lipschitz modulus. By providing a non-linear robustness analysis, our framework also allows us to study non-Lipschitz, non-differentiable losses (square-root, entropy) where traditional methods fail.  }
    \label{fig:adv-score}
\end{figure}

\begin{example}[Figure~\ref{fig:adv-score} - Left] \label{exam:activation}
    The adversarial score \(F\) of some common layer functions are given as follows.
\begin{itemize}
    \item[\(\operatorname{(a)}\)] Saturating activation (Sigmoid, Tanh) where \(f\colon x \mapsto (\sigma(x_1),\dots,\sigma(x_n)) \). If \(r=1\) then \(F(t) = n\left[\sigma\left(\frac{t}{2n}\right) - \sigma\left(\frac{-t}{2n}\right)\right]\); if \(r=2\) then \(F(t) = \sqrt{n}\left[\sigma\left(\frac{t}{2\sqrt{n}}\right) - \sigma\left(\frac{-t}{2\sqrt{n}}\right)\right]\); if \(r=\infty\) then \(F(t) = \sigma\left(\frac{t}{2}\right) - \sigma\left(\frac{-t}{2}\right) \). In all cases, \(F(t)< \operatorname{Lip}_{f} t\).
    \item[\(\operatorname{(b)}\)] Softmax \(f(x) = \frac{1}{\inprod{\exp(x),\bm{1}_n}}\exp(x) \): then \(F_{\operatorname{Softmax}}(t) = F_{\operatorname{Sigmoid}}(t) < \operatorname{Lip}_{f} t\).
    \item[\(\operatorname{(c)}\)]  Normalization \(f(x) = \frac{x }{\sqrt{\norm{x}_2^2/n+c}} \) when \(r=2\): then \(F(t) = \frac{t}{\sqrt{t^2/4n + c}} < \operatorname{Lip}_{f} t = t/\sqrt{c} \).     
    \item[\(\operatorname{(d)}\)] A-map with bounded input \(f(X) = \operatorname{Softmax}(XWX^T) \) where \(X\in\mathbb{R}^{n\times \tilde{n}}\), \(r=2\) and \(\sup_{X\in\mathcal{X}}\norm{X}_2\leq c \): then \(F_{\operatorname{A-map}}(t) = \mathcal{C}(F_{\operatorname{Softmax}}(2c\norm{W}_2t + \norm{W}_2t^2))< \operatorname{Lip}_{f} t \).
    \item[\(\operatorname{(e)}\)] Other Lipschitz activations (\(\operatorname{ReLU}\), \(\operatorname{Softplus}\)); Cross-Entropy loss \((\log \circ \operatorname{Softmax})\); Margin loss \cite{carlini2017towards};  Centering \(f(x) = x - \operatorname{Mean}(x)\);  Linear layer \(f(x) = Wx + b\): then  \(F(t) = \operatorname{Lip}_{f} t \).
\end{itemize}
\end{example}

\subsubsection{Classification}\label{sec:classification}
Consider an \(m\)-classification problem with feature space \( \mathcal{X}\subseteq\mathbb{R}^{n} \) and label space  \( \mathcal{Y}\subseteq\Delta^{m}\coloneqq\left\{y\in\mathbb{R}^{m}\mid \sum_{j=1}^{m}y_j=1, y\geq0 \right\} \). To fit a network \( f_{\theta}\colon\mathcal{X}\rightarrow\mathbb{R}^{m} \) that predicts \(y\) from the state \( f_{\theta}(x) \), one often considers the loss function given by

\begin{equation}\label{eq:loss-classification}
	\bm{l}(x,y;\theta) =  \inprod{y,f_{\theta}(x)}.
\end{equation}
We shall show that the adversarial score of \(\bm{l}\) can be calculated directly from the adversarial score  \(F_{\theta}\) of the network \(f_{\theta}\) as studied in Lemma~\ref{lem:layer-rule}, see proof in Appendix~\ref{proof:classification}.

\begin{proposition}[Adversarial Score in Classification]
\label{prop:classification}
   Given Assumption~\ref{asmp:Zdl} where \(d(z',z) = \norm{x'-x}_{r} + \kappa\norm{y'-y}_{1}\) and a network \(f_{\theta}\colon\mathcal{X}\rightarrow\mathbb{R}^{m}\), define the loss function by \(\bm{l}(x,y;\theta) =  \inprod{y,f_{\theta}(x)} \). Suppose that  \(F_{\theta} \) is an adversarial score of \(f_{\theta}\) (Definition~\ref{def:adv-score}). 
   \begin{enumerate}
       \item[$\operatorname{(a)}$] If \(\kappa=\infty \), then  \(\mathcal{A}_{\theta}(t)\coloneqq F_{\theta}(t) \) is an adversarial score of \(\bm{l}\).
       \item[$\operatorname{(b)}$] If \(\kappa\in(0,\infty)\) and there exists   \(M\in(0,\infty)\) such that \(\norm{f_{\theta}(x)}_{\infty}\leq M \) for any \(x\in\mathcal{X}\), then  \(\mathcal{A}_{\theta}(t)\coloneqq\sup_{\tau\in[0,t]} \left\{  F_{\theta}(t-\tau) + M\kappa^{-1}\tau \right\}\) is an adversarial score of \(\bm{l}\). 
   \end{enumerate}
\end{proposition}

Proposition~\ref{prop:classification} justifies that the classification model with loss \( \bm{l}(x,y;\theta) = \inprod{y,f_{\theta}(x)} \) is always robust with respect to its feature \( x \), and robust with respect to its label if the output \( f_{\theta}(x) \) is bounded by \( M \). This covers several practical scenarios such as when \(f_{\theta}\) is continuous and \(\mathcal{X}\) is compact (images' pixels or waves' signals), then \(M<\infty\); or when last layer is the Softmax layer then, \(M=1\).

\subsubsection{Regression}\label{sec:regression}

Consider the regression problem with  feature space \( \mathcal{X}\subseteq\mathbb{R}^{n} \) and label space  \( \mathcal{Y}\subseteq\mathbb{R} \). To fit a network \( f_{\theta}\colon\mathcal{X}\rightarrow\mathbb{R}^{m} \) to predict \( y \), one often considers the loss function given by
\begin{equation}\label{eq:loss-regression}
	\bm{l}(x,y;\theta) =  \gamma\left( \abs{y - f_{\theta}(x)} \right),
\end{equation}
where \( \gamma\colon[0,\infty)\rightarrow[0,\infty) \). Example~\ref{exam:gamma} lists common regression loss functions (possibly non-differentiable, non-convex, or non-Lipschitz). One can observe that \( \Gamma(t) < \operatorname{Lip}_{\gamma} t \) in many scenarios, providing evidence of our competitiveness and versatility  compared to traditional certificates. 

\begin{example}
    [Figure~\ref{fig:adv-score} - Right] \label{exam:gamma} The adversarial score \(\Gamma\) of some common regression losses \(\gamma\colon[0,\infty)\rightarrow[0,\infty)\) are given as follows.
    \begin{enumerate}
        \item[\(\operatorname{(a)}\)]  For Holder loss \(\gamma(t) =  c t^{\alpha} \), if \(\alpha \in(0,1)\) then \(\Gamma(t) = c t^{\alpha}<\operatorname{Lip}_{\gamma}t=\infty\); if \(\alpha =1\) then \(\Gamma(t) =\operatorname{Lip}_{\gamma} t=t\); if \(\alpha\in(1,\infty)\) then \(\Gamma(t)=\operatorname{Lip}_{\gamma}t=\infty\).
        \item[\(\operatorname{(b)}\)] If \(\gamma\) is the Huber loss \( \left(\gamma(t)=  \frac{t^2}{2} \text{ if }t\in[0,c],\, c{t} - \frac{c^2}{2} \text{ otherwise} \right)\)  then \(\Gamma(t) =\operatorname{Lip}_{\gamma} t=ct\).
        \item[\(\operatorname{(c)}\)] If \(\gamma\) is the truncated loss \citep{yang2014robust} \(\left( \gamma(t)= \min \left\{\frac{1}{2}c^2 , \frac{1}{2}t^{2}\right\}\right)\) then \( \Gamma(t) = \frac{2tc - t^2}{2} \text{  if } t\in[0,c],\, \frac{c^2}{2}  \text{ if }t\in(c,\infty)\) and \(\Gamma(t)< \operatorname{Lip}_{\gamma} t=ct\).
        \item[\(\operatorname{(d)}\)] If \(\gamma\) is the robust loss \citep{barron2019general} \(\left(\gamma(t)= \frac{1}{2}\frac{c^2 t^2}{ac^2+t^2}\text{ where } a = \frac{27}{256}\right)\) then \(\Gamma(t) = \gamma(s_t+t)-\gamma(s_t) \) where \(s_t = \frac{1}{6}\left( \sqrt{3t^2 + 6\sqrt{t^4+4ac^2t^2+16a^2c^2}-4ac^2  }-3t \right) \). Moreover, \(\Gamma(t)< \operatorname{Lip}_{\gamma} t=ct\).
        \item[\(\operatorname{(e)}\)] If \(\gamma\) is the entropy-like loss \(\left(\gamma(t)=  -t\log(t) \text{ if }t\in[0,e^{-1}],\, e^{-1} \text{ otherwise}\right)\), then \(\Gamma(t) = \gamma(t) < \operatorname{Lip}_{\gamma}t = \infty\).
    \end{enumerate}
\end{example}

We show that the adversarial score of the regression model \eqref{eq:loss-regression} now can be calculated from the adversarial score  \(F_{\theta}\) of the network \(f_{\theta}\) and  the  adversarial score of the given \(\gamma\). 

\begin{proposition}[Adversarial Score in Regression]
\label{prop:regression}
    Given Assumption~\ref{asmp:Zdl} where \(d(z',z) = \norm{x'-x}_{r} + \kappa\norm{y'-y}_{1}\) and a network \(f_{\theta}\colon\mathcal{X}\rightarrow\mathbb{R}^{m}\), define the loss function by  \(\bm{l}(x,y;\theta) = \gamma\left( \abs{y - f_{\theta}(x)} \right) \).  Suppose that  \(F_{\theta} \) and \(\Gamma\) are an adversarial scores of \(f_{\theta}\) and \(\gamma\), respectively. 
    \begin{enumerate}
       \item[$\operatorname{(a)}$] If \(\kappa=\infty \), then  \(\mathcal{A}_{\theta}(t)\coloneqq \Gamma(F_{\theta}(t)) \) is an adversarial score of \(\bm{l}\).
       \item[$\operatorname{(b)}$] If \(\kappa\in(0,\infty)\) then \(\mathcal{A}_{\theta}(t) = \Gamma \left(\sup_{\tau\in[0,t]} \left\{  F_{\theta}(t-\tau) + \kappa^{-1}\tau \right\} \right)\)  is an adversarial score of \(\bm{l}\). 
   \end{enumerate}
\end{proposition}

Similar to Proposition~\ref{prop:classification},  Proposition~\ref{prop:regression} also reveals that  the label sensitivity \(\kappa\) in the cost function \(d(z',z) = \norm{x'-x}_{r} + \kappa\norm{y'-y}_{1}\) dictates the robustness of the entire network by balancing the impact of feature noise against label shifts. Specifically, a finite $ \kappa $ effectively couples these two sources of uncertainty, whereas $ \kappa = \infty $ simplifies the certificate to a focus on feature stability only.

\subsection{Generalization Bound under Wasserstein Shift} \label{sec:GB}
Having established tractable robustness certificates for deep networks, we now turn to statistical learning theory to understand how these models generalize to unseen data. The following corollary is an immediate consequence of Theorem~\ref{thm:main0}.
\begin{corollary}
    [The Worst-case Generalization Bound via Concave Complexity]\label{cor:det-gen-bound} Suppose that  \(\mathcal{W}_{p}(\mathbb{P}_{\rm true}, \mathbb{P}_N) \leq\epsilon \) for some \(p\in[1,\infty]\) and denote the worst-case generalization bound as \(\operatorname{GB}_{p}(\epsilon)\triangleq \sup_{\theta\in\Theta} \mathbb{E}_{\mathbb{P}_{\rm true}}[\bm{l}(Z;\theta)] - \mathbb{E}_{\mathbb{P}_{N}}[\bm{l}(Z;\theta)].\)    Define the concave complexity of the loss class \(\mathcal{L} = \{\bm{l}(\cdot;\theta) \mid \theta\in\Theta \} \) as
    \begin{equation}\label{eq:CC}
        \operatorname{(CC)} \quad \hat{\mathfrak{C}}_{\mathcal{Z}_N}(\mathcal{L},\epsilon) \triangleq \sup_{\theta\in\Theta} \operatorname{cc}_1(\epsilon) = \sup_{\theta\in\Theta} \mathcal{C}_{\Delta_{\theta}^{\max}}(\epsilon).
    \end{equation}
    where \(\operatorname{cc}_1 \) is given in Theorem~\ref{thm:main0}. Then
    \begin{equation}\label{eq:det-gen-gap}
         \operatorname{GB}_{p}(\epsilon) \leq \hat{\mathfrak{C}}_{\mathcal{Z}_N}(\mathcal{L},\epsilon). 
    \end{equation}
    Notably, under the setting of Example~\ref{exam:unbounded}, equality holds for \(p=\alpha=1\). If \(cc_1(\epsilon)\) is not available, we can relax it by any valid adversarial score \(\sup_{\theta\in\Theta} \operatorname{cc}_1(\epsilon) \leq\sup_{\theta\in\Theta} \mathcal{A}_{\theta}(\epsilon)  \) (Definition~\ref{def:adv-score}).
\end{corollary}

\begin{figure}[h!]
    \centering
    \includegraphics[width=0.95\linewidth]{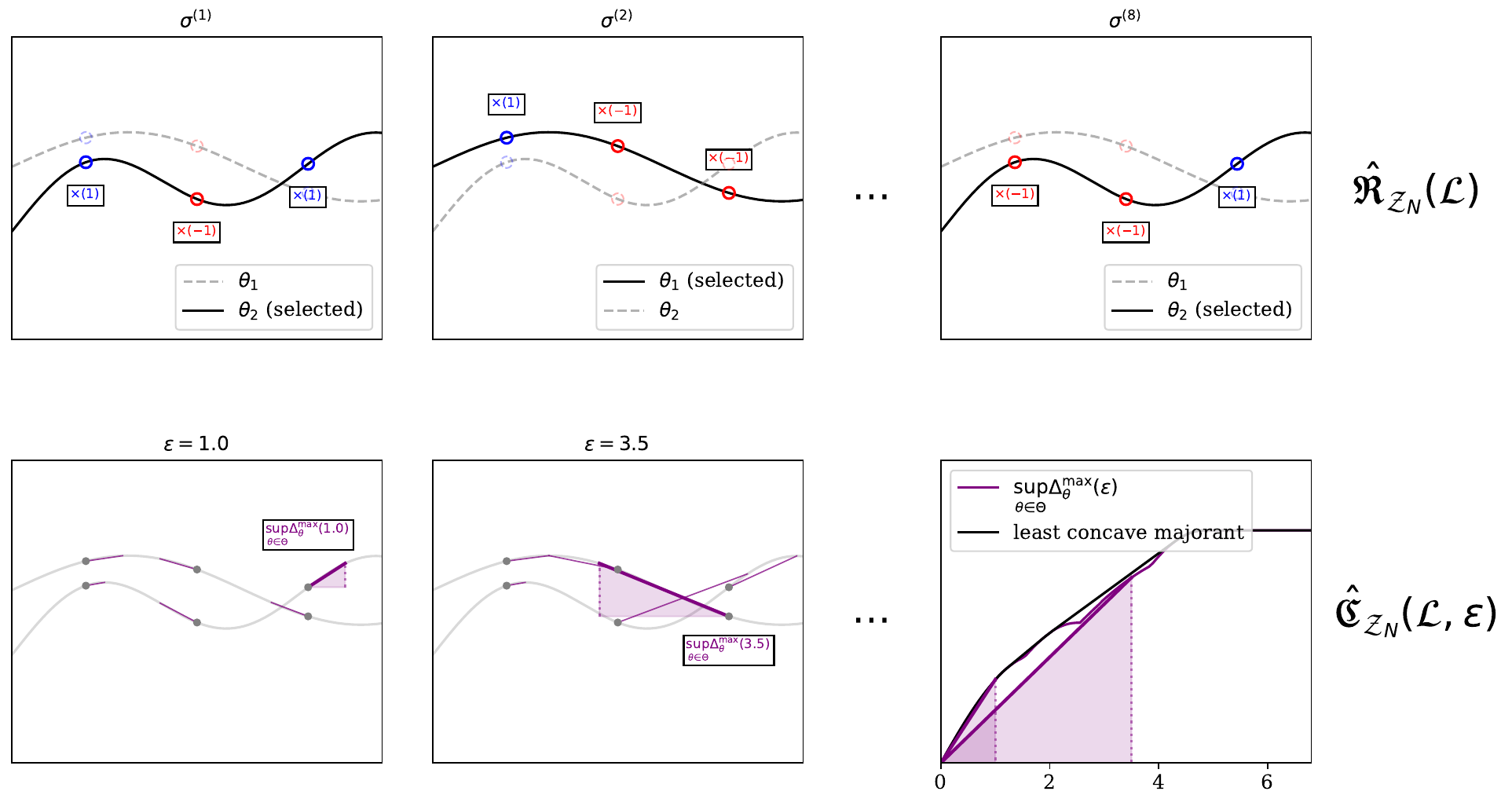}
    \caption{Intuition of  \(\hat{\mathfrak{R}}_{\mathcal{Z}_N}(\mathcal{L})=\frac{1}{8}\sum_{k=1}^{8}\sup_{\theta \in \Theta} \frac{1}{3} \sum_{i=1}^{3} \sigma_i^{(k)} \boldsymbol{l}(Z^{(i)}; \theta)\) and \(\hat{\mathfrak{C}}_{\mathcal{Z}_N}(\mathcal{L},\epsilon)=\mathcal{C}_{\sup_{\theta \in \Theta} \Delta_{\theta}^{\max}}(\epsilon)\). CC calculates in average how well \(\mathcal{L}\) could fit random noise. CC calculates at worst how fluctuating \(\mathcal{L}\) could be.}
    \label{fig:RC-CC}
\end{figure}

Recall that the empirical  Rademacher complexity (RC)  \eqref{eq:RC}  measures \underline{on average} how well a loss class fits random noise, yielding a statistical generalization bound \eqref{eq:stat-gen-bound} driven by sample size \(\operatorname{GB} \leq \operatorname{const}\times\hat{\mathfrak{R}}_{\mathcal{Z}_N}(\mathcal{L}) + \operatorname{conf}(\delta)\). In contrast, our proposed concave complexity (CC)  \(\hat{\mathfrak{C}}_{\mathcal{Z}_N}(\mathcal{L},\epsilon)\) \eqref{eq:CC} measures \underline{the worst} loss fluctuation (Figure~\ref{fig:RC-CC}), yielding the worst-case generalization bound governed by the transport budget \(\epsilon\) rather than \(N\). Based on standard Wasserstein concentration inequalities \cite{fournier2015rate,weed2019sharp}, the required budget in high-dimensional settings ($n > 2p$) scales as $\epsilon = \mathcal{O}(N^{-1/n})$. This reveals a fundamental theoretical trade-off: \eqref{eq:CC} eliminates structural dependencies in \eqref{eq:RC} but replaces standard RC rate $\mathcal{O}(1/\sqrt{N})$  with the $\mathcal{O}(N^{-1/n})$  embedded within $\epsilon$. We shall show in the later part of this section that \eqref{eq:det-gen-gap} remains beneficial for analyzing many real-life over-parameterized networks.

\subsubsection{Properties of Concave Complexity}

Interestingly, the proposed concave complexity \(\hat{\mathfrak{C}}_{\mathcal{Z}_N}\) satisfies several properties analogous to Rademacher complexity \(\hat{\mathfrak{R}}_{\mathcal{Z}_N}\).

\begin{proposition}\label{prop:calculus}
    Given \(0<\epsilon<\epsilon'<\infty\), \(b,c>0\) and \(\mathcal{L}\subseteq\mathcal{L}'\), then 
    \begin{itemize}
        \item[$\operatorname{(a)}$]  \(\hat{\mathfrak{C}}_{\mathcal{Z}_N}(\mathcal{L},\epsilon) \leq\hat{\mathfrak{C}}_{\mathcal{Z}_N}(\mathcal{L},\epsilon') \) and  \(\hat{\mathfrak{C}}_{\mathcal{Z}_N}(\mathcal{L},\epsilon+\epsilon') \leq\hat{\mathfrak{C}}_{\mathcal{Z}_N}(\mathcal{L},\epsilon)+\hat{\mathfrak{C}}_{\mathcal{Z}_N}(\mathcal{L},\epsilon') \). \hfill (subadditivity)
        \item[$\operatorname{(b)}$]  \(\hat{\mathfrak{C}}_{\mathcal{Z}_N}(\mathcal{L},\epsilon) \leq\hat{\mathfrak{C}}_{\mathcal{Z}_N}(\mathcal{L}',\epsilon) \) and \(\hat{\mathfrak{C}}_{\mathcal{Z}_N}(c\cdot\mathcal{L}+b,\epsilon) = c\hat{\mathfrak{C}}_{\mathcal{Z}_N}(\mathcal{L},\epsilon) \). \hfill (positively affine)
        \item[$\operatorname{(c)}$]  \(\hat{\mathfrak{C}}_{\mathcal{Z}_N}(\operatorname{conv}(\mathcal{L}),\epsilon) =\hat{\mathfrak{C}}_{\mathcal{Z}_N}(\mathcal{L},\epsilon) \). \hfill  (invariant under convex hull)
    \end{itemize}
\end{proposition}

\noindent\textit{Proof.} (a) For any \(\theta\), one has that \(\mathcal{C}_{\Delta_{\theta}^{\max}}\) is concave, non-decreasing (by Lemma~\ref{lem:univariate-majorant}) and \(\mathcal{C}_{\Delta_{\theta}^{\max}}(t)\geq \mathcal{C}_{f^{\max}}(0)\geq \Delta_{\theta}^{\max}(0)\geq 0\). Thus, \(\mathcal{C}_{\Delta_{\theta}^{\max}}\) is subadditive. Since supremum is also subadditive, so is \(\hat{\mathfrak{C}}_{\mathcal{Z}_N}\). (b) As the maximal rate of \(c\cdot \bm{l} +b\) is \(c\cdot\Delta_{\theta}^{\max}  \) and   \( \mathcal{C}_{c\cdot\Delta_{\theta}^{\max} } =c\mathcal{C}_{\Delta_{\theta}^{\max}}\).\\
(c) Let \(\bar{\bm{l}}= \sum_{j} \alpha_j \bm{l}_j \in \operatorname{conv}(\mathcal{L})\) where \(\sum \alpha_j = 1, \alpha_j \ge 0\). Then the individual rate function of \(\bar{\bm{l}}\) is given by \(\Delta_{\bar{\bm{l}}}(\hat{z}, t) = \sup\limits_{z'\in\mathcal{Z}}  \sum\limits_{j} \alpha_j (\bm{l}_j(z')-\bm{l}_j(\hat{z}))\). Since supremum is subadditive, \(\Delta_{\bar{\bm{l}}}(\hat{z}, t) \leq \sum\limits_{j} \alpha_j \Delta_{\bm{l}_j}(\hat{z}, t)\). 
By Lemma~\ref{lem:calculation}, for any \(\hat{z}\in\mathcal{Z}_N\), \(\mathcal{C}_{\Delta_{\bar{\bm{l}}}(\hat{z}, \cdot) } \leq \mathcal{C}_{\sum\limits_{j} \alpha_j \Delta_{\bm{l}_j}(\hat{z}, \cdot)} \leq \sum\limits_{j} \alpha_j \mathcal{C}_{\Delta_{\bm{l}_j}(\hat{z}, \cdot)} \leq \mathcal{C}_{\Delta_{\bm{l}_j}^{\max}(\cdot)}\).
Therefore, \(\mathcal{C}_{\Delta_{\bar{\bm{l}}}^{\max} }\leq \mathcal{C}_{\Delta_{\bm{l}_j}^{\max}}\). Take the supremum over all \(\bar{\bm{l}}\in\operatorname{conv}(\mathcal{L}) \), we have \(\hat{\mathfrak{C}}_{\mathcal{Z}_N}(\operatorname{conv}(\mathcal{L}),\epsilon) \leq \hat{\mathfrak{C}}_{\mathcal{Z}_N}(\mathcal{L},\epsilon)\). By part (b), one also has \(\hat{\mathfrak{C}}_{\mathcal{Z}_N}(\operatorname{conv}(\mathcal{L}),\epsilon) \geq\hat{\mathfrak{C}}_{\mathcal{Z}_N}(\mathcal{L},\epsilon)\).
\hfill\Halmos

To this end, we illustrate a contraction lemma, which mirrors the Ledoux-Talagrand Lemma \citep{ledoux2013probability}, allowing for the analysis of composite functions.

\begin{proposition}[Contraction Lemma] \label{prop:contraction} Let \(\mathcal{L}\coloneqq \{\bm{l}_{\theta}= \ell\circ f_{\theta} \mid f_{\theta} \in \mathcal{F} \} \) be a class of loss functions where \(-f_{\theta}\in\mathcal{F} \) for any \(\theta\) and  \(\ell \) is a \(\operatorname{Lip}_{\ell}\)-Lipschitz univariate function. Then
    \begin{equation*}
       \hat{\mathfrak{C}}_{\mathcal{Z}_N}(\mathcal{L},\epsilon) \leq \operatorname{Lip}_{\ell}\times \mathcal{C}_{ \hat{\mathfrak{C}}_{\mathcal{Z}_N}(\mathcal{F},\cdot)}(\epsilon).
    \end{equation*}
\end{proposition}

\noindent\textit{Proof.} If \(f_{\theta}(z')\geq f_{\theta}(\hat{z})\) then \(\bm{l}(z';\theta) - \bm{l}(\hat{z};\theta) =\ell \circ f_{\theta}(z')-\ell \circ f_{\theta}(\hat{z}) \leq \operatorname{Lip}_{\ell}\times (f_{\theta}(z')-f_{\theta}(\hat{z})) \leq  \operatorname{Lip}_{\ell}\times \Delta_{f_{\theta}}(\hat{z},t) \). Otherwise, \(\bm{l}(z';\theta) - \bm{l}(\hat{z};\theta) \leq  \operatorname{Lip}_{\ell}\times \Delta_{-f_{\theta}}(\hat{z},t) \). Therefore 
\begin{equation*}
    \Delta_{\bm{l}_{\theta}}(\hat{z},t) \leq  \operatorname{Lip}_{\ell}\times \max\{\Delta_{f_{\theta}}(\hat{z},t),\Delta_{-f_{\theta}}(\hat{z},t)\}.
\end{equation*}
Taking the maximum over all (finite) samples \(\hat{z}\in\mathcal{Z}_N\),
\begin{equation}\label{eq:contraction-1}
    \Delta_{\bm{l}_{\theta}}^{\max}(t) \leq  \operatorname{Lip}_{\ell}\times \max\{\Delta_{f_{\theta}}^{\max}(t),\Delta_{-f_{\theta}}^{\max}(t)\}\leq  \operatorname{Lip}_{\ell}\times \max\{\mathcal{C}_{\Delta_{f_{\theta}}^{\max}}(t),\mathcal{C}_{\Delta_{-f_{\theta}}^{\max}}(t)\}.
\end{equation}
Moreover, 
\begin{equation}\label{eq:contraction-2}
    \max\{\mathcal{C}_{\Delta_{f_{\theta}}^{\max}}(t),\mathcal{C}_{\Delta_{-f_{\theta}}^{\max}}(t)\} \leq \sup_{f_{\theta}\in\mathcal{F}} \mathcal{C}_{\Delta_{f_{\theta}}^{\max}}(t) = \hat{\mathfrak{C}}_{\mathcal{Z}_N}(\mathcal{F},t).
\end{equation}
From \eqref{eq:contraction-1} and \eqref{eq:contraction-2}, one has \(\mathcal{C}_{\Delta_{\bm{l}_{\theta}}^{\max}} \leq \operatorname{Lip}_{\ell}\times\mathcal{C}_{ \hat{\mathfrak{C}}_{\mathcal{Z}_N}(\mathcal{F},\cdot)}\) 
The proof is completed by taking supremum over all \(\theta\). \hfill\Halmos

\begin{example} \label{exam:concave-complexity}
    Consider linear hypothesis \(f_{\theta}(x)= \inprod{\theta,x}\) and the cost function  \( d(z',z) = \norm{x'-x} + \infty\abs{y'-y} \). Denote \(\mathcal{F}\coloneqq \{ f_{\theta} \mid \theta\in\Theta\} \). By Example~\ref{exam:unbounded}, \(\operatorname{cc}_1(\epsilon) = \norm{\theta}_{*} \epsilon\) and thus  \(\hat{\mathfrak{C}}_{\mathcal{Z}_N}(\mathcal{F},\epsilon)  \leq   \sup_{\theta\in\Theta}\norm{\theta}_{*}\times \epsilon\).    Let  \(\mathcal{L}\coloneqq \{ \bm{l} = \ell \circ f_{\theta} \mid \theta\in\Theta\} \), if \(\ell\) is \(\operatorname{Lip}_{\ell}\)-Lipschitz  then
    \begin{equation*}
       \hat{\mathfrak{C}}_{\mathcal{Z}_N}(\mathcal{L},\epsilon)  \leq \operatorname{Lip}_{\ell}\times  \sup_{\theta\in\Theta}\norm{\theta}_{*}\times \epsilon.
    \end{equation*}
    In general, if \(\mathcal{A}_{\theta}\) is an adversarial score (Definition~\ref{def:adv-score}) of \(\bm{l}(\cdot;\theta) = \ell\circ  f_{\theta}\), then 
    \begin{equation*}
       \hat{\mathfrak{C}}_{\mathcal{Z}_N}(\mathcal{L},\epsilon)  \leq \mathcal{C}_{\sup_{\theta\in\Theta}\mathcal{A}_{\theta}}(\epsilon).
    \end{equation*}
\end{example}

For linear models, the proposed concave complexity $\hat{\mathfrak{C}}_{\mathcal{Z}_N}$  decouples the complexity of $\mathcal{L}$ into three factors: the Lipschitz constant $\operatorname{Lip}_{\ell}$ of the loss function, the maximal weight norm $\sup_{\theta\in\Theta}\|\theta\|_{*}$ of the hypothesis class, and the uncertainty radius $\epsilon$.   Compared to the standard empirical Rademacher complexity bounds $\hat{\mathfrak{R}}_{\mathcal{Z}_N} \leq \operatorname{Lip}_{\ell}  \sup_{\theta\in\Theta}\norm{\theta}_{*}/\sqrt{N} $ for  linear models \cite{bartlett2002rademacher,koltchinskii2002empirical,kakade2008complexity}, our geometric formulation  eliminates the restrictive assumptions such as bounded feature space $\mathcal{X}$ or a bounded/differentiable loss $\ell$. In general, because \(\mathcal{A}_{\theta}\) is calculated via the layer rule (Lemma~\ref{lem:layer-rule}),  $\hat{\mathfrak{C}}_{\mathcal{Z}_N}$ remains independent of the width \(h\) and depth \(K\) of the  architecture, highlighting the difference with Rademacher complexity  bounds  for neural networks \cite{bartlett2002rademacher,golowich2018size,neyshabur2018pac}. Standard RC bounds decay rapidly as \(1/\sqrt{N}\rightarrow0\), but suffer from network structural inflation.  For example, in GPT-3 model ($K=96, h=n=12288$) with $N=10^9$ samples, traditional RC bounds carry an uninformative term $\mathcal{O}(h K^{3/2}/\sqrt{N}) \approx 365$.  While we sacrifice the asymptotic speed of \(1/\sqrt{N}\), CC  absorbs $\epsilon \approx \mathcal{O}(N^{-1/n})$ and provides a strict geometric ceiling on how much the risk can degrade under any adversarial data shift. This highlights CC as a specialized tool for large architectures where \(N\) is finite.

\subsubsection{Adversarial Complexity Gaps}
We now extend our framework to compare the standard loss class \(\mathcal{L}\) with the class of worst-case losses \(\tilde{\mathcal{L}}_{\epsilon}=\{\tilde{\bm{l}}_{\epsilon}(\cdot;\theta) \mid \theta\in\Theta \} \) where 
\begin{equation*}
    \tilde{\bm{l}}_{\epsilon}(z;\theta)= \sup_{z' : d(z', z) \le \epsilon} \bm{l}(z'; \theta).
\end{equation*}
The gap between complexity of \(\tilde{\mathcal{L}}\) and \(\mathcal{L}\) measures the added complexity of adversarial learning. While prior work estimates this gap by exploiting specific structural properties of the loss \(\mathcal{L}\) (such as linearity or boundedness), we aim to bound it using the geometric certificates from Theorem~\ref{thm:main0}. 

\begin{theorem}[Adversarial Complexity Gaps]\label{thm:ACG}
     Given Notation~\ref{nota:main} where \(\mu_i=\frac{1}{N}\), denote the  class of individual rates \(\Delta_{\theta}\) \eqref{eq:ind-rate} as \({\Upsilon}_{\epsilon}\coloneqq \{z\mapsto \Delta_{\theta}(z,\epsilon) \mid \theta\in\Theta \} \),   then 
     \begin{equation}\label{eq:ARC-RC}
         \abs{\hat{\mathfrak{R}}_{\mathcal{Z}_{N}}(\tilde{\mathcal{L}}_{\epsilon})- \hat{\mathfrak{R}}_{\mathcal{Z}_{N}}(\mathcal{L})} \leq \hat{\mathfrak{R}}_{\mathcal{Z}_{N}}({\Upsilon}_{\epsilon}) \leq \sup_{\theta\in\Theta}\operatorname{lb}_{\infty}(\epsilon)= \sup_{\theta\in\Theta}\sum_{i=1}^{N}\frac{1}{N} \Delta_{\theta}(Z^{(i)},\epsilon),
     \end{equation}   
     where the first inequality is from prior works \cite{yin2019rademacher,awasthi2020adversarial}.   In addition, if \(\Delta_{\theta}(z,\epsilon) = \Delta_{\theta}^{\max}(\epsilon)\) for any \(z\in\mathcal{Z}_N \), then one can replace the last term with a tighter bound \(\frac{1}{\sqrt{N}} \sup_{\theta\in\Theta}\Delta_{\theta}^{\max}(\epsilon)\).     On the other hand, for any \(t>0\),  
     \begin{equation}\label{eq:ACC-CC}
        0 \leq \hat{\mathfrak{C}}_N(\tilde{\mathcal{{L}}}_{\epsilon},t) \leq \hat{\mathfrak{C}}_N(\mathcal{L},t) +\hat{\mathfrak{C}}_N({\Upsilon}_{\epsilon},t). 
    \end{equation}
\end{theorem}

We first observe that both ARC-RC and ACC-CC gaps are directly bounded by the complexity of \(\Upsilon_{\epsilon}\) (the class of individual rate functions). In addition, \eqref{eq:ARC-RC} reveals that ARC-RC gap somewhat constrains to the point-wise perturbation \(p=\infty\) via \(\operatorname{lb}_{\infty}\).  In contrast, the ACC-CC gap ties to \(\hat{\mathfrak{C}}_N\) \eqref{eq:CC} and \(\operatorname{cc}_1\) \eqref{eq:det-gen-gap}. Thus it covers all \(\mathcal{W}_p\)-distributional perturbations, up to \(p=1\).   In the remainder of this section, we illustrate Theorem~\ref{thm:ACG} for linear models (Example~\ref{exam:ACG-linear}), MLPs (Example~\ref{exam:ACG-MLP}), and compare our results with existing works.

\begin{example}[Adversarial complexity gap for linear models] \label{exam:ACG-linear}
    Under the setting of linear hypothesis in  Example~\ref{exam:concave-complexity}, we have \({\Upsilon}_{\epsilon} =  \{z\mapsto \Delta_{\theta}(z,\epsilon) \mid \theta\in\Theta \}  = \{ z\mapsto\epsilon\norm{\theta}_* \mid \theta\in\Theta \}\).     Thus
    \begin{equation}\label{eq:ARC-RC-linear}
        \abs{\hat{\mathfrak{R}}_{\mathcal{Z}_{N}}(\tilde{\mathcal{F}}_{\epsilon})- \hat{\mathfrak{R}}_{\mathcal{Z}_{N}}(\mathcal{F})} \leq \frac{ \epsilon }{\sqrt{N}}\sup_{\theta\in\Theta}\norm{\theta}_{*}.
    \end{equation}
    Specifically, if \(\Theta= \Theta_*\coloneqq \{\theta\in\mathbb{R}^n\mid \norm{\theta}_* \leq c  \}  \)  then \(\abs{\hat{\mathfrak{R}}_{\mathcal{Z}_{N}}(\tilde{\mathcal{F}}_{\epsilon})- \hat{\mathfrak{R}}_{\mathcal{Z}_{N}}(\mathcal{F})}\leq  \frac{ \epsilon }{\sqrt{N}}c  \). 
    On the other hand, \(\Delta_{\Delta_{\theta}(\cdot,\epsilon)} = 0 \) for any \(\theta\in\Theta\), hence \(\hat{\mathfrak{C}}_N(\Upsilon_{\epsilon},t)=0\) and for any \(t>0\),
    \begin{equation}\label{eq:ACC-CC-linear}
        \hat{\mathfrak{C}}_N(\tilde{\mathcal{{L}}}_{\epsilon},t) \leq  \hat{\mathfrak{C}}_N(\mathcal{L},t) = t\sup_{\theta\in\Theta}\norm{\theta}_{*}.
    \end{equation}
    
\end{example}

In Example~\ref{exam:ACG-linear}, \eqref{eq:ARC-RC-linear} says that the adversarial linear model \(\tilde{\mathcal{{L}}}_{\epsilon}\) \underline{can be} worse  than the empirical linear model \(\mathcal{L}\), proportionally with the adversarial budget \(\epsilon\) and the bound of learning weight \(c\). In contrast, \eqref{eq:ACC-CC-linear}  says that \(\tilde{\mathcal{{L}}}_{\epsilon}\) is \underline{never} worse  than \(\mathcal{L}\), which  reflects more accurately the reality that adversarial perturbations merely shift the entire loss class without inflating its complexity. We now revisit some works on  \eqref{eq:ARC-RC-linear} using different techniques.
\begin{itemize}
    \item[$\operatorname{(a)}$] In \citet[Thm 2]{yin2019rademacher}, authors consider   the \(\infty\)-adversarial attack (i.e., \(d(z',z) = \norm{x'-x}_{r=\infty} + \infty\abs{y'-y}\)) and show that if  \(\Theta= \{\theta\in\mathbb{R}^n\mid \norm{\theta}_{s} \leq  W   \}  \) then \(\hat{\mathfrak{R}}_{\mathcal{Z}_{N}}(\tilde{\mathcal{L}}_{\epsilon})- \hat{\mathfrak{R}}_{\mathcal{Z}_{N}}(\mathcal{L})\leq   \frac{ \epsilon }{\sqrt{N}} W n^{1-1/s} \). We can obtain this bound by noting in \eqref{eq:ARC-RC-linear} that \(\Theta \subseteq  \Theta_* = \{\theta\in\mathbb{R}^n\mid \norm{\theta}_{1} \leq c =  W n^{1-1/s}  \} \).
    \item[$\operatorname{(b)}$] In \citet[Thm 4]{awasthi2020adversarial}, authors consider arbitrary attack \(d(z',z) = \norm{x'-x}_{r} + \infty\abs{y'-y}\) and show that if  \(\Theta= \{\theta\in\mathbb{R}^n\mid \norm{\theta}_{s} \leq  W  \}  \) then \(\hat{\mathfrak{R}}_{\mathcal{Z}_{N}}(\tilde{\mathcal{L}}_{\epsilon})- \hat{\mathfrak{R}}_{\mathcal{Z}_{N}}(\mathcal{L})\leq  \frac{ \epsilon }{\sqrt{N}} W \max\{ n^{1-1/r-1/s},1\}  \). We can obtain this bound by noting in \eqref{eq:ARC-RC-linear} that \(\Theta \subseteq \Theta_*= \{\theta\in\mathbb{R}^n\mid \norm{\theta}_{1/(1-1/r)} \leq c = W \max\{ n^{1-1/r-1/s},1\}  \} \). Therefore, ARC-RC gap is dimension-free if \( 1/r+1/s\geq 1\).
\end{itemize}

\begin{example}
    [Adversarial complexity gap for MLPs] \label{exam:ACG-MLP} Given Assumption~\ref{asmp:Zdl} where \(\kappa=\infty\), then \(\mathcal{A}_{\theta}\) is an adversarial score (Lipschitz) of MLP \(f_{\theta} \colon x \mapsto  f (W_{K} f( \cdots f(W_1 x)\cdots))\)
    where \(\mathcal{A}_{\theta}(t) = t\operatorname{Lip}_f^{K}   \sup_{\theta\in\Theta}\Pi_{k=1}^{K}\norm{W_k}_r\). 
    By Theorem~\ref{thm:main0} and Theorem~\ref{thm:ACG},  
    \begin{equation}\label{eq:ARC-RC-MLP}
        \abs{\hat{\mathfrak{R}}_{\mathcal{Z}_{N}}(\tilde{\mathcal{L}}_{\epsilon})- \hat{\mathfrak{R}}_{\mathcal{Z}_{N}}(\mathcal{L})} \leq  \mathcal{A}_{\theta}(\epsilon) .
    \end{equation}
\end{example}

Similar to Example~\ref{exam:concave-complexity}, \eqref{eq:ARC-RC-MLP} successfully decouples the complexity of the network from its architecture, eliminating dependencies on both width $h$ and depth $K$, while trading the traditional asymptotic rate of $1/\sqrt{N}$ for the geometric transport budget $\epsilon \approx \mathcal{O}(N^{-1/n})$.  Specifically, in \citet[Thm 7]{awasthi2020adversarial} and \citet{xiao2022adversarial}, besides the standard Lipschitz and weight norm bounds, the ARC-RC gap carries aggressive scaling factors such as $\frac{(\text{diam}(\mathcal{Z}_N) + \epsilon)}{\sqrt{N}} \max\{n^{1-1/s-1/r}, 1\} \times m \sqrt{n}$ (one-layer network) or $\frac{(\text{diam}(\mathcal{Z}_N) + \epsilon)}{\sqrt{N}} \times \sqrt{K \log(K)}h$. In contrast, our bound relies solely on $\epsilon \approx \mathcal{O}(N^{-1/n})$ and completely bypasses the data diameter \(\text{diam}(\mathcal{Z}_N)\). For example, in GPT-3 model mentioned above,   $\frac{\sqrt{K \log(K)}h}{\sqrt{N}} = \frac{\sqrt{96 \log(96)} \times 12288}{\sqrt{10^9}} \approx 8.13$ and ours $\epsilon \approx \mathcal{O}(N^{-1/n}) \approx 0.998$.

\section{Numerical Experiments}\label{sec:exp}

\subsection{Traffic Data Regression} \label{sec:madrid}

\begin{figure}[H]
    \centering
    \includegraphics[width=0.95\linewidth]{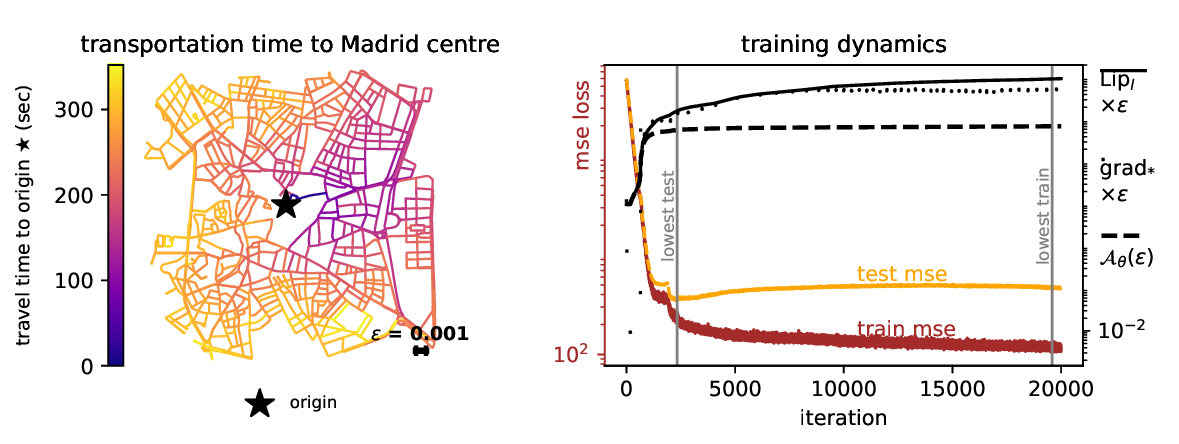}
    \caption{Left: Transportation time heatmap from location to origin. Right: Training dynamics of losses and certificates at \(\epsilon=10^{-3}\).}
    \label{fig:madrid-1}
\end{figure}

\textbf{(Figure~\ref{fig:madrid-1} - Left).} We obtain the Madrid road network from the osmnx open-source library, which provides a graph of nodes and edges. For \(1000\) random nodes, let $X^{(i)} \in \mathbb{R}^2$ represents  geospatial coordinates and $Y^{(i)} \in \mathbb{R}$ represents the shortest travel time from \(X^{(i)}\) to the city center ($\star$ in Figure~\ref{fig:madrid-1}). We train a hypothesis $f_\theta\colon \mathbb{R}^2\rightarrow\mathbb{R} $ (a multi-layer perceptron with 2 layers, 16 neurons each and Tanh activation) using absolute deviation loss \(\bm{l}(Z;\theta) = \abs{Y-f_{\theta}(X)} \). We set the parameter \(\kappa=10^{-4}\) and \(r=2\).

\textbf{Training Dynamic (Figure~\ref{fig:madrid-1} - Right).} During the training process, we monitor training loss, testing losses and values of three certificates: Lipschitz certificate \(\operatorname{Lip}\times\epsilon\) \citep{blanchet2019robust,gao2022wasserstein}, gradient-based certificate \(\operatorname{grad}_{*}\times\epsilon\) where \(\operatorname{grad}_{*}= \left( \mathbb{E}_{\mathbb{P}_N} [ \norm{\nabla_x \bm{l} (Z;\theta)}_{*}^q ] \right) ^{1/q}\) \citep{bartl2021sensitivity,bai2024wasserstein}, and  our proposed adversarial score \(\mathcal{A}_{\theta}(\epsilon)\) at \(\epsilon=10^{-3}\). We can see that all three certificates increase as training loss decreases, indicating higher sensitivity to input noise. 
However, our adversarial score is more stable and less volatile than the  $\operatorname{grad}_*$  certificate, and significantly tighter than the Lipschitz bound.

\begin{figure}[h]
    \centering
    \includegraphics[width=0.8\linewidth]{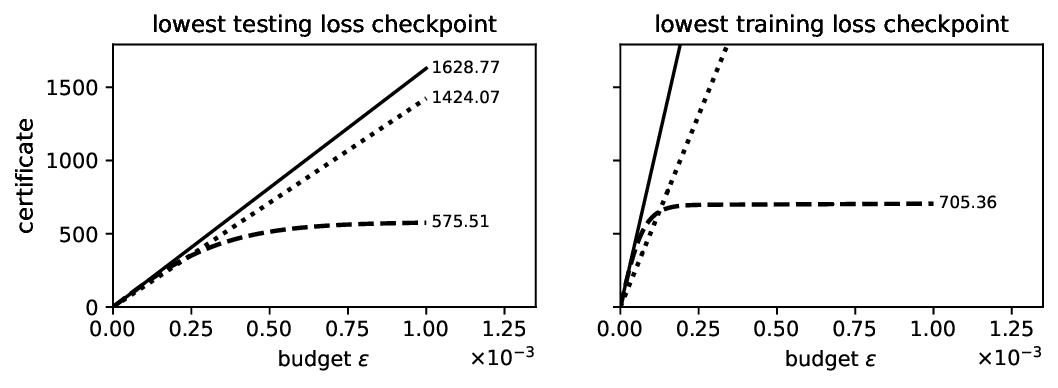}
    \caption{Dynamics of \(\operatorname{Lip}\times\epsilon\) \citep{blanchet2019robust,gao2022wasserstein} (solid),  \(\operatorname{grad}_{*}\times\epsilon\) \citep{bartl2021sensitivity,bai2024wasserstein} (dotted), and \(\mathcal{A}_{\theta}(\epsilon)\) (dashed) at two checkpoints in Figure~\ref{fig:madrid-1}.}
    \label{fig:madrid-2}
\end{figure}

\textbf{Budget Dynamic (Figure~\ref{fig:madrid-2}).} We analyze these three certificates across a noise budget range of \(\epsilon\in[0,10^{-3}] \) at the checkpoints for the lowest training and testing losses. First, we observe that our proposed \(\mathcal{A}_{\theta}(\epsilon)\)  are strictly tighter than the existing Lipschitz certificate. Besides, the  $\operatorname{grad}_*$  certificate serves only as a first-order estimation when \(\epsilon\) is small rather than a theoretical upper bound, i.e., it can  underestimate or overestimate the true risk. Furthermore, our non-linear certificate effectively captures the behavior of the Tanh activation: it magnifies small noise but saturates as the noise level increases.

\subsection{Generalization Capability of Adversarial Learning} \label{sec:mnist}

We evaluate our theoretical findings on the MNIST dataset, building upon the experimental design of \citet{yin2019rademacher}. The training set is subsampled to $N = 1000$ images. We consider the standard cross-entropy loss and $f_{\theta}$ is a Convolutional Neural Network (CNN) + Fully Connected Layer (FCL) of CNN-32 + CNN-64 + FCL-1024 + FCL-10. This corresponds to depth $K=4$, width $h=1024$, and  $n=784$. From this baseline, we conduct three independent sets of experiments (a), (b) and (c) by varying \ $h$, $K$ and $n$.

\begin{figure}[h!]
    \centering
    
    % --- Subfigure (a) ---
    \begin{subfigure}[b]{0.9\textwidth}
        \centering
        \includegraphics[width=\textwidth]{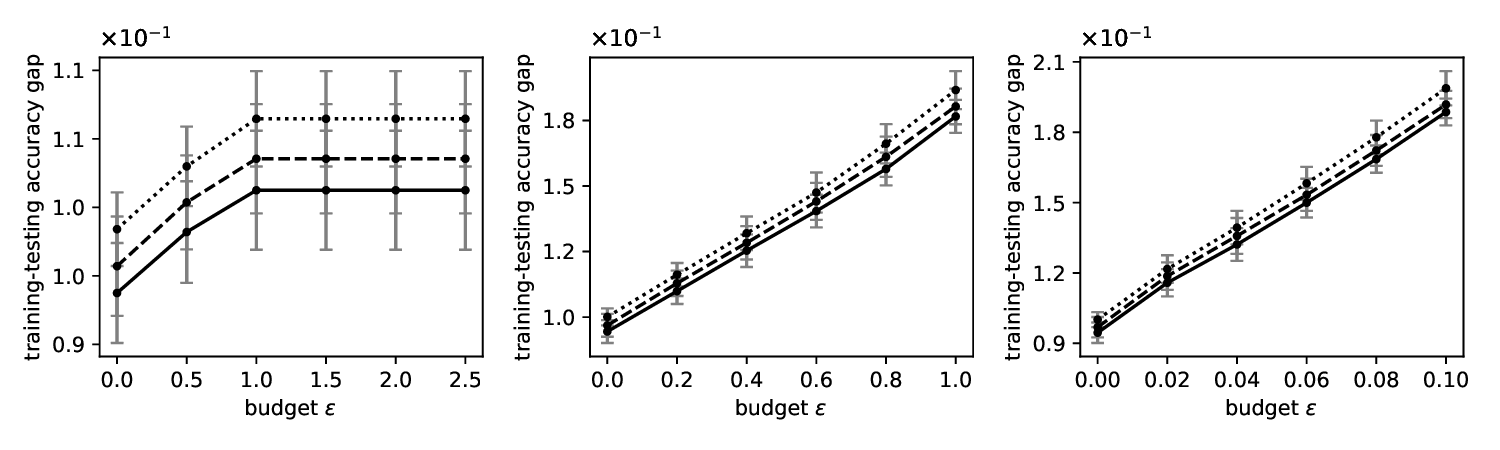}
        \caption{Varying width $h=512$ (dotted), $h=1024$ (dashed) and $h=512$ (solid), fixing dimension $n=784$ and depth $K=4$.}
        \label{fig:mnist-cnn-width}
    \end{subfigure}
    
    % --- Subfigure (b) ---
    \begin{subfigure}[b]{0.9\textwidth}
        \centering
        \includegraphics[width=\textwidth]{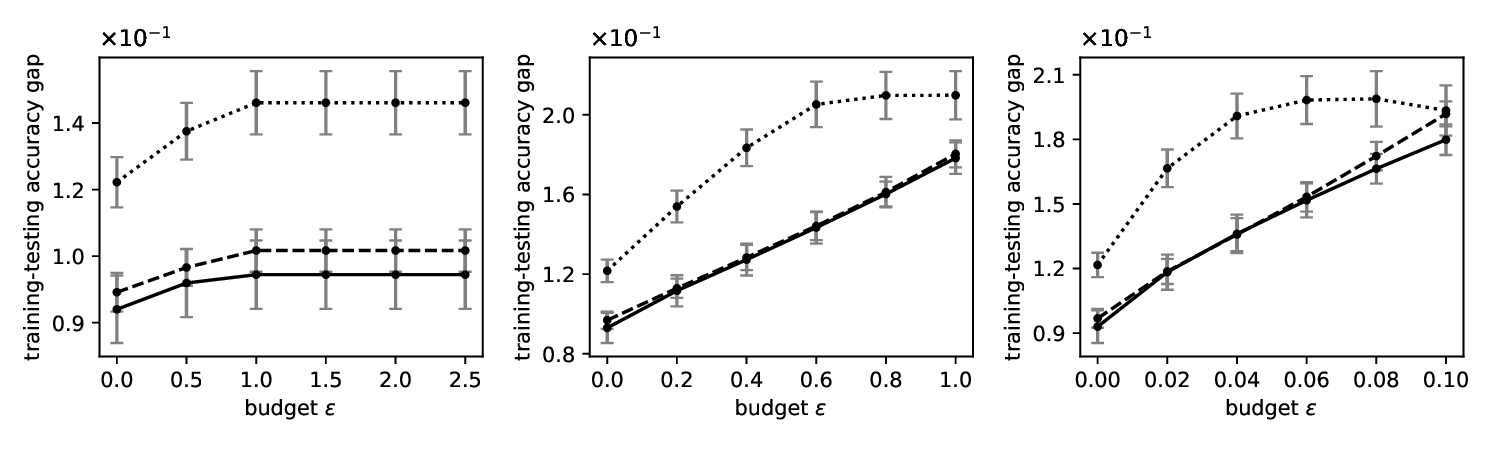}
        \caption{Varying depth $K=2$ (dotted), $K=4$ (dashed) and $K=6$ (solid), fixing dimension $n=784$ and width $h=1024$.}
        \label{fig:mnist-cnn-depth}
    \end{subfigure}

    % --- Subfigure (c) ---
    \begin{subfigure}[b]{0.9\textwidth}
        \centering
        \includegraphics[width=\textwidth]{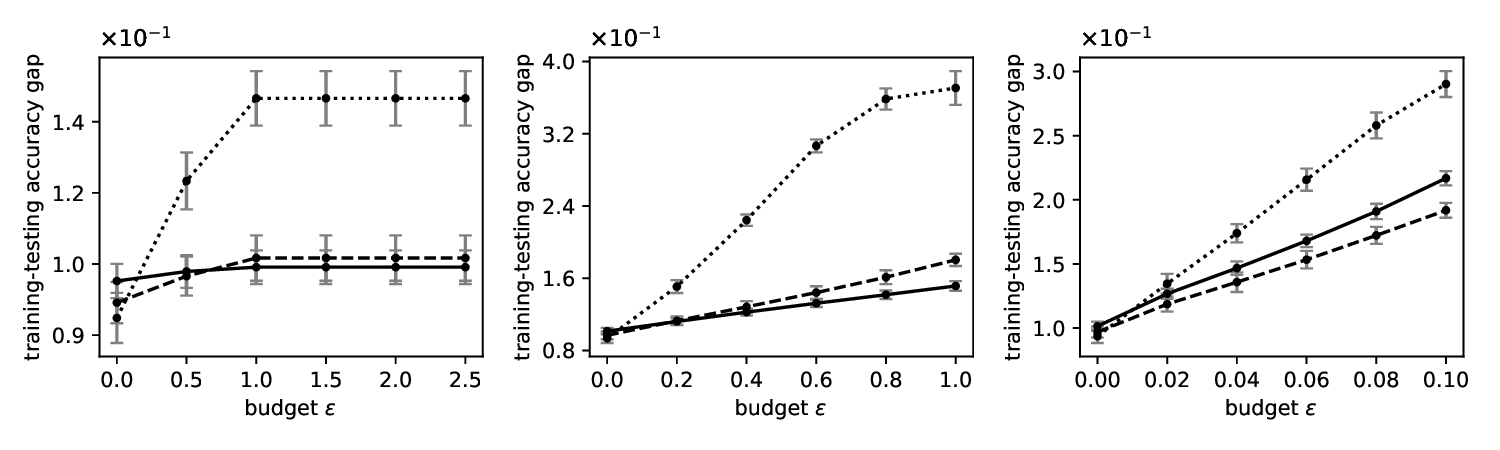}
        \caption{Varying dimension $n=196$ (dotted), $n=784$ (dashed), and $n=3136$ (solid), fixing width $h=1024$ and depth $K=4$.}
        \label{fig:mnist-cnn-dim}
    \end{subfigure}
    
    \caption{Generalization capability of CNN on MNIST. From left to right: $r=1$, $r=2$ and $r=\infty$.}
    \label{fig:mnist-cnn}
\end{figure}

\textbf{Adversarial Training} To find the optimal learning weights, we utilize the adversarial method. The adversarial perturbation is constrained by the \(r\)-norm where \(r=\{1,2,\infty\}\). Given an input \(z=(x,y)\), the adversarial example \(\tilde{z}= (\tilde{x},y)\) is generated as \(\tilde{x} = \text{clip}_{[0,1]}\left( x + \epsilon \cdot \Phi (\nabla_x \bm{l}(z;\theta)) \right)\), where \( \Phi(x)=\operatorname{sign}(x_{j_{\max}}) \bm{e}_{j_{\max}} \text{ if } r = 1,\) \( \Phi(x)=x / \norm{x}_2  \text{ if } r = 2\) ($\operatorname{grad}_*$), and \( \Phi(x)=\operatorname{sign}(x)  \text{if } r = \infty\) with \(j_{\max}\coloneqq \arg\max_j \abs{x_j} \) (FGSM). For each budget \(\epsilon\), we conduct 10 independent runs,  record the mean  and  standard deviation of the training-testing accuracy gap, and report them  in Figure~\ref{fig:mnist-cnn}.

\textbf{Analysis} Results in  Figure~\ref{fig:mnist-cnn} provides numerical validation for our above theoretical analysis. In contrast to traditional generalization bounds that scale heavily with network size, our empirical analysis confirms that increasing in depth ($K$) and width ($h$) do not cause the generalization gap to blow up. Besides, a smaller input dimension ($n=196$) actually yields a larger generalization gap compared to a  higher dimension ($n=3136$), aligning with our Example~\ref{exam:ACG-linear}+\ref{exam:ACG-MLP} that the ambient dimension $n$ is not an isolated  penalty on the complexity gap, but is rather  absorbed into the geometry of the optimal transport budget  \(\epsilon\approx\mathcal{O}(N^{-1/n})\). Finally, the trajectory of generalization gaps against \(\epsilon\) is clearly concave in many cases, reflecting our theoretical geometric certificate \(\operatorname{cc}\) in Theorem~\ref{thm:main0}.

\section{Conclusion}\label{sec:con}

In this paper, we proposed a novel framework using geometric  certificates to establish tight distributionally robust risk bounds, applicable even to non-Lipschitz and non-differentiable losses. This approach yields a worst-case generalization bounds and introduces a tractable adversarial score for layer-wise deep network analysis. Our comprehensive experiments numerically validated these findings, confirming that our proposed certificates are strictly tighter and more stable than traditional Lipschitz or gradient-based methods. This work   opens several  challenges for the broader community. Key theoretical directions include mitigating the curse of dimensionality within the optimal transport budget and developing scalable approximations for the NP-hard Exact CD condition. In addition, leveraging the proposed adversarial score to design robust and trustable neural networks also presents a highly promising topic for future research.

\begin{appendices}

\section{Proofs of Theorems and Propositions}
\subsection{Proof of Theorem~\ref{thm:main0}} \label{proof:main0}
The cases of \(\operatorname{lb}_{\infty}\) and \(\operatorname{cc}_{\infty}\) are immediate results by Lemma~\ref{lem:p-infinity}. We now consider \(p\in[1,\infty) \).

\noindent\textbf{Lower Bound.}  
For any \(\tilde{Z}^{(i)}\in\mathcal{Z}\) and \(\eta_i\in[0,1]\) where \( i = 1,\dots,N\), let  \(\tilde{\mathbb{P}}\in\mathcal{P}(\mathcal{Z})\) and \(\tilde{\pi}\in\Pi(\tilde{\mathbb{P}},\mathbb{P}_N) \) be defined as
\begin{equation*}
	\begin{array}{l}
	    \tilde{\mathbb{P}} \coloneqq \sum\limits_{i=1}^{N}\mu_i (1-\eta_i) \bm{\chi}_{\{Z^{(i)}\}} + \mu_i\eta_i \bm{\chi}_{\{\tilde{Z}^{(i)}\}},\quad\,
     \tilde{\pi}\coloneqq \sum\limits_{i=1}^{N}\mu_i (1-\eta_i) \bm{\chi}_{\{\left(Z^{(i)},Z^{(i)}\right)\}}+ \mu_i\eta_i \bm{\chi}_{\{\left(\tilde{Z}^{(i)},Z^{(i)}\right)\}}.
	\end{array}
\end{equation*}
Then the loss expectation with respect to \(\tilde{\pi}\) is given by
\begin{equation} \label{eq:WE}
	\begin{array}{rl}
		\mathbb{E}_{\tilde{\mathbb{P}}} [\bm{l} (Z;\theta) ]  
  &=\sum\limits_{i=1}^{N}\mu_i (1-\eta_i) \bm{l}(Z^{(i)};\theta) 
  + \mu_i\eta_i \bm{l}(\tilde{Z}^{(i)};\theta)\\
		&= \mathbb{E}_{{\mathbb{P}_N}} [\bm{l} (Z;\theta) ]+ \sum\limits_{i=1}^{N}\mu_i\eta_i \left(\bm{l}(\tilde{Z}^{(i)};\theta) -\bm{l}(Z^{(i)};\theta)   \right),
	\end{array}
\end{equation}
and
 \begin{equation*}       \begin{array}{ll}
        \mathcal{W}_{p}\left(\tilde{\mathbb{P}},\mathbb{P}_N\right)&\leq			 \left(\int_{\mathcal{Z}\times\mathcal{Z}}d^{p}(\tilde{z},z)\mathrm{d}\tilde{\pi}(\tilde{z},z)\right)^{1/p}
    = \left(\sum_{i=1}^{N}\mu_i\eta_i d^{p}\left(\tilde{Z}^{(i)},Z^{(i)} \right)\right)^{1/p}.
    \end{array}
    \end{equation*}
Hence the following optimization problem yields a lower bound of \(\mathcal{R}_{p}(\epsilon) = \sup_{\mathbb{P} \colon \mathcal{W}_{p}(\mathbb{P},\mathbb{P}_N) \leq \epsilon }  \ \mathbb{E}_{\mathbb{P}} [\bm{l} (Z;\theta) ]\).
\begin{equation}\label{eq:opt-1}
	\begin{array}{lll}
		\sup & \mathbb{E}_{{\mathbb{P}_N}} [\bm{l} (Z;\theta) ] 
		+  \sum\limits_{i=1}^{N}\mu_i\eta_i \left( \bm{l}(\tilde{Z}^{(i)};\theta) -\bm{l}(Z^{(i)};\theta) \right)\\
		& \text{such that }\eta_i \in [0,1], \,  \tilde{Z}^{(i)}\in\mathcal{Z}, \, i = 1,\dots,N, \\
		& \text{and }\left(\sum_{i=1}^{N}\mu_i\eta_i d^{p}\left(\tilde{Z}^{(i)},Z^{(i)} \right)\right) \leq \epsilon^{p}.
	\end{array} 
\end{equation}
 Let \(t_i=\frac{\epsilon}{\eta_i^{1/p}}\), or equivalently, \( \eta_i = \epsilon^{p}/t_i^{p} \). Then  \(t_i\in[\epsilon,\infty) \) implies \(\eta_i\in(0,1]\), and  \(\eqref{eq:opt-1}\geq\eqref{eq:opt-2} \) where 
\begin{equation}\label{eq:opt-2}
	\begin{array}{lll}
		\sup & \mathbb{E}_{{\mathbb{P}_N}} [\bm{l} (Z;\theta) ]   
		+   \sum\limits_{i=1}^{N}\mu_i \frac{\epsilon^{p}}{t_i^{p}} \left( \bm{l}(\tilde{Z}^{(i)};\theta) -\bm{l}(Z^{(i)};\theta)   \right)\\
		& \text{such that }t_i \in [\epsilon,\infty), \,  \tilde{Z}^{(i)}\in\mathcal{Z},\,  i = 1,\dots,N\\
  &\text{and } \sum_{i=1}^{N}\mu_i\frac{\epsilon^{p}}{t_i^{p}} d^{p}\left(\tilde{Z}^{(i)},Z^{(i)} \right) \leq \epsilon^{p}.
	\end{array} 
\end{equation}
Let \(\rho\) be an arbitrary positive scalar. For every \(i=1,\dots,N\), by definition of the individual rate \(\Delta_{\theta}(Z^{(i)},\epsilon)\) \eqref{eq:ind-rate}, there exists \(\tilde{Z}_{\rho}^{(i)}\in\mathcal{Z}\) such that \(d\left(\tilde{Z}_{\rho}^{(i)},Z^{(i)} \right)  \leq t_i\) and
\begin{equation*}
	\bm{l}(\tilde{Z}_{\rho}^{(i)};\theta) -\bm{l}(Z^{(i)};\theta)\geq  \Delta_{\theta}(Z^{(i)},t_i) - \rho.
\end{equation*}
This implies that \(\sum_{i=1}^{N}\mu_i\frac{\epsilon^{p}}{t_i^{p}} d^{p}\left(\tilde{Z}_{t_i,\rho}^{(i)},Z^{(i)} \right) \leq \sum_{i=1}^{N}\mu_i\frac{\epsilon^{p}}{t_i^{p}} t_i^{p}= \epsilon^{p} \). Therefore, \( \left\{t_i\in [\epsilon,\infty), \tilde{Z}^{(i)}=\tilde{Z}_{t_i,\rho}^{(i)}\right\}_{i=1}^{N} \)  is a feasible solution of \eqref{eq:opt-2} and hence,
\begin{equation*}
	\begin{array}{rl}
		\eqref{eq:opt-2} 
		\geq  \mathbb{E}_{{\mathbb{P}_N}} [\bm{l} (Z;\theta) ]  +  \sum\limits_{i=1}^{N}\mu_i\sup\limits_{t_i\in[\epsilon,\infty)}   
		\left\{ \frac{\epsilon^{p} \Delta_{\theta}(Z^{(i)},t_i)}{t_i^{p}}  -  \frac{\epsilon^{p}}{t_i^{p}}\rho\right\} 
		\geq \mathbb{E}_{{\mathbb{P}_N}} [\bm{l} (Z;\theta) ]  + \sum\limits_{i=1}^{N}\mu_i s^{(i)}(\epsilon) -  \rho,
	\end{array}
\end{equation*} 
where \(s^{(i)}(\epsilon)= \sup_{t_i\in[\epsilon,\infty)}  
\left\{ \frac{\epsilon^{p} \Delta_{\theta}(Z^{(i)},t_i)}{t_i^{p}}  \right\}  = \sup_{u\in[\epsilon^p,\infty)}  
\left\{ \frac{\epsilon^{p} \Delta_{\theta}(Z^{(i)},u^{1/p})}{u}  \right\}\). By Lemma~\ref{lem:calculation}, \(s^{(i)}(\epsilon) = \mathcal{S}_{f^{(i)}}(\epsilon^p)\). This holds true for every \(\rho>0\). Therefore,
\begin{equation*}
	\textstyle\mathcal{R}_{p}(\epsilon)  \geq  \eqref{eq:opt-1} \geq \eqref{eq:opt-2}  \geq  \mathbb{E}_{{\mathbb{P}_N}} [\bm{l} (Z;\theta) ] 
	+    \sum\limits_{i=1}^{N}  \mu_i \mathcal{S}_{f^{(i)}}( \epsilon^p).
\end{equation*} 

\noindent\textbf{Upper Bound.}  By the  definition of the maximal rate \(\Delta_{\theta}^{\max}\) in \eqref{eq:max-rate}, we have that  \[\sup_{z'\in\mathcal{Z},\hat{z}\in\mathcal{Z}_N} \left\{ \bm{l}(z';\theta) - \bm{l}(\hat{z};\theta)  \mid d(z',\hat{z}) \leq  t \right\} = \Delta_{\theta}^{\max}(t)\] for any \(t\geq0\). Since \(\mathcal{C}_{f^{\max}}(t)\) is the least concave majorant of \(\Delta_{\theta}^{\max}(t^{1/p}) \), we have that \(\Delta_{\theta}^{\max}(t^{1/p})\leq \mathcal{C}_{f^{\max}}(t)\) and thus 
 \begin{equation*} 
     \sup_{z'\in\mathcal{Z},\hat{z}\in\mathcal{Z}_N}         \left\{ \bm{l}(z';\theta) - \bm{l}(\hat{z};\theta)   \mid  d^{p}(z',\hat{z}) \leq  t \right\} \leq \mathcal{C}_{f^{\max}}(t).
 \end{equation*}
 This implies that for any \(z'\in\mathcal{Z},\hat{z}\in\mathcal{Z}_N\),
 \begin{equation}\label{eq:Hd}
     \bm{l}(z';\theta) - \bm{l}(\hat{z};\theta) \leq \mathcal{C}_{f^{\max}}(d^{p}(z',\hat{z})).
 \end{equation}
Let \(\rho>0\) be an arbitrary scalar. For any \(  \tilde{\mathbb{P}} \in \mathcal{P}(\mathcal{Z})  \) such that \(\mathcal{W}_{p}(\tilde{\mathbb{P}},\mathbb{P}_N) \leq \epsilon \), by the definition of $\mathcal{W}_{p}$, there exists \(\tilde{\pi}\in\Pi(\tilde{\mathbb{P}},\mathbb{P}_N ) \)  such that 
	\begin{equation}\label{eq:Wp-delta}
		\textstyle\left(\int_{\mathcal{Z}\times\mathcal{Z}_N} d^{p}(\tilde{z},z) \mathrm{d}\tilde{\pi}(\tilde{z},z)\right)^{1/p} \leq \epsilon + \rho.
	\end{equation}
Note that \(\mathbb{E}_{\tilde{\mathbb{P}}} [\bm{l}(Z;\theta)] =\int_{\mathcal{Z}} \bm{l}(\tilde{z};\theta)\mathrm{d}\tilde{\mathbb{P}}(\tilde{z}) =\int_{\mathcal{Z}\times\mathcal{Z}_N} \bm{l}(\tilde{z};\theta)\mathrm{d}\tilde{\pi}(\tilde{z},z) \) and \[\int_{\mathcal{Z}\times\mathcal{Z}_N} \bm{l}({z};\theta)\mathrm{d}\tilde{\pi}(\tilde{z},z)= \int_{\mathcal{Z}} \bm{l}({z};\theta)\mathrm{d}{\mathbb{P}_N}(z)=\mathbb{E}_{{\mathbb{P}_N}} [\bm{l}(Z;\theta)] \].

\noindent Therefore, \(\mathbb{E}_{\tilde{\mathbb{P}}} [\bm{l}(Z;\theta)]
			=\mathbb{E}_{{\mathbb{P}_N}} [\bm{l}(Z;\theta)]   +\int_{\mathcal{Z}\times\mathcal{Z}_N} \left(\bm{l}(\tilde{z};\theta)-\bm{l}({z};\theta)\right)\mathrm{d}\tilde{\pi}(\tilde{z},z)\) and

	\begin{equation*}
		\begin{array}{rll}
	            & \int\nolimits_{\mathcal{Z}\times\mathcal{Z}_N} \left(\bm{l}(\tilde{z};\theta)-\bm{l}({z};\theta)\right)\mathrm{d}\tilde{\pi}(\tilde{z},z)
			&\leq  \int\nolimits_{\mathcal{Z}\times\mathcal{Z}_N} \mathcal{C}_{f^{\max}}\left(d^{p}(\tilde{z},z)\right)\mathrm{d}\tilde{\pi}(\tilde{z},z) \\
			\leq& \mathcal{C}_{f^{\max}}\left(\int\nolimits_{\mathcal{Z}\times\mathcal{Z}_N} d^{p}(\tilde{z},z)\mathrm{d}\tilde{\pi}(\tilde{z},z) \right) 
			&\leq  \mathcal{C}_{f^{\max}}\left((\epsilon+\rho)^{p}\right),
		\end{array}
	\end{equation*} 
	where the first inequality follows from \eqref{eq:Hd}; the second equality follows from the fact that the least concave majorant \(\mathcal{C}_{f^{\max}}\) is concave; the last inequality follows from \eqref{eq:Wp-delta} and the fact that \(\mathcal{C}_{f^{\max}}\) is non-decreasing (Lemma~\ref{lem:univariate-majorant}). This means that for any \(\rho>0\) and any \(  \tilde{\mathbb{P}} \in \mathcal{P}(\mathcal{Z})  \) such that \(\mathcal{W}_{p}(\tilde{\mathbb{P}},\mathbb{P}_N) \leq \epsilon \), we have shown that \(\mathbb{E}_{\tilde{\mathbb{P}}} [\bm{l}(Z;\theta)] \leq \hat{\mathcal{R}} + \mathcal{C}_{f^{\max}}\left((\epsilon+\rho)^{p}\right)\). Thus
	\begin{equation*}
		\begin{array}{ll}
		    \mathcal{R}_{p}(\epsilon)      &\leq  \hat{\mathcal{R}} + \mathcal{C}_{f^{\max}}\left((\epsilon+\rho)^{p}\right).
		\end{array}
	\end{equation*}
	Since \(\mathcal{C}_{f^{\max}}\) is concave on \((0,\infty)\), it is also continuous of \((0,\infty)\) and by letting \(\rho\rightarrow0\), we have the desired conclusion.

\subsection{Proof of Proposition~\ref{prop:ExactCD}} \label{proof:ExactCD} 
We first calculate the worst-case zero-one loss, denoted as $\tilde{\bm{l}}_{\epsilon}(Z^{(i)};\theta)$, for each empirical point at radius $\epsilon$.
\begin{equation*}
    \begin{aligned}
        \tilde{\bm{l}}_{\epsilon}(Z^{(i)};\theta) &\coloneqq \sup_{d(z',Z^{(i)}) \leq  \epsilon} \bm{l}(z';\theta) = \bm{l}(Z^{(i)};\theta) + \Delta_{\theta}(Z^{(i)},\epsilon) \\
        &= \begin{cases}
            1 & \text{ if } \exists x  \text{ s.t. } d_{\mathcal{X}}(x,X^{(i)})\leq \epsilon \text{ and } f_{\theta}(x) \ne Y^{(i)}, \\
            0 & \text{ if } f_{\theta}(x) = Y^{(i)} \; \forall x \text{ s.t. } d_{\mathcal{X}}(x,X^{(i)})\leq \epsilon.
        \end{cases}
    \end{aligned}
\end{equation*}

\noindent\textbf{Exact CD implies Robust.} Suppose that $\mathbb{P}_N$ is $(\epsilon,\delta)$-Exact CD with respect to $f_\theta$. This means there exist subsets $\Omega_k$ satisfying the $\delta$-coverage and $\epsilon$-immunity properties. By Lemma 1, the total robust risk is:
\begin{equation*}
    \begin{array}{c}
       \hat{\mathcal{R}}_N + \operatorname{lb}_{\infty}(\epsilon) = \sum\limits_{i=1}^N \mu_i \tilde{\bm{l}}_{\epsilon}(Z^{(i)};\theta) 
    =  \sum\limits_{i \colon X^{(i)} \in \Omega_{Y^{(i)}}} \mu_i \tilde{\bm{l}}_{\epsilon}(Z^{(i)};\theta) + \sum\limits_{i \colon X^{(i)} \notin \Omega_{Y^{(i)}}} \mu_i \tilde{\bm{l}}_{\epsilon}(Z^{(i)};\theta).
    \end{array}
\end{equation*}
For any $X^{(i)} \in \Omega_{Y^{(i)}}$, the $\epsilon$-immunity property guarantees that its entire $\epsilon$-neighborhood is classified as $Y^{(i)}$. Thus, $\tilde{\bm{l}}_{\epsilon}(Z^{(i)};\theta) = 0$, and the first term vanishes. Since $\tilde{\bm{l}}_{\epsilon}(Z^{(i)};\theta)$ is upper bounded by $1$, the second term is bounded by the mass of the points outside the coverages. By the $\delta$-coverage property:
\begin{equation*}
    \begin{array}{c}
         \hat{\mathcal{R}}_N + \operatorname{lb}_{\infty}(\epsilon) \leq 0 + \sum\limits_{i \colon X^{(i)} \notin \Omega_{Y^{(i)}}} \mu_i = 1 - \sum\limits_{k=1}^m \sum\limits_{i \colon X^{(i)} \in \Omega_k} \mu_i \leq 1 - (1 - \delta) = \delta.
    \end{array}
\end{equation*}

Therefore, $f_\theta$ is $(\epsilon,\delta)$-robust.

\noindent\textbf{Robust implies Exact CD.} Conversely, suppose that $f_{\theta}$ is $(\epsilon,\delta)$-robust, meaning $\hat{\mathcal{R}}_N + \operatorname{lb}_{\infty}(\epsilon) \leq \delta$. We construct the subsets $\Omega_k$ as the collection of empirical points of class $k$ that have a zero worst-case loss:
$$
\Omega_k \coloneqq \{X^{(i)} \colon Y^{(i)} = k \text{ and } \tilde{\bm{l}}_{\epsilon}(Z^{(i)};\theta) = 0\}.
$$
Equivalently, $f_\theta(x) = k$ for any $x\in\Omega^{+\epsilon}_k$. Hence $\Omega_k^{+\epsilon} \subseteq \mathcal{D}_{f_\theta,k}$ ($\epsilon$-immunity). Now note that $\tilde{\bm{l}}_{\epsilon}(Z^{(i)};\theta)$ is a zero-one function; thus the probability mass of the robust points is the complement of the mass of the vulnerable points:
\begin{equation*}
    \begin{array}{c}
        \sum\limits_{k=1}^m \sum\limits_{i \colon X^{(i)} \in \Omega_k} \mu_i = \sum\limits_{i \colon \tilde{\bm{l}}_{\epsilon}(Z^{(i)};\theta) = 0} \mu_i = 1 - \sum\limits_{i \colon \tilde{\bm{l}}_{\epsilon}(Z^{(i)};\theta) = 1} \mu_i. 
    \end{array}
\end{equation*}
Besides, the total robust risk is exactly the mass of the vulnerable points:
\begin{equation*}
    \begin{array}{c}
         \hat{\mathcal{R}}_N + \operatorname{lb}_{\infty}(\epsilon) = \sum\limits_{i=1}^N \mu_i \tilde{\bm{l}}_{\epsilon}(Z^{(i)};\theta) = \sum\limits_{i \colon \tilde{\bm{l}}_{\epsilon}(Z^{(i)};\theta) = 1} \mu_i \leq \delta.
    \end{array}
\end{equation*}
Therefore, we obtain the $\delta$-coverage property as:

\begin{equation*}
    \begin{array}{c}
         \sum\limits_{k=1}^m \sum\limits_{i \colon X^{(i)} \in \Omega_k} \mu_i = 1 - \big(\hat{\mathcal{R}}_N + \operatorname{lb}_{\infty}(\epsilon)\big) \geq 1 - \delta. \hfill\Halmos
    \end{array}
\end{equation*}

\subsection{Proof of Adversarial Complexity Gap Theorem~\ref{thm:ACG}} \label{proof:ACG}

\noindent\textbf{ARC-RC Gap.} 
Since \(\tilde{\bm{l}}_{\epsilon}(Z^{(i)};\theta) = \bm{l}(Z^{(i)};\theta) + \Delta_{\theta}(Z^{(i)},\epsilon)\) and \(\sup(A+B)\leq \sup A + \sup B\), 
\begin{equation*}
\begin{array}{rll}
   &\hat{\mathfrak{R}}_{\mathcal{Z}_{N}}(\tilde{\mathcal{L}}_{\epsilon})= \mathbb{E}_{\sigma} \left[ \sup\limits_{\theta\in\Theta} \left( \frac{1}{N}   \sum\limits_{i=1}^{N} \sigma_{i} \tilde{\bm{l}}(Z^{(i)};\theta) \right) \right]    &= \mathbb{E}_{\sigma} \left[ \sup\limits_{\theta\in\Theta} \left( \frac{1}{N}   \sum\limits_{i=1}^{N} \sigma_{i} \left( \bm{l}(Z^{(i)};\theta) + \Delta_{\theta}(Z^{(i)},\epsilon) \right) \right) \right]  \\
   \leq & \hat{\mathfrak{R}}_{\mathcal{Z}_{N}}(\mathcal{L})+ \mathbb{E}_{\sigma} \left[ \sup\limits_{\theta\in\Theta} \left( \frac{1}{N}   \sum\limits_{i=1}^{N} \sigma_{i}\Delta_{\theta}(Z^{(i)},\epsilon)\right) \right]&= \hat{\mathfrak{R}}_{\mathcal{Z}_{N}}(\mathcal{L})+ \hat{\mathfrak{R}}_{\mathcal{Z}_{N}}({\Upsilon}_{\epsilon}).
\end{array}
\end{equation*}
Similarly, as  \(\sup(A+B)\geq \sup A + \inf B = \sup A - \sup (-B)\) and \(\sigma \sim -\sigma \), 
\begin{equation*}
\begin{array}{rll}
   &\hat{\mathfrak{R}}_{\mathcal{Z}_{N}}(\tilde{\mathcal{L}}_{\epsilon})= \mathbb{E}_{\sigma} \left[ \sup\limits_{\theta\in\Theta} \left( \frac{1}{N}   \sum\limits_{i=1}^{N} \sigma_{i} \tilde{\bm{l}}(Z^{(i)};\theta) \right) \right]    = \mathbb{E}_{\sigma} \left[ \sup\limits_{\theta\in\Theta} \left( \frac{1}{N}   \sum\limits_{i=1}^{N} \sigma_{i} \left( \bm{l}(Z^{(i)};\theta) + \Delta_{\theta}(Z^{(i)},\epsilon) \right) \right) \right]  \\
   \geq & \hat{\mathfrak{R}}_{\mathcal{Z}_{N}}(\mathcal{L})- \mathbb{E}_{\sigma} \left[ \sup\limits_{\theta\in\Theta} \left( \frac{1}{N}   \sum\limits_{i=1}^{N} -\sigma_{i}\Delta_{\theta}(Z^{(i)},\epsilon)\right) \right] 
   =  \hat{\mathfrak{R}}_{\mathcal{Z}_{N}}(\mathcal{L})- \mathbb{E}_{\sigma} \left[ \sup\limits_{\theta\in\Theta} \left( \frac{1}{N}   \sum\limits_{i=1}^{N} \sigma_{i}\Delta_{\theta}(Z^{(i)},\epsilon)\right) \right] \\
   = & \hat{\mathfrak{R}}_{\mathcal{Z}_{N}}(\mathcal{L})- \hat{\mathfrak{R}}_{\mathcal{Z}_{N}}({\Upsilon}_{\epsilon}).
\end{array}
\end{equation*}
Therefore, \(\abs{\hat{\mathfrak{R}}_{\mathcal{Z}_{N}}(\tilde{\mathcal{L}}_{\epsilon})- \hat{\mathfrak{R}}_{\mathcal{Z}_{N}}(\mathcal{L})}\leq \hat{\mathfrak{R}}_{\mathcal{Z}_{N}}({\Upsilon}_{\epsilon}) \). The inequality is completed by noting that \[\frac{1}{N} \sum_{i=1}^N \sigma_i \Delta_\theta(Z^{(i)}, \epsilon) \le \frac{1}{N} \sum_{i=1}^N \Delta_\theta(Z^{(i)}, \epsilon) = \operatorname{lb}_{\infty}(\epsilon)\].  

Now, assume that the rate is data-independent, i.e., \(\Delta_{\theta}(Z^{(i)},\epsilon) = \Delta_{\theta}^{\max}(\epsilon)\) for any \(i\).  Then

\begin{equation*}
\begin{array}{rll}
    \hat{\mathfrak{R}}_{\mathcal{Z}_{N}}({\Upsilon}_{\epsilon}) &= \mathbb{E}_{\sigma} \left[ \sup\limits_{\theta\in\Theta} \left( \frac{1}{N}   \sum\limits_{i=1}^{N} \sigma_{i}\Delta_{\theta}^{\max}(\epsilon)\right) \right] &= \mathbb{E}_{\sigma} \left[ \sup\limits_{\theta\in\Theta}\Delta_{\theta}^{\max}(\epsilon) \left( \frac{1}{N}   \sum\limits_{i=1}^{N} \sigma_{i}\right) \right] \\
    & \leq \sup\limits_{\theta\in\Theta}\Delta_{\theta}^{\max}(\epsilon)  \mathbb{E}_{\sigma} \left[  \abs{\frac{1}{N}   \sum\limits_{i=1}^{N} \sigma_{i}} \right] &\leq \frac{1}{\sqrt{N}}\sup\limits_{\theta\in\Theta}\Delta_{\theta}^{\max}(\epsilon) ,
\end{array}
\end{equation*}
where the last inequality follows Khintchine's inequality \citep{haagerup1981best}.

\noindent\textbf{ACC-CC Gap.} Recall that  rates of \(\bm{l}\) and \(\tilde{\bm{l}}\) are given by
\begin{equation*}
    \begin{array}{rl}
         \Delta_{\theta}(\hat{z},t) &= \sup\limits_{z'\in\mathcal{Z}} \left\{ \bm{l}(z';\theta) - \bm{l}(\hat{z};\theta) \mid d(z',\hat{z}) \leq  t \right\}, \\
         \tilde{\Delta}_{\theta}(\hat{z},t) &= \sup\limits_{z'\in\mathcal{Z}} \left\{ \tilde{\bm{l}}_{\epsilon}(z';\theta) - \tilde{\bm{l}}_{\epsilon}(\hat{z};\theta) \mid d(z',\hat{z}) \leq  t \right\}.
    \end{array}
\end{equation*}
One can rewrite the quantity of the second supremum  as
\begin{equation*}
    \begin{array}{rll}
         &\tilde{\bm{l}}_{\epsilon}(z';\theta) - \tilde{\bm{l}}_{\epsilon}(\hat{z};\theta)
        &= \sup_{u \colon d(u, z') \leq \epsilon} \bm{l}(u; \theta) - \sup_{v \colon d(v, \hat{z}) \leq \epsilon}
        \bm{l}(v; \theta) \\
        = & \Delta_{\theta}(z',\epsilon)+\bm{l}(z';\theta) - \Delta_{\theta}(\hat{z},\epsilon)-\bm{l}(\hat{z};\theta)
        &= \left[\Delta_{\theta}(z',\epsilon) - \Delta_{\theta}(\hat{z},\epsilon) \right] + \left[\bm{l}(z';\theta)-\bm{l}(\hat{z};\theta) \right].
    \end{array}
\end{equation*}
Taking supremum over all \(z'\) such that \(d(z',\hat{z}) \leq  t\), then the left-hand-side is the rate \(\tilde{\Delta}_{\theta} \) of \(\tilde{\bm{l}}\), the first term in the right-hand-side is the rate of \(\Delta_{\theta}(\cdot,\epsilon)\), and the second term is the rate of \(\bm{l}\). Using \(\sup(A+B)\leq \sup(A)+\sup(B) \), we have 
\begin{equation*}
    \tilde{\Delta}_{\theta}(\hat{z},t) \leq \Delta_{\Delta_{\theta}(\cdot,\epsilon)}(\hat{z},t) +  \Delta_{\theta}(\hat{z},t).
\end{equation*}
Taking the maximum over all \(\hat{z}\in\mathcal{Z}_N\), then supremum over all \(\theta\in\Theta\), we have \(\sup_{\theta\in\Theta}\tilde{\Delta}_{\theta}^{\max}(t)  \leq \sup_{\theta\in\Theta} \Delta_{\Delta_{\theta}(\cdot,\epsilon)}^{\max}(t) + \sup_{\theta\in\Theta}  \Delta_{\theta}^{\max}(t)\). 
Taking least concave majorant and using Lemma~\ref{lem:calculation} and using notation \({\Upsilon}_{\epsilon}\coloneqq \{z\mapsto \Delta_{\theta}(z,\epsilon) \mid \theta\in\Theta \} \), we obtain 
\begin{equation*}
    \hat{\mathfrak{C}}_N(\tilde{\mathcal{{L}}}_{\epsilon},t) \leq \hat{\mathfrak{C}}_N({\Upsilon}_{\epsilon},t) +\hat{\mathfrak{C}}_N(\mathcal{L},t). 
\end{equation*}
The proof is completed. \hfill\Halmos

\subsection{Proof of Classification Proposition~\ref{prop:classification}} \label{proof:classification}

\noindent\textit{Proof.} To show that \(\mathcal{A}_{\theta}\) is an adversarial score, we need to verify that \(\mathcal{A}_{\theta}\) is concave and for any \(t\geq0\),
\begin{equation*}
\begin{array}{c}
     \sup_{z'\in\mathcal{Z},\hat{z}\in\mathcal{Z}_N}       \left\{ \bm{l}(z';\theta) - \bm{l}(\hat{z};\theta) \mid d(z',\hat{z}) \leq  t \right\} \leq \mathcal{A}_{\theta}(t).
\end{array}
  \end{equation*}
Note that \(y\in\epsilon^{m}\) is a probability vector, thus \(\norm{y}_{s}\leq 1\) for any \(s\in[1,\infty]\). 
  \begin{enumerate}
      \item[(a)] Since \(\kappa=\infty\), \(d(z',z)\leq t\) if and only if \(y'=y\) and \( \norm{x'-x}_{r}\leq t \).   In that case, we have \( \bm{l}(z';\theta) - \bm{l}(\hat{z};\theta) = \inprod{y, f_{\theta}(x')-f_{\theta}(x)}\leq \norm{y}_{1/(1-1/r)}\norm{f_{\theta}(x')-f_{\theta}(x)}_r\leq F_{\theta}(\norm{x'-x}_r)\leq  F_{\theta}(t) \). As \(F_{\theta}\)  is an adversarial score of \(f_{\theta}\), it is concave, and therefore \(F_{\theta}\) is an adversarial score of \(\bm{l}\).
      \item[(b)] We can decompose the loss difference into two components by \(\bm{l}(z';\theta) - \bm{l}(\hat{z};\theta) =\inprod{y', f_{\theta}(x')-f_{\theta}(x)}+\inprod{y'-y,f_{\theta}(x)}\). For any \(z,z'\) such that \(d(z',z) \leq t\), it is equivalent to \(\norm{x'-x}_{r} \leq t- \kappa\norm{y'-y}_{1} = t - \tau\) where \(\tau\coloneqq\kappa\norm{y'-y}_{1}\in[0,t]\). Thus the first component \(\inprod{y', f_{\theta}(x')-f_{\theta}(x)} \leq \norm{y'}_{1/(1-1/r)}\norm{f_{\theta}(x')-f_{\theta}(x)}_r\leq F_{\theta}(\norm{x'-x}_{r})\leq F_{\theta}(t-\tau)  \) and the second component \(\inprod{y'-y,f_{\theta}(x)} \leq \norm{y'-y}_{1}\norm{f_{\theta}(x)}_{\infty}\leq  M\norm{y'-y}_{1}= M\kappa^{-1}\tau  \). Hence,
      \begin{equation*}
          \bm{l}(z';\theta) - \bm{l}(\hat{z};\theta) 
              \leq  F_{\theta}(t-\tau) + M\kappa^{-1}\tau .
      \end{equation*}
      Therefore, \(\bm{l}(z';\theta) - \bm{l}(\hat{z};\theta) \leq \sup_{\tau\in[0,t]} \left\{  F_{\theta}(t-\tau) + M\kappa^{-1}\tau \right\} = \mathcal{A}_{\theta}(t) \) wherever \(d(z',z) = \norm{x'-x}_{r} + \kappa\norm{y'-y}_{1} \leq t\). Since \(F_{\theta}(t)\) is concave,  \(\mathcal{A}_{\theta}(t)\) is also concave (see Lemma~\ref{lem:concave-function}), and \(\mathcal{A}_{\theta}(t)\) is an adversarial score of \(\bm{l}\).
      \hfill\Halmos
  \end{enumerate}

\subsection{Proof of Regression Proposition~\ref{prop:regression}}\label{proof:regression}
   
\noindent\textit{Proof.} We need to verify that \(\mathcal{A}_{\theta}\) is concave and for any \(t\geq0\), that is, 
\begin{equation*}
\sup\limits_{z'\in\mathcal{Z},\hat{z}\in\mathcal{Z}_N}       \left\{ \bm{l}(z';\theta) - \bm{l}(\hat{z};\theta) \mid d(z',\hat{z}) \leq  t \right\} \leq \mathcal{A}_{\theta}(t).
  \end{equation*}
Let \(u'=  {y' - f_{\theta}(x')} \) and \(u= {y - f_{\theta}(x)} \), then  \( \bm{l}(z';\theta) - \bm{l}(\hat{z};\theta) = \gamma\left( \abs{u'} \right) - \gamma\left(\abs{u}\right)\). Since \(\Gamma\) is an adversarial score of \(\gamma\), one has  \( \gamma\left( \abs{u'} \right) - \gamma\left(\abs{u}\right) \leq \Gamma(\abs{\abs{u'}-\abs{u}}) \leq   \Gamma (\abs{u'-u} ) \). 
  \begin{enumerate}
      \item[(a)] Since \(\kappa=\infty\), \(d(z',z)\leq t\) if and only if \(y'=y\) and \( \norm{x'-x}_{r}\leq t \).   In that case, we have \(\abs{u'-u} = \abs{f_{\theta}(x')-f_{\theta}(x)} \leq F_{\theta}(\norm{x'-x}_r)\leq  F_{\theta}(t) \) and thus \(\bm{l}(z';\theta) - \bm{l}(\hat{z};\theta)  \leq \Gamma (\abs{u'-u} ) \leq \Gamma (F_{\theta}(t))\). As  both \(\Gamma\) and \(F_{\theta}\) are non-decreasingly concave,  \(\mathcal{A}_{\theta} = \Gamma\circ F_{\theta}\) is also concave and therefore it is an adversarial score of \(\bm{l}\) at \(\mathcal{Z}_N\).
      
      \item[(b)] For any \(z,z'\) such that \(d(z',z) \leq t\), it is equivalent to \(\norm{x'-x}_{r} \leq t- \kappa\norm{y'-y}_{1} = t - \tau\) where \(\tau\coloneqq\kappa\norm{y'-y}_{1}\in[0,t]\). Hence,
      \begin{equation*}
      \begin{array}{ll}
         \abs{u'-u}\leq   \abs{f_{\theta}(x')-f_{\theta}(x)}+\abs{y'-y} 
              \leq  F_{\theta}(t-\tau) + \kappa^{-1}\tau. 
      \end{array}
      \end{equation*}
      Therefore, \(\bm{l}(z';\theta) - \bm{l}(\hat{z};\theta) \leq  \Gamma (\abs{u'-u} )  \leq \Gamma\left( \sup_{\tau\in[0,t]} \left\{  F_{\theta}(t-\tau) + \kappa^{-1}\tau \right\} \right) = \mathcal{A}_{\theta}(t) \) wherever \(d(z',z) = \norm{x'-x}_{r} + \kappa\norm{y'-y}_{1} \leq t\). By Lemma~\ref{lem:concave-function},  \(\mathcal{A}_{\theta}(t)\) is  concave  and thus
      \(\mathcal{A}_{\theta}(t)\) is an adversarial score of \(\bm{l}\) at \(\mathcal{Z}_N\).
      \hfill\Halmos
  \end{enumerate}

\section{Technical Lemmas and Proofs}

\subsection{Proof of Lemma~\ref{lem:p-infinity}} \label{proof:p-infinity} We shall prove that the point-wise RO is \(\mathcal{W}_{p=\infty}\)DRO. Let \(\rho>0\) be an arbitrary scalar. For any \(  \tilde{\mathbb{P}} \in \mathcal{P}(\mathcal{Z})  \) such that \(\mathcal{W}_{\infty}(\tilde{\mathbb{P}},\mathbb{P}_N) \leq \epsilon \), by the definition of $\mathcal{W}_{\infty}$, there exists \(\tilde{\pi}\in\Pi(\tilde{\mathbb{P}},\mathbb{P}_N ) \)  such that \(	\operatorname{ess.sup}_{\tilde{\pi}}(d) < \epsilon+\rho \). (Recall that  the essential supremum is defined as \(\operatorname{ess.sup}_{\tilde{\pi}}(d) \coloneqq \inf \left\{a \in\mathbb{R} \mid \tilde{\pi} \left(\left\{(\tilde{z},\hat{z}) \colon d(\tilde{z},z) > a  \right\}\right) = 0  \right\}\).) 
It means that \(\tilde{\pi} (\tilde{A}_{\epsilon+\rho}) = 1 \) where \(\tilde{A}_{\epsilon+\rho} \coloneqq \left\{(\tilde{z},z) \colon d(\tilde{z},z) <  \epsilon + \rho   \right\}\). Since the second marginal of \(\tilde{\pi}\) is \(\mathbb{P}_N\), one has
\begin{equation*}
    \tilde{\pi} (\tilde{A}_{\epsilon+\rho} \cap \mathcal{Z}\times\mathcal{Z}_N )=1.
\end{equation*}
Let \(B_{d,\epsilon+\rho}^{(i)} \coloneqq \{\tilde{z} \colon (\tilde{z}, Z^{(i)}) \in \tilde{A}_{\rho} \} = \{\tilde{z} \colon d(\tilde{z},Z^{(i)}) < \epsilon+\rho \} \) be the \((\epsilon+\rho)\)-ball centered at \(Z^{(i)}\). Then one has the following disjoint partition
\begin{equation*}
    \begin{array}{c}
         \tilde{A}_{\epsilon+\rho} \cap \mathcal{Z}\times\mathcal{Z}_N = \bigsqcup_{i=1}^{N}B_{d,\epsilon+\rho}^{(i)} \times \{Z^{(i)} \}. 
    \end{array}
\end{equation*}
As \(\tilde{\pi} (\tilde{A}_{\epsilon+\rho} \cap \mathcal{Z}\times\mathcal{Z}_N )=1\), it shows that 
\begin{equation*}
    \begin{array}{c}
         \mathbb{E}_{\tilde{\mathbb{P}}}[\bm{l}(\tilde{Z};\theta)] = \mathbb{E}_{\tilde{\pi}}[\bm{l}(\tilde{Z};\theta)] \leq \sum_{i=1}^{N}\mu_i \sup_{\tilde{z}\in B_{d,\epsilon+\rho}^{(i)}} \bm{l}(\tilde{z};\theta).
    \end{array}
\end{equation*}
In summary, for any \(  \tilde{\mathbb{P}} \in \mathcal{P}(\mathcal{Z})  \) such that \(\mathcal{W}_{\infty}(\tilde{\mathbb{P}},\mathbb{P}_N) \leq \epsilon \), we have shown that \(\mathbb{E}_{\tilde{\mathbb{P}}}[\bm{l}(\tilde{Z};\theta)] \leq \sum_{i=1}^{N}\mu_i \sup_{\tilde{z}\in B_{d,\epsilon+\rho}^{(i)}} \bm{l}(\tilde{z};\theta) \). Therefore, \(\mathcal{R}_p(\epsilon)\leq \sum_{i=1}^{N}\mu_i \sup_{\tilde{z}\in B_{d,\epsilon+\rho}^{(i)}} \bm{l}(\tilde{z};\theta)\).

On the other hand, for any collection of point-wise attack \( \{\tilde{Z}^{(i)}\}_{i=1}^{N} \) satisfying \(\tilde{Z}^{(i)} \in B_{d,\epsilon}^{(i)} \), define \(\tilde{\mathbb{P}} \coloneqq \sum_{i=1}^N \mu_i \bm{\chi}_{\tilde{Z}^{(i)}}\). Then \(\mathcal{W}_{\infty}(\tilde{\mathbb{P}},\mathbb{P}_N) \leq \epsilon \) and therefore,
\begin{equation*}
     \begin{array}{c}
          \sum_{i=1}^{N}\mu_i \sup_{\tilde{z}\in B_{d,\epsilon}^{(i)}} \bm{l}(\tilde{z};\theta) \leq \sup_{\mathbb{P} \colon \mathcal{W}_{\infty}(\mathbb{P},\mathbb{P}_N) \leq \epsilon }   \mathbb{E}_{\mathbb{P}} [\bm{l} (Z;\theta) ] = \mathcal{R}_p(\epsilon).
     \end{array}
\end{equation*}

\subsection{Properties of concave function.}

\begin{lemma}
    [Three-slope lemma]  \label{lem:3chords} \citep{roberts1974convex}
Let \(\Gamma: I \to \mathbb{R}\) be a univariate function defined on an interval \(I \subseteq \mathbb{R}\). 
The function \(\Gamma\) is concave if and only if for any \(t_1, t_2, t_3 \in I\) such that \(t_1 < t_2 < t_3\), \(\frac{\Gamma(t_2) - \Gamma(t_1)}{t_2 - t_1} \geq \frac{\Gamma(t_3) - \Gamma(t_1)}{t_3 - t_1} \geq \frac{\Gamma(t_3) - \Gamma(t_2)}{t_3 - t_2}.\)
\end{lemma}

\begin{lemma}\label{lem:concave-function}
    Suppose that \(\varphi, \varphi_2\colon[0,\infty)\rightarrow[0,\infty)  \) are non-decreasingly concave. 
    \begin{itemize}
        \item \(\varphi\circ\varphi_2 \) is also non-decreasingly concave.
        \item \(\phi(t)\coloneqq \sup_{\tau\in[0,t]}\left\{ \varphi(t-\tau) + c\tau \right\} \) is non-decreasingly concave for any \(c>0\).
    \end{itemize}  
\end{lemma}

\noindent\textit{Proof.} The first item follows \citet[Theorem 5.1]{rockafellar-convex-analysis}. To prove the second item, fix arbitrary \(0< t_1\leq t_2<\infty\) and \(\eta\in[0,1]\). One has  \(\phi(t_1)=\sup_{\tau\in[0,t_1]}\left\{ \varphi(t_1-\tau) + c\tau \right\} \leq \sup_{\tau\in[0,t_2]}\left\{ \varphi(t_1-\tau) + c\tau \right\}  \leq \sup_{\tau\in[0,t_2]}\left\{ \varphi(t_2-\tau) + c\tau \right\}=\phi(t_2) \) since \(\varphi\) is non-decreasing, thus \(\phi\) is non-decreasing. In addition,
\begin{equation*}
    \begin{array}{rl}
    &\eta\phi(t_1) +(1-\eta)\phi(t_2)\\
       = &\eta\sup_{\tau_1\in[0,t_1]}\left\{\varphi(t_1-\tau_1) + c\tau_1\right\} 
        + (1-\eta)\sup_{\tau_2\in[0,t_2]}\left\{\varphi(t_2-\tau_2) + c\tau_2\right\} \\
        =&         \sup\limits_{\tau_1\in[0,t_1],\tau_2\in[0,t_2] } \left\{   \eta\left(\varphi(t_1-\tau_1) + c\tau_1\right) 
            + (1-\eta)\left(\varphi(t_2-\tau_2) + c\tau_2\right)\right\} 
        % \begin{array}{l}
        %     \eta\left(\varphi(t_1-\tau_1) + c\tau_1\right) \\
        %     + (1-\eta)\left(\varphi(t_2-\tau_2) + c\tau_2\right)
        % \end{array}
        \\
        \leq&         \sup\limits_{\tau_1\in[0,t_1],\tau_2\in[0,t_2] }   \left\{ 
        \varphi\left( \eta (t_1-\tau_1) + (1-\eta)(t_2-\tau_2)\right)
        + c\left(\eta\tau_1+(1-\eta)\tau_2\right)
        % \begin{array}{l}
        %     \varphi\left( \eta (t_1-\tau_1) + (1-\eta)(t_2-\tau_2)\right)\\
        % + c\left(\eta\tau_1+(1-\eta)\tau_2\right)
        % \end{array}
        \right\} \\
        \leq& \sup\limits_{\tau\in[0,\eta t_1+(1-\eta)t_2]}\left\{  \varphi(\eta t_1 + (1-\eta)t_2 - \tau) + c\tau\right\}
        = \phi(\eta t_1+ (1-\eta)t_2).
    \end{array}
\end{equation*} Here the first inequality follows from the concavity of  \(\varphi\), and the second inequality follows from letting \(\tau = \eta\tau_1 +  (1-\eta)\tau_2\). Thus \(\phi\) is concave.

\subsection{Univariate least concave majorant}
\begin{lemma}
    \label{lem:univariate-majorant}
Suppose that \(\gamma\colon[0,\infty)\rightarrow[0,\infty) \) is non-decreasing. Let \(\Gamma\) be the least concave majorant (Definition~\ref{def:least-majorant}) of \(\gamma\). Furthermore, define the rate function 
\begin{equation*}
    \begin{array}{c}
         \Delta_{\gamma}(t) = \sup_{s\geq 0}\{ \gamma(s+t)-\gamma(s)\},
    \end{array}
\end{equation*}
and let \(\mathcal{C}_{\gamma}\) and \(\mathcal{C}_{\Delta_\gamma}\) be the least concave majorant of \(\gamma\) and \(\Delta_{\gamma}\). Then the following properties hold.
\begin{itemize}
    \item[$\operatorname{(a)}$] The least concave majorant \(\mathcal{C}_{\gamma}\) of \(\gamma\) is also non-decreasing on \([0,\infty)\).
    \item[$\operatorname{(b)}$] If \(\gamma\) is \(L\)-Lipschitz, then \(\gamma(t) - \gamma(0) \leq \Delta_{\gamma}(t)\leq \mathcal{C}_{\Delta_\gamma}(t) \leq Lt \) for any \(t\geq0\).
    \item[$\operatorname{(c)}$] If \(\gamma\) is concave, then \(\Delta_{\gamma}(t)=\mathcal{C}_{\Delta_\gamma}(t) = \gamma(t) - \gamma(0)\) for any \(t\geq0\).
    \item[$\operatorname{(d)}$] If \(\Delta_{\gamma}\) is concave, then obviously \(\Delta_{\gamma}(t)=\mathcal{C}_{\Delta_\gamma}(t)\) for any \(t\geq0\).  
\end{itemize}
\end{lemma}

\noindent\textit{Proof.} \begin{itemize}
    \item[$\operatorname{(a)}$] Suppose that \(\mathcal{C}_{\gamma}\) is \textit{not} non-decreasing on \([0,\infty)\). That is to say, there exists \(0\leq t_1 < t_2 <\infty \) such that \(\mathcal{C}_{\gamma}(t_1) > \mathcal{C}_{\gamma}(t_2) \). Since \(\mathcal{C}_{\gamma}\) is concave, by Lemma~\ref{lem:3chords}, for any \(t_3 >t_2\) one has \(\frac{\mathcal{C}_{\gamma}(t_3)-\mathcal{C}_{\gamma}(t_1)}{t_3-t_1} \leq \frac{\mathcal{C}_{\gamma}(t_2) - \mathcal{C}_{\gamma}(t_1)}{t_2-t_1}<0. \) Let \(a\coloneqq \frac{\mathcal{C}_{\gamma}(t_2) - \mathcal{C}_{\gamma}(t_1)}{t_2-t_1}\) and \(b\coloneqq - at_1 + \mathcal{C}_{\gamma}(t_1)\). Then \(a<0\) and \(\mathcal{C}_{\gamma}(t_3) \leq a t_3 + b\) for any \(t_3>t_2\). This implies that \(\lim_{t_3\rightarrow+\infty}\mathcal{C}_{\gamma}(t_3)\leq -\infty\), which contradicts the fact that \(\mathcal{C}_{\gamma}(t) \geq \gamma(t)\geq0 \). 
\item[$\operatorname{(b)}$] 
The first inequality holds by choosing \(s=0\). For any \(s, t \geq 0\), we have \(\gamma(s+t) - \gamma(s) \leq L|(s+t) - s| = Lt\). Taking the supremum over \(s \geq 0\), we obtain the second inequality. By definition of least concave majorant, \(\Delta_{\gamma}(t) \leq \mathcal{C}_{\Delta_\gamma}(t)\). Finally, since \(h(t) = Lt\) is concave and \(\Delta_{\gamma} \leq h\), the lowest concave upper bound \(\mathcal{C}_{\Delta_\gamma}\) must satisfy \(\mathcal{C}_{\Delta_\gamma}(t) \leq Lt\).   
\item[$\operatorname{(c)}$] If \(\gamma\) is concave, the function \(s \mapsto \gamma(s+t) - \gamma(s)\) is non-increasing in \(s\) for any fixed \(t \geq 0\). Therefore, the supremum defining \(\Delta_{\gamma}\) is achieved at \(s=0\), yielding \(\Delta_{\gamma}(t) = \gamma(t) - \gamma(0)\), thus \(\Delta_{\gamma}\) is also concave and \(\mathcal{C}_{\Delta_\gamma}(t) = \Delta_{\gamma}(t) = \gamma(t) - \gamma(0)\).
% \item[$\operatorname{(d)}$] 
% If \(\Delta_{\gamma}\) is already concave, then it is the least concave majorant of itself.
\hfill\Halmos
\end{itemize}

\end{appendices}

\bibliography{references} % if more than one, comma separated

%%%%%%%%%%%%%%%%%
\end{document}